\documentclass[11pt]{amsart}
\usepackage{amsmath}
\usepackage{amssymb}
\usepackage{extarrows}
\usepackage{mathabx}
\usepackage[all]{xy}
\usepackage{amsbsy}
\usepackage{graphicx}
\usepackage{appendix}
\usepackage{float}
\usepackage{subfigure}
\usepackage[numbers,sort&compress]{natbib}
\usepackage{graphicx}
\usepackage{amsthm}
\usepackage{geometry}
\usepackage{fancyhdr}
\usepackage{color} 
\usepackage{overpic}
\usepackage{mathabx}
\usepackage{tikz-cd}
\usepackage{enumerate}
\usepackage{mathrsfs}
\usepackage{colonequals}
\usepackage{microtype}
\usepackage{hyperref}

\newtheorem{thm}{Theorem}[section]
\newtheorem{cor}[thm]{Corollary}
\newtheorem{prop}[thm]{Proposition}
\newtheorem{lem}[thm]{Lemma}

\theoremstyle{definition}
\newtheorem{defn}[thm]{Definition}
\newtheorem{cons}[thm]{Construction}
\newtheorem{exmp}[thm]{Example}
\newtheorem{conj}[thm]{Conjecture}
\newtheorem*{fact}{Fact}

\newtheorem*{conv}{Convention}

\newtheorem*{org}{Organization}
\newtheorem*{ack}{Acknowledgement}

\theoremstyle{remark}
\newtheorem{rem}[thm]{Remark}

\numberwithin{equation}{section}

\newcommand{\beq}{\begin{equation*}\begin{aligned}}
\newcommand{\eeq}{\end{aligned}\end{equation*}}
\newcommand{\bpf}{\begin{proof}}
\newcommand{\epf}{\end{proof}}
\newcommand{\bthm}{\begin{thm}}
\newcommand{\ethm}{\end{thm}}
\newcommand{\bprop}{\begin{prop}}
\newcommand{\eprop}{\end{prop}}
\newcommand{\bcor}{\begin{cor}}
\newcommand{\ecor}{\end{cor}}
\newcommand{\blem}{\begin{lem}}
\newcommand{\elem}{\end{lem}}
\newcommand{\bdefn}{\begin{defn}}
\newcommand{\edefn}{\end{defn}}
\newcommand{\bcons}{\begin{cons}}
\newcommand{\econs}{\end{cons}}
\newcommand{\bexmp}{\begin{exmp}}
\newcommand{\eexmp}{\end{exmp}}
\newcommand{\brem}{\begin{rem}}
\newcommand{\erem}{\end{rem}}
\newcommand{\bfa}{\begin{fact}}
\newcommand{\efa}{\end{fact}}
\newcommand{\benu}{\begin{enumerate}[(1)]}
\newcommand{\eenu}{\end{enumerate}}

\newcommand{\bdia}{\begin{displaymath}\xymatrix}
\newcommand{\edia}{\end{displaymath}}

\newcommand{\shi}{\underline{\rm SHI}}

\newcommand{\shib}{{\mathbf{SHI}}}

\newcommand{\shg}{\mathbf{SH}}

\newcommand{\khii}{\underline{\rm KHI}}

\newcommand{\shiib}{\underline{\mathbf{SHI}}}

\newcommand{\deq}{\colonequals}

\newcommand{\spin}{{\rm Spin}^c}

\newcommand{\al}{\alpha}
\newcommand{\be}{\beta}
\newcommand{\ga}{\gamma}
\newcommand{\Ga}{\Gamma}

\newcommand{\p}{\prime}
\newcommand{\pp}{{\prime\prime}}

\newcommand{\gr}{\chi_{\rm gr}}

\newcommand{\en}{\chi_{\rm en}}
\newcommand{\aand}{~{\rm and}~}

\newcommand{\intg}{\mathbb{Z}}
\newcommand{\ft}{{\mathbb{F}_2}}

\newcommand{\posi}{\mathbb{N}_+}

\newcommand{\ra}{\rightarrow}

\newcommand{\xra}{\xrightarrow}

\DeclareMathOperator{\cok}{Coker}
\DeclareMathOperator{\im}{Im}
\DeclareMathOperator{\ke}{Ker}

\DeclareMathOperator{\cone}{Cone}

\begin{document}

\title{SU(2) representations and a large surgery formula}


\author{Zhenkun Li}
\address{Academy of Mathematics and Systems science\\Chinese Academy of Science}
\email{zhenkun@amss.ac.cn}

\author{Fan Ye}
\address{Department of Mathematics, Harvard University}
\email{fanye@math.harvard.edu}

\begin{abstract}
A knot $K\subset S^3$ is called $SU(2)$-abundant if it satisfies two conditions: first, for all but finitely many $r\in\mathbb{Q}\backslash\{0\}$, there exists an irreducible representation $\pi_1(S^3_r(K))\to SU(2)$; second, any slope $r=u/v\neq 0$ for which $S^3_r(K)$ admits no irreducible $SU(2)$ representation must satisfy $\Delta_K(\zeta^2)= 0$ for some $u$-th root of unity $\zeta$. We show that if a nontrivial knot $K\subset S^3$ is not $SU(2)$-abundant then it is a prime knot whose Alexander polynomial $\Delta_K(t)$ has coefficients restricted to $\{-1,0,1\}$. This implies, in particular, that all hyperbolic alternating knots are $SU(2)$-abundant. Our proof hinges on a large surgery formula connecting instanton knot homology $KHI(S^3,K)$ and framed instanton homology $I^\sharp(S^3_n(K))$ for integers $n$ satisfying $|n|\ge 2g(K)+1$. Using this technique, we derive several interesting results in instanton Floer homology: for any Berge knot $K$, the spaces $KHI(S^3,K)$ and $\widehat{HFK}(S^3,K)$ have identical dimension; for any dual knot $K_r\subset S^3_r(K)$ of a Berge knot $K$ with $r> 2g(K)-1$, we prove $\dim_\mathbb{C}KHI(S^3_r(K),K_r)=|H_1(S^3_r(K);\mathbb{Z})|$; and for any genus-one alternating knot $K$ and any $r\in\mathbb{Q}\backslash\{0\}$, the spaces $I^\sharp(S^3_r(K))$ and $\widehat{HF}(S_r^3(K))$ have equal dimension.
\end{abstract}
\maketitle
\tableofcontents

\section{Introduction}
The fundamental group stands as the most significant invariant of a 3-manifold. However, directly analyzing the fundamental group often proves challenging. A productive alternative approach to studying the fundamental group involves studying homomorphisms from the fundamental group to more tractable groups (e.g. $SU(2),SL(2,\mathbb{C}),SL(2,\mathbb{R})$), which yields computationally accessible invariants. Notable examples include the Casson invariant \cite{casson90} and the Casson-Lin invariant \cite{lin92}, both constructed using $SU(2)$ representations, as well as the A-polynomial \cite{apolynomial}, which is derived from the $SL(2,\mathbb{C})$ character variety.

In this paper, we study $SU(2)$ representations of a 3-manifold $Y$, i.e.\, homomorphisms from the fundamental group $\pi_1(Y)$ to $SU(2)$. For a knot $K$ in $S^3$, let $\Delta_K(t)\in\mathbb{Z}[t,t^{-1}]$ denote its symmetrized Alexander polynomial with conditions \begin{equation}\label{eq: alex conditions}
    \Delta_K(t)=\Delta_K(t^{-1})\aand \Delta_K(1)=1.
\end{equation}
For a knot $K$ in a closed 3-manifold $Y$, we write $Y(K)=Y\backslash{\rm int}N(K)$ for the knot complement and $Y_r(K)$ for the manifold obtained from $Y$ by a Dehn surgery along $K$ with slope $r$ with respect to a chosen framing of the knot $K$. If $K\subset Y$ is (integrally) null-homologous, then we always use the canonical Seifert framing.
\bdefn
An $SU(2)$ representation $\rho:\pi_1(Y)\to SU(2)$ is called \textbf{abelian} if the image ${\rm im}(\rho)$ is contained in an abelian subgroup of $SU(2)$. An $SU(2)$ representation is called \textbf{irreducible} if it is not abelian. A knot $K\subset S^3$ is called \textbf{$SU(2)$-abundant} if the following two conditions hold:
\benu
\item For all but finitely many $r\in\mathbb{Q}\backslash\{0\}$, the manifold $S^3_r(K)$ has an irreducible $SU(2)$ representation.
\item For any $r=u/v\neq 0$ such that $S^3_r(K)$ has only abelian $SU(2)$ representations, there is some $u$-th root of unity $\zeta$ such that $\Delta_K(\zeta^2)= 0$.
\eenu
\edefn
\brem\label{rem: nondegenerate}
The first condition implies that $K$ is not \textbf{$SU(2)$-averse} in the sense of \cite{sivek20su2}. Note that if $b_1(Y) = 0$, then an $SU(2)$ representation of $Y$ has abelian image if and only if it has cyclic image. The second condition corresponds to some nondegenerate condition as in \cite[Corollary 4.8]{baldwin2018stein}. By \cite[Remark 1.6]{baldwin2019lspace}, when $u$ is a prime power, $\Delta_K(\zeta^2)\neq 0$ for any $K$ and any $u$-th root of unity $\zeta$. Moreover, rationals with prime power numerators are dense in $\mathbb{Q}$. 
\erem
Suppose $K\subset S^3$ is a nontrivial knot and $r\in\mathbb{Q}$. It is already known that if $|r|\le 2$ \cite[Theorem 1]{kronheimer04su2} or if $|r|$ is sufficiently large \cite[Corollary 1.2]{sivek20su2}, then $S_r^3(K)$ has an irreducible $SU(2)$ representation. There are many other closed 3-manifolds with irreducible $SU(2)$ representations; see \cite{kronheimer04witten,jianfeng16su2,Zentner17su2,Zentner18sl2c,baldwin2018stein,toroidal21su2,baldwin21splicing,Zentner21menagerie,yixiesu2}.

In this paper, we provide some sufficient conditions for $SU(2)$-abundant knots.

\bthm\label{thm: main SU2}
If a nontrivial knot $K\subset S^3$ is not $SU(2)$-abundant, then it satisfies the following conditions.

\benu
\item There exist $k\in\mathbb{N}_+$ and integers $n_k>n_{k-1}>\dots>n_1>n_0=0$ such that $$\pm \Delta_K(t)=(-1)^k+\sum_{j=1}^k(-1)^{k-j}(t^{n_j}+t^{-n_j}).$$
\item The Seifert genus satisfies $g(K)=n_k=n_{k-1}+1$.
\item $K$ is a prime knot, i.e.\, it is not a connected sum of two knots.

\eenu
\ethm
\brem\label{rem: determinant genus}
By term (1) and term (2) in Theorem \ref{thm: main SU2}, we have \begin{equation}\label{eq: det genus}\det(K)=|\Delta_K(-1)|\le 2k+1\le 2g(K)+1.\end{equation}
\erem
\brem\label{rem: more obstructions}
In \cite[Theorem 1.5]{baldwin2019lspace} and \cite[Corollary 1.7, and Proposition 5.4]{baldwin2020concordance}, Baldwin and Sivek proved that a nontrivial knot $K$ is $SU(2)$-abundant unless $K$ is both fibered and strongly quasi-positive (up to mirror), the 4-ball genus $g_4(K)$ equals $g(K)$, and the slope $r$ with no irreducible $SU(2)$ representations satisfies $|r|\ge 2g(K)-1$. It is worth mentioning that by techniques developed in this paper, it is possible to provide alternative proofs of those results. 
\erem
From classification results in \cite{Ozsvath2005,baker18montesinos,vafaee21braid}, we have the following corollary.
\bcor\label{cor: abundant knot su2}
The following knots are $SU(2)$-abundant.
\benu
\item Hyperbolic alternating knots, i.e.\, alternating knots that are not torus knots $T(2,2n+1)$.
\item Montesinos knots (including all pretzel knots), except torus knots $T(2,2n+1)$, pretzel knots $P(-2, 3, 2n + 1)$ for $n\in\mathbb{N}_+$, and their mirrors.
\item Knots that are closures of 3-braids, except twisted torus knots $K(3,q;2,p)$ with $pq>0$ and their mirrors, where $K(3,q;2,p)$ is the closure of a 3-braid made up of a $(3, q)$ torus braid with $p$ full twist(s) on two adjacent strands.
\eenu
\ecor

The proof of Theorem \ref{thm: main SU2} is based on instanton knot homology $KHI(Y,K)$ \cite{kronheimer2010knots} and framed instanton homology $I^\sharp(Y)$ \cite{kronheimer2011khovanov}, which are vector spaces over $\mathbb{C}$ for a knot $K$ in a closed 3-manifold $Y$. There are relative $\mathbb{Z}_2$-gradings on $KHI(Y,K)$ and $I^\sharp(Y)$. Furthermore, a Seifert surface $S$ of $K$ induces a $\mathbb{Z}$-grading on $KHI(Y,K)$ \cite{kronheimer2010instanton,li2019direct,li2019decomposition}, which we write as$$KHI(Y,K)=\bigoplus_{i\in\mathbb{Z}}KHI(Y,K,S,i).$$

We write $(-Y,K)$ for the induced knot in the manifold $-Y$ obtained from $Y$ by reversing the orientation and call it the \textbf{mirror} of $K$ or $(Y,K)$. For a knot $K$ in $S^3$, we write $\widebar{K}$ for the mirror of $K$, i.e.\, $(S^3,\widebar{K})=(-S^3,K)$. We write $-K$ for the knot with reverse orientation, which is different from $\widebar{K}$. Then we have canonical isomorphisms
\begin{equation}\label{eq: mirror iso}
    KHI(-Y,K,S,i)\cong {\rm Hom}_\mathbb{C}(KHI(Y,K,S,-i),\mathbb{C}) \aand I^\sharp(-Y)\cong {\rm Hom}_\mathbb{C}(I^\sharp(Y),\mathbb{C}).
\end{equation}
\bdefn
A rational homology sphere $Y$ is called an \textbf{instanton L-space} if $\dim_\mathbb{C}I^\sharp(Y)=|H_1(Y;\mathbb{Z})|$. A knot $K$ in an instanton L-space $Y$ is called an \textbf{instanton L-space knot} if a nontrivial surgery on it also gives an instanton L-space. We call $K$ a \textbf{positive instanton L-space knot} if a positive surgery on it also gives an instanton L-space.
\edefn
\brem\label{rem: su2 abundant}
It follows directly from (\ref{eq: mirror iso}) that $Y$ is an instanton L-space if and only if $-Y$ is an instanton L-space.  Since $S_r^3(\widebar{K})=-S_{-r}^3(K)$, a positive surgery on $K$ gives an instanton L-space if and only if a negative surgery on $\widebar{K}$ gives an instanton L-space. By \cite[Theorem 1.1]{sivek20su2} and \cite[Corollary 4.8]{baldwin2018stein}, if $K\subset S^3$ is not $SU(2)$-abundant, then $K$ is an instanton L-space knot. By \cite[Theorem 1.15]{baldwin2019lspace} and passing to the mirror if necessary, we can further assume that for any sufficiently large integer $n$, the manifold $S_n^3(K)$ is an instanton L-space.
\erem
 
The following theorem is the main theorem of this paper.

\bthm\label{thm: main 2}
If $K\subset S^3$ is an instanton L-space knot, then $K$ is a prime knot and there exists $k\in\mathbb{N}$ and integers $$n_k>n_{k-1}>\dots>n_1>n_0=0> n_{-1}>\dots>n_{1-k}>n_{-k}\text{ with }n_{-j}=-n_j$$ such that $$\dim_\mathbb{C}KHI(S^3,K,S,i)=\begin{cases}1&\text{if }i=n_j\text{ for }j\in[-k,k],\\0 &\text{else}, \end{cases},$$where the $\mathbb{Z}_2$-gradings of the generators of $KHI(S^3,K,S,n_j)\cong\mathbb{C}$ are alternating with respect to $j$.
\ethm
We prove Theorem \ref{thm: main SU2} by Theorem \ref{thm: main 2}.
\bpf[Proof of Theorem \ref{thm: main SU2}]
By Remark \ref{rem: su2 abundant}, if $K\subset S^3$ is not $SU(2)$-abundant, then $K$ is an instanton L-space knot. Then Theorem \ref{thm: main 2} applies to $K$, and we obtain term (3). Since the space in the top $\mathbb{Z}$-grading of $KHI(S^3,K)$ is one-dimensional, it follows from \cite[Section 7]{kronheimer2010knots} that $K$ is fibered. Then by \cite[Theorem 1.7]{baldwin2018khovanov}, we know that $\dim_\mathbb{C}KHI(S^3,K,S,g(K)-1)\ge 1$, and Theorem \ref{thm: main 2} forces equality to hold. Thus, terms (1) and (2) follow from$$\sum_{i\in\mathbb{Z}}\chi(KHI(S^3,K,S,i))\cdot t^i=\pm \Delta_K(t)$$\cite{Lim2009,kronheimer2010instanton}, where the sign ambiguity is due to the relative $\mathbb{Z}_2$-grading.
\epf

Theorem \ref{thm: main 2} is an instanton analog of \cite[Theorem 1.2]{Ozsvath2005} in Heegaard Floer theory due to Ozsv\'{a}th and Szab\'{o}. The key step to prove Theorem \ref{thm: main 2} is to establish an instanton version of the large surgery formula in Heegaard Floer theory. We will explain more details about this strategy in Subsection \ref{subsec: Constructions parallel to Heegaard Floer theory}. Here we state more applications of techniques developed in this paper.

First, we can compare instanton knot homology of an instanton L-space knot $K\subset Y$ to the knot Floer homology $\widehat{HFK}(Y,K)$ introduced in \cite{ozsvath2004holomorphicknot,Rasmussen2003}, which verifies more examples of \cite[Conjecture 7.24]{kronheimer2010knots}. The main inputs are a generalization of Theorem \ref{thm: main 2}, results about Heegaard Floer theory from \cite{Ozsvath2005,Rasmussen2017}, and the equation of graded Euler characteristics from \cite{LY2021}\begin{equation}\label{eq: chi same}
    \gr(KHI(Y,K))=\gr(\widehat{HFK}(Y,K))\in\mathbb{Z}[H]/\pm H,
\end{equation}
where $H=H_1(Y(K);\mathbb{Z})/{\rm Tors}$.
\bdefn[{\cite{ozsvath2004holomorphic,Ozsvath2005}}]
A rational homology sphere $Y$ is called an \textbf{(Heegaard Floer) L-space} if $\dim_{\mathbb{F}_2}\widehat{HF}(Y)=|H_1(Y;\mathbb{Z})|$. A knot $K$ in an L-space $Y$ is called an \textbf{(Heegaard Floer) L-space knot} if a nontrivial surgery on it also gives an L-space. 
\edefn

\bthm\label{main: L space knot}
Suppose $K\subset Y$ is a knot with $H_1(Y(K);\mathbb{Z})\cong \mathbb{Z}$ and suppose the meridian of $K$ represents $q$ times the generator of $H_1(Y(K);\mathbb{Z})$. Suppose $K$ is both an L-space knot and an instanton L-space knot such that $Y_{u/v}(K)$ is an instanton L-space. If $\gcd(q,v)=1$, then we have
\begin{equation}\label{eq: hfk=khi}
    \dim_\mathbb{C}KHI(Y,K)=\dim_{\mathbb{F}_2}\widehat{HFK}(Y,K).
\end{equation}Moreover, when properly fixing the gradings associated to the Seifert surface $S$ of $K$, we have \begin{equation}\label{eq: dim equal s}
    \dim_\mathbb{C}KHI(Y,K,S,i)=\dim_{\mathbb{F}_2}\widehat{HFK}(Y,K,S,i)\le 1 \text{ for any }i\in\mathbb{Z}.
\end{equation}
\ethm
\brem
When $H_1(Y(K);\mathbb{Z})$ has torsion, we can still decompose $KHI(Y,K)$ along elements in $H_1(Y(K);\mathbb{Z})$ as in \cite{LY2021enhanced}. However, since this decomposition is not canonical and adapting the proofs to this case is subtle, we leave the discussion in this case to the future. That is why we assume $H_1(Y(K);\mathbb{Z})\cong \mathbb{Z}$ and $\gcd(q,v)=1$. 
\erem
\brem
From \cite{lidman2020framed,baldwin2020concordance}, for a knot $K\subset S^3$ that is both an L-space knot and an instanton L-space knot, we have $\dim_\mathbb{C}I^\sharp(S_r^3(K))=\dim_{\mathbb{F}_2}\widehat{HF}(S_r^3(K))$ for any $r\in\mathbb{Q}$.
\erem
From \cite[Corollary 1.3]{alfieri2020framed}, a Seifert fibered space is an L-space if and only if it is an instanton L-space. In particular, closed 3-manifolds with elliptic geometry are both (Heegaard Floer) L-spaces and instanton L-spaces \cite[Proposition 2.3]{Ozsvath2005} (or equivalently, with finite fundamental group by the Geometrization theorem; see \cite{elliptic}). In particular, $S^3$, the Poincar\'{e} sphere $\Sigma(2,3,5)$, and all lens spaces $L(p,q)$ are both (Heegaard Floer) L-spaces and instanton L-spaces. From \cite{ozsvath2005double,scaduto2015instanton}, double-branched covers of Khovanov-thin knots (in particular, all quasi-alternating knots) are also both L-spaces and instanton L-spaces. Note that when $Y$ is an integral homology sphere in Theorem \ref{main: L space knot}, then we have $q=1$ and hence $\gcd(q,v)=1$ for any $v$. Thus, we have the following corollary.

\bcor
Suppose $K$ is a knot in $Y=S^3$ or the Poincar\'{e} sphere $\Sigma(2,3,5)$. If there is some $r\in\mathbb{Q}\backslash\{0\}$ such that $Y_r(K)$ is a Seifert fibered L-space or a double-branched cover of a Khovanov-thin knot. Then (\ref{eq: hfk=khi}) and (\ref{eq: dim equal s}) hold.
\ecor
\brem
There are many examples of knots in $S^3$ and $\Sigma(2,3,5)$ that admit lens space surgeries, such as Berge's knots \cite{Berge2018b} in $S^3$, Tange's knots \cite[Theorem 4.1]{Tange09} in $\Sigma(2,3,5)$, Hedden's knots \cite{Hedden2007} in $\Sigma(2,3,5)$ dual to $T_R$ and $T_L$ in lens spaces (see also \cite{Rasmussen2007,baker14knot,baker20lspace}), Baker's tunnel number two knots \cite{baker20lspace} in $\Sigma(2,3,5)$. There are also other twist families of knots admitting Seifert fibered L-space surgeries \cite{motegi16twisted,baker19twisted}.
\erem

Second, we can relate the knot in the following definition to the framed instanton homology of large surgeries on it. The main input is the large surgery formula introduced in Subsection \ref{subsec: Constructions parallel to Heegaard Floer theory}. The analog in Heegaard Floer theory was proved in \cite[Section 3]{Rasmussen2017}.

\bdefn
A knot $K$ in an instanton L-space $Y$ is called an \textbf{instanton Floer simple knot} if $\dim_\mathbb{C}KHI(Y,K)=\dim_\mathbb{C}I^\sharp(Y)=|H_1(Y;\mathbb{Z})|.$
\edefn

\bthm[]\label{thm: floer simple khi}
Suppose $K\subset Y$ is a knot with $H_1(Y(K);\mathbb{Z})\cong\mathbb{Z}$. Suppose the basis of $H_1(\partial Y(K);\mathbb{Z})$ is induced by the meridian of $K$. Then $K$ is an instanton Floer simple knot if and only if, for any $r\in\mathbb{Q}$ with $|r|$ sufficiently large, the manifold $Y_r(K)$ is an instanton L-space.
\ethm
\brem
In \cite[Theorem 1.10]{LY2020}, we proved Theorem \ref{thm: floer simple khi} for simple knots in lens spaces without assuming $H_1(Y(K);\mathbb{Z})\cong\mathbb{Z}$. The technique there is different from the ones in this paper.
\erem
\brem
From \cite[Theorem 1.15]{baldwin2019lspace}, we know that if $K\subset S^3$ is a positive instanton L-space knot, then $S_r^3(K)$ is an instanton L-space if and only if $r\ge 2g(K)-1$. Hence, we can apply Theorem \ref{thm: floer simple khi} to the dual knot $K_r\subset S_r^3(K)$ of a Berge knot $K$ with $r> 2g(K)-1$ to obtain that $K_r$ is an instanton Floer simple knot.
\erem

Third, we can make some calculations for manifolds obtained from surgeries on genus-one knots. If $K\subset S^3$ with $g(K)=1$, we can use the large surgery formula introduced in Subsection \ref{subsec: Constructions parallel to Heegaard Floer theory} to compute $I^\sharp(S_r^3(K))$ when $|r|$ sufficiently large (indeed $|r|\ge 2g(K)+1=3$ is large enough). Furthermore, we can compute $I^\sharp(S_r^3(K))$ for any slope $r$ by the concordance invariant $\nu^{\sharp}(K)$ defined by Baldwin and Sivek \cite{baldwin2020concordance}. In particular, we have the following theorem.

\bthm\label{thm: alternating}
Suppose $K$ is a genus-one alternating knot. Then for any $r\in\mathbb{Q}\backslash\{0\}$, we have
$$\dim_{\mathbb{C}}I^{\sharp}(S^3_r(K))=\dim_{\ft}\widehat{HF}(S^3_r(K)).$$
\ethm
\brem
For genus-one Khovanov-thin knots (in particular, genus-one quasi-alternating knots \cite[Corollary 1.6]{kronheimer2011khovanov}), we can also fix the value $\dim_\mathbb{C}I^{\sharp}(S^3_r(K))$ up to the mirror of $K$; see Section \ref{sec: Dehn surgery along genus-one knots} for more details.
\erem

After the first announcement of this paper, several results have been developed using the results and ideas from this paper:
\begin{itemize}
	\item In \cite{BSLY2024}, the large surgery formula developed in the current paper was utilized to established the main result that the fundamental group of the $3$-surgery of any non-trivial knots in $S^3$ admits an irreducible $SU(2)$-representation.
	\item In \cite{LY2022surgeryformula1,LY2022surgeryformula2}, a general surgery formula was developed based on the large surgery formula in the current paper.
	\item In \cite{LY2024torsion}, a variation of the large surgery formula in the current paper was used to derive that the unreduced singular instanton Floer homology of any unknotting-number-one knot in $S^3$ admits $2$-torsion.
\end{itemize}

\subsection{A large surgery formula in instanton theory}\label{subsec: Constructions parallel to Heegaard Floer theory}\quad

In this subsection, we sketch the idea of the proof of Theorem \ref{thm: main 2} and introduce a large surgery formula relating $KHI(S^3,K)$ and $I^\sharp(S^3_n(K))$ for any integer $n$ that satisfies $|n|\ge 2g(K)+1$. By Remark \ref{rem: su2 abundant}, we may assume $S_n^3(K)$ is an instanton L-space for any sufficiently large integer $n$.  However, to apply the proof of \cite[Theorem 1.2]{Ozsvath2005}, we need to recover (at least partially) the following structures in instanton theory.
\bfa\label{fact: HF construction}
Suppose $K$ is a knot in $S^3$ and $n\in\mathbb{N}_+$. We have the following structures in Heegaard Floer theory \cite{ozsvath2004holomorphic,ozsvath2004holomorphicknot,Rasmussen2003}:
\benu
\item A decomposition of $\widehat{HF}(S_n^3(K))$ associated to $\spin(S_n^3(K))\cong \mathbb{Z}_n$:$$\widehat{HF}(S_n^3(K))=\bigoplus_{[s]\in\mathbb{Z}_n}\widehat{HF}(S_n^3(K),[s]).$$
\item The filtration on the Heegaard Floer chain complex $\widehat{CF}(S^3)$ associated to $K$, which induces a spectral sequence from $\widehat{HFK}(S^3,K)$ to $\widehat{HF}(S^3)$.
\item The large surgery formula computing $\widehat{HF}(S_n(K),[s])$ for any large integer $n$ and $[s]\in\mathbb{Z}_n$ from the filtrations associated to $K$ and $-K$.
\item The differential $D$ on the doubly-graded Heegaard Floer chain complex $CFK^\infty(S^3,K)$, in particular the fact that $D^2=0$.
\eenu
\efa
Since we will use bypass maps based on contact geometry throughout the paper, it is more convenient to use manifolds with reverse orientations. For technical reasons, we replace the notation $KHI$ with $\khii$. The constructions below can be generalized to a rationally null-homologous knot in a closed 3-manifold. For simplicity, we only discuss the constructions for a knot $K$ in an integral homology sphere $Y$ and deal with the general case in the main body of the paper. Suppose $S$ is a Seifert surface of $K$.

The analogy of term (1) can be found in \cite[Section 4]{LY2020}. We write the decomposition as $$I^\sharp(-Y_{-n}(K))=\bigoplus_{[s]\in\mathbb{Z}_{n}}I^\sharp(-Y_{-n}(K),[s]).$$

Since there is no explicit construction of the chain complex of $\khii(Y,K)$, it is hard to construct the filtration directly. Fortunately, it is possible to recover the spectral sequence and then lift the spectral sequence to a filtered chain complex by an algebraic construction. For the analog of term (2), we construct two spectral sequences from $\khii(-Y,K)$ to $I^\sharp(-Y)$ using two types of bypass maps, and construct two filtered differentials $d_+$ and $d_-$ on $\khii(-Y,K)$ with $$H(\khii(-Y,K),d_+)\cong H(\khii(-Y,K),d_-)\cong I^\sharp(-Y).$$

For the analog of term (3), we need to introduce the bent complex (cf.\ Construction \ref{cons: bent complex} and Construction \ref{cons: bent complex, dual}) as follows. 

For any integer $s$, the \textbf{bent complex} and the \textbf{dual bent complex} are the chain complexes \[A_s=A_s(-Y,K)\deq(\khii(-Y,K),d_s)\aand A_s^\vee=A_s^\vee(-Y,K)\deq(\khii(-Y,K),d_s^\vee),\]respectively, where for any element $x\in \khii(-Y,K,S,k)$, 
\[
d_s(x)=\begin{cases}
d_+(x)&k>0,\\
d_+(x)+d_-(x)&k=0,\\
d_-(x)&k<0,
 \end{cases}\aand d_s^\vee(x)=\begin{cases}
d_-(x)&k>0,\\
d_+(x)+d_-(x)&k=0,\\
d_+(x)&k<0.
 \end{cases}\] 
Since $d_+\circ d_+=d_-\circ d_-=0$, we have $d_s\circ d_s=d_s^\vee\circ d_s^\vee=0$. Hence, we can consider the homologies $H(A_s)$ and $H(A^\vee_s)$. The proof of the following theorem is purely algebraic. The main ingredient is the octahedral axiom for a triangulated category.

\bthm[Large surgery formula]\label{thm: large surgery formula, main}
For a fixed integer $n$ satisfying $|n|\ge 2g(K)+1$, suppose $$s_{min}=-|n|+1+g(K)\aand s_{max}=|n|-1-g(K).$$For any integer $s^\p$, suppose $[s^\p]$ is the image of $s^\p$ in $\mathbb{Z}_{|n|}$. For any integer $s\in [s_{min},s_{max}]$, we have
$$
    I^\sharp(-Y_{-n}(K),[s-s_{min}])\cong \begin{cases}
    H(A_{-s})&\text{ if } n>0,\\
    H(A_{-s}^\vee)&\text{ if } n<0.
    \end{cases}
$$
\ethm

We do not know how to construct the analog of the term (4). However, the proof of \cite[Theorem 1.2]{Ozsvath2005} only uses the fact that $D^2=0$ on some subcomplexes of $CFK^\infty(S^3,K)$. Thus, to obtain a proof of Theorem \ref{thm: main 2}, we only need some weaker vanishing results. Since the precise statement is too technical, we only state some byproducts in the next subsection, which are of independent interests for contact geometry.

\subsection{Instanton contact element and Giroux torsion}\quad

For a contact 3-manifold $(N,\xi)$ with convex boundary and dividing set $\Ga$ on $\partial N$, Baldwin and Sivek \cite{baldwin2016instanton} constructed an instanton contact element $\theta(N,\Ga,\xi)$ that lives in a version of sutured instanton homology $\shi(-N,-\Ga)$ \cite{baldwin2015naturality}. Suppose $(Y,\xi^\p)$ is a closed contact 3-manifold, and suppose $(Y(1),\delta,\xi^\p|_{Y(1)})$ is obtained from $(Y,\xi^\p)$ by removing a 3-ball. Then Baldwin and Sivek defined$$\theta(Y,\xi^\p)\deq \theta(Y(1),\delta,\xi^\p|_{Y(1)})\in \shi(-Y(1),-\delta)=I^\sharp(Y).$$ We have the following theorems for the instanton contact element.
\bthm\label{thm: contact support}
Suppose $(N,\xi)$ is a contact 3-manifold with convex boundary and dividing set $\Ga$ on $\partial N$. Suppose $S$ is an admissible surface (cf.\ Definition \ref{defn_2: admissible surfaces}) in $(N,\Ga)$, and suppose $S_+$ and $S_-$ are the positive region and the negative region of $S$ with respect to $\xi$, respectively. We write the $\mathbb{Z}$-grading associated to $S$ as $$\shi(-N,-\Ga)=\bigoplus_{i\in\mathbb{Z}}\shi(-N,-\Ga,S,i).$$Then the instanton contact element $\theta(N,\Ga,\xi)$ lives in $$\shi(-N,-\Ga,S,\frac{\chi(S_+)-\chi(S_-)}{2}).$$
\ethm
\bdefn\label{defn: giroux torsion}
A contact closed 3-manifold $(Y, \xi)$ has \textbf{Giroux torsion} if there is an embedding of $(T^2 \times [0,1], \eta_{2\pi})$ into $(Y,\xi)$, where $(x, y, t)$ are coordinates on $T^2 \times  [0, 1] \cong  \mathbb{R}^2/\mathbb{Z}^2 \times [0, 1]$ and $$\eta_{2\pi} = \ke(cos(2\pi t)dx-sin(2\pi t)dy).$$ 
\edefn
\bthm[]\label{thm: SHI giroux torsion}
If a closed contact 3-manifold $(Y, \xi)$ has Giroux torsion, then its instanton contact element $\theta(Y, \xi) \in I^\sharp(-Y)$ vanishes.
\ethm
\brem
There is a contact element in Heegaard Floer theory, constructed by Ozsv\'{a}th and Szab\'{o} \cite{OS05contact} for closed contact 3-manifolds, and extended by Honda, Kazez, and Mati\'{c} \cite{honda2009contact} for contact 3-manifolds with convex boundary. The analog of Theorem \ref{thm: contact support} in Heegaard Floer theory holds by definition of the contact element. The analog of Theorem \ref{thm: SHI giroux torsion} in Heegaard Floer theory was first conjectured by Ghiggini \cite[Conjecture 8.3]{girouxtorsion0}, and then proved by Ghiggini, Honda, and Van Horn-Morris \cite{girouxtorsion1}. More proofs can be found in \cite{girouxtorsion2,girouxtorsion3}.
\erem

\begin{org}
The paper is organized as follows. In Section \ref{sec: Algebraic preliminaries}, we collect some algebraic results about spectral sequences and the triangulated category, which are used in the proof of the large surgery formula. In Section \ref{sec: Large surgery formula}, we constructed differentials $d_+$ and $d_-$ on $\khii(Y,K)$ for a rationally null-homologous knot $K$ in a closed 3-manifold $Y$ and prove a generalization of Theorem \ref{thm: large surgery formula, main}. In Section \ref{sec: Vanishing results for contact elements}, we prove some vanishing results about contact elements and cobordism maps associated to contact structures. In particular, we prove Theorem \ref{thm: contact support} and Theorem \ref{thm: SHI giroux torsion}. In Section \ref{sec: Knots with instanton L-space surgeries}, we use results in former sections to prove a generalization of Theorem \ref{thm: main SU2}. Moreover, we prove Theorem \ref{main: L space knot} and Theorem \ref{thm: floer simple khi}. In Section \ref{sec: Dehn surgery along genus-one knots}, we study surgeries on genus-one knots in $S^3$ and prove Theorem \ref{thm: alternating}. In Section \ref{sec:Examples of SU(2)-abundant knots}, we provide examples of $SU(2)$-abundant knots and prove Corollary \ref{cor: abundant knot su2}. In Section \ref{sec: future}, we discuss some further directions of techniques introduced in this paper and make some conjectures.

\end{org}


\begin{conv}
If it is not mentioned, all manifolds are smooth, oriented, and connected. All contact structures are oriented and positively co-oriented. Classical homology groups and cohomology groups are defined with $\mathbb{Z}$ coefficients (while the instanton Floer homology in this paper is defined with $\mathbb{C}$ coefficients). We write $\mathbb{Z}_n$ for $\mathbb{Z}/n\mathbb{Z}$ and $\mathbb{F}_2$ for the field with two elements.

A knot $K\subset Y$ is called \textbf{null-homologous} if it represents the trivial homology class in $H_1(Y;\mathbb{Z})$, while it is called \textbf{rationally null-homologous} if it represents the trivial homology class in $H_1(Y;\mathbb{Q})$.

For any compact 3-manifold $M$, we write $-M$ for the manifold obtained from $M$ by reversing the orientation. For any surface $S$ in a compact 3-manifold $M$ and any suture $\ga\subset \partial M$, we write $S$ and $\ga$ for the same surface and suture in $-M$, without reversing their orientations. For a knot $K$ in a 3-manifold $Y$, we write $(-Y,K)$ for the induced knot in $-Y$ with induced orientation, called the \textbf{mirror knot} of $K$. The corresponding balanced sutured manifold is $(-Y(K),-\ga_K)$.

\end{conv} 
\begin{ack}
The authors would like to thank John A. Baldwin, Paolo Ghiggini, Ko Honda, Wenyuan Li, Ciprian Manolescu, Linsheng Wang, and  Yi Xie for valuable discussions. The authors are grateful to Ian Zemke for pointing out the proof of Proposition \ref{prop: vanishing cobordism}. The second author would like to thank his supervisor, Jacob Rasmussen, for patient guidance and helpful comments and thank his parents for support and constant encouragement. The second author is also grateful to Yi Liu for inviting him to BICMR at Peking University.
\end{ack}

\section{Algebraic preliminaries}\label{sec: Algebraic preliminaries}

In this section, we collect some algebraic results from homological algebra. All vector spaces are finite-dimensional and defined over a fixed field.

\subsection{Unrolled exact couples}\label{subsec: Unrolled exact couples}\quad

In this subsection, we explain the construction of the spectral sequence from an unrolled exact couple \cite{Boardman99} and describe the relationship between the spectral sequence and the filtered chain complex.

\bdefn\label{defn: unrolled exact couples}
An \textbf{unrolled exact couple} $(E^s,A^s)$ is a diagram of graded vector spaces and homomorphisms of the form
\begin{equation*}\label{eq: unrolled couples}
\xymatrix@R=6ex{
\cdots\ar[r]&A^{s+2}\ar[rr]^{i}&&A^{s+1}\ar[rr]^{i}\ar[dl]_{j}&&A^{s}\ar[rr]^{i}\ar[dl]_{j}&&A^{s-1}\ar[dl]_{j}\ar[r]&\cdots
\\&\cdots&E^{s+1}\ar[ul]_{k}&&E^s\ar[ul]_{k}&&E^{s-1}\ar[ul]_{k}&\cdots&
}
\end{equation*}
in which each triangle $$\cdots\ra A^{s+1}\ra A^s\ra E^s \ra A^{s+1}\ra \cdots$$ is a long exact sequence. An unrolled exact couple is called \textbf{bounded} by an interval $[s_1,s_2]$ if $E^s=0$ for $s\not\in [s_1,s_2]$. A morphism between two unrolled exact couples $(E^s,A^s)$ and $(\bar{E}^s,\bar{A}^s)$ consists of maps $f^s:E^s\to \bar{E}^s$ and $g^s:A^s\to \bar{A}^s$ that make all squares commute.

Suppose $(E^s,A^s)$ is an unrolled exact couple. For any integers $s$ and $r$, define$$\ke^r A^s=\ke (i^{(r)}:A^{s}\to A^{s-r})\aand \im^r A^s=\im (i^{(r)}:A^{s+r}\to A^s),$$where $i^{(r)}$ denotes the $r$-fold iterate of $i$. There are subgroups of $E^s$:
$$0=B_1^s\subset B_2^s\subset\cdots \subset \im j=\ke k\subset\cdots\subset Z_2^s\subset Z_1^s=E^s,$$where $$B_r^s=j(\ke^{r-1}A^s)\aand Z_r^s=k^{-1}(\im^{r-1} A^{s+1}).$$We call $B_r^s$ and $Z_r^s$ the $r$-th \textbf{boundary subgroup} and the $r$-th \textbf{cycle subgroup} of $E^s$, respectively. We call the quotient $$E_r^s=Z_r^s/B_r^s$$ the $s$-component of the $r$-th \textbf{page}. Note that $E_1^s=E^s$. If the unrolled exact couple is bounded by $[s_1,s_2]$, then we call the direct sum $$E_r=\bigoplus_{s_1}^{s_2} E_r^s$$the $r$-th \textbf{page}.
\edefn

\brem\label{rem: finite}
If the unrolled exact couple $(E^s,A^s)$ is bounded by $[s_1,s_2]$, then for any integers $r_1,r_2> s_2-s_1$ and any integer $s$, we have $$B^s_{r_1}=B^s_{r_2},Z^s_{r_1}=Z^s_{r_2}, E^s_{r_1}=E^s_{r_2}= E^s_\infty, \aand E_{r_1}=E_{r_2}=E_\infty.$$
\erem

\bprop[{\cite[Section 0]{Boardman99}}]\label{prop: spectral sequence}
Suppose $(E^s,A^s)$ is an unrolled exact couple. For any integers $s$ and $r$, there exists a well-defined map $$d_r^s: E_r^s\to E_r^{s+r}$$ induced by $j\circ (i^{(r-1)})^{-1}\circ k$ such that $$d_r^{s+r}\circ d_r^s=0\aand \ke d_r^s /\im d_{r}^{s-r}\cong E_{r+1}^s.$$ Equivalently, the set $\{(E_r^s,d_r^s)\}_{r\ge 1}$ forms a spectral sequence. Moreover, a morphism between two unrolled exact couples induces a map between the corresponding spectral sequences.
\eprop
Boardman studied the convergence of the spectral sequence in Proposition \ref{prop: spectral sequence} carefully, while we only need the special case for bounded unrolled exact couples.

\bthm[{\cite[Theorem 6.1]{Boardman99}}]\label{thm: convergence}
Suppose $(E^s,A^s)$ is an unrolled exact couple bounded by $[s_1,s_2]$. Then by exactness we have $$A^{s_1}\cong A^{s_1-1}\cong A^{s_1-2}\cong\cdots\aand A^{s_2+1}\cong A^{s_2+2}\cong A^{s_2+3}\cong\cdots$$Consider the spectral sequence $\{(E_r,d_r)\}_{r\ge 1}$ from Proposition \ref{prop: spectral sequence}, where we omit the superscript $s$ to denote the direct sum of all $s$-components. Then we have the following results.
\benu
    \item If $A^{s_1}=0$, then $\{(E_r,d_r)\}_{r\ge 1}$ converges to $G=A^{s_2+1}$ with filtration $F^sG=\ke^{s_2+1-s}A^{s_2+1}$, and we have $F^sG/F^{s+1}G\cong E_{\infty}^s$.
    \item If $A^{s_2+1}=0$, then $\{(E_r,d_r)\}_{r\ge 1}$ converges to $G=A^{s_1}$ with filtration $F^sG=\im^{s-s_1}A^{s_1}$, and we have $F^sG/F^{s+1}G\cong E_{\infty}^s$.
\eenu
\ethm

It is well-known that a filtered chain complex can induce a spectral sequence. Conversely, we may construct a filtered chain complex from a spectral sequence. However, a priori we may lose information when passing a filtered chain complex to a spectral sequence, so the reverse procedure is not always canonical. When fixing an inner product on the first page or equivalently fixing a basis, we have the following canonical construction. 
\bcons\label{cons: recover filtered chain complex}
Suppose $(E^s,A^s)$ is an unrolled exact couple bounded by $[s_1,s_2]$ and suppose $\{(E_r,d_r)\}_{r\ge 1}$ is the spectral sequence from Proposition \ref{prop: spectral sequence}. Fix an inner product on $E^s_1=E^s$ for all integers $s$. For simplicity, we omit the superscript $s$ and consider the direct sum $E$ of all $E^s$. 

For any subgroup $X$ of $E$, there is a canonical isomorphism $E/X\cong X^\perp$, where $X^\perp$ is the orthogonal complement of $X$ under the fixed inner product. From Definition \ref{defn: unrolled exact couples} and Remark \ref{rem: finite}, there are subgroups of $E$:
$$0=B_1\subset B_2\subset\cdots B_{s_2-s_1+1}\subset Z_{s_2-s_1+1}\subset\cdots\subset Z_2\subset Z_1=E.$$For $p=1,\dots,s_2-s_1$, define $B_p^\p$ as the orthogonal complement of $B_p$ in $B_{p+1}$, define $Z_{p}^\p$ as the orthogonal complement of $Z_{p+1}$ in $Z_{p}$, and define $E_{\infty}^\p$ as the orthogonal complement of $B_{s_2-s_1+1}^\p$ in $Z_{s_2-s_1+1}^\p$. Then we have
$$\begin{aligned}
    E_r&=Z_r/B_r\cong \bigoplus_{p=r}^{s_2-s_1}(B_p^\p\oplus Z_p^\p)\oplus E_{\infty}^\p, \\\ke d_r&=Z_{r+1}/B_r\cong \bigoplus_{p=r+1}^{s_2-s_1}(B_p^\p\oplus Z_p^\p)\oplus E_{\infty}^\p\oplus B_r^\p,\\
    \im d_r&=B_{r+1}/B_r\cong B_r^\p
\end{aligned}$$
Hence we can lift $d_r:E_r\to E_r$ to a map $$d_r^\p=I\circ d_r\circ P:E\to E,$$where $P$ and $I$ are the projection and the inclusion, respectively. The only nontrivial part of $d_r^\p$ is from $Z_r^\p$ to $B_r^\p$, so for any $r_1,r_2\in\{1,\dots,s_2-s_1\}$, we have $d_{r_1}^\p\circ  d_{r_2}^\p=0$. Hence, the summation $$d=\sum_{r=1}^{s_2-s_1}d_r^\p$$ is a differential on $E$, i.e.\, $d^2=0$. Moreover, we have $$H(E,d)\cong E_\infty^\p\cong E_{s_2-s_1+1}\cong E_\infty.$$
It is straightforward to check that the filtration $F^sE=\bigoplus_{p\ge s}E^p$ on $(E,d)$ induces the spectral sequence $\{(E_r,d_r)\}_{r\ge 1}$.
\econs

\subsection{The octahedral axiom}\quad

It is well-known that the derived category of an abelian category is a triangulated category (for example, see \cite[Proposition 10.2.4]{homologicalalgebra94}). In particular, the derived category of the category of vector spaces is triangulated. Graded vector spaces can be regarded as objects in the derived category with trivial differentials. The following theorem is a special case of the octahedral axiom of the triangulated category. 

\bthm\label{thm: octahedral axiom}
Suppose $X,Y,Z,X^\p,Y^\p,Z^\p$ are graded vector spaces satisfying long exact sequences$$\begin{aligned}
        \cdots&\to X\xra{f} Y\xra{h} Z^\p\to X\{1\}\to\cdots\\
        \cdots&\to Y\xra{g} Z\to X^\p\xra{l} Y\{1\}\to\cdots\\
        \cdots&\to X\xra{g\circ f} Z\xra{j} Y^\p\to X\{1\}\to\cdots
    \end{aligned}$$where $X\{1\}$ denotes the grading shift of $X$ by $1$; so do $Y\{1\}$ and $Z\{1\}$. Then we have the fourth long exact sequence$$\cdots\to Z^\p\xra{\psi} Y^\p\xra{\phi} X^\p \xra{h\{1\}\circ l} Z^\p\{1\}\to\cdots$$such that the following diagram commutes
    \begin{equation*}\label{eq: octahedral}
    \xymatrix@R=6ex{
    &&&Z^\p\ar[dr]^{\psi}\ar[rr]&&X\{1\}\ar[ddr]^{f\{1\}}&&
    \\
    &&Y\ar[dr]^{g}\ar[ur]&&Y^\p\ar[ddrr]^{\phi}\ar[ur]&&&
    \\
 &&&Z\ar[drrr]\ar[ur]^{j}&&&Y\{1\}\ar[dr]^{h\{1\}}&
    \\
    X\ar[uurr]^{f}\ar[urrr]^{g\circ f}&&&&&&X^\p\ar[r]^{ h\{1\}\circ l}\ar[u]^{l}&Z^\p\{1\}
    }
    \end{equation*}
    
      
    where the arrows come from four long exact sequences.
\ethm
\bpf[Sketch of the proof]
We regard graded vector spaces as chain complexes with trivial differentials. By the long exact sequences in the assumption, we know that $Z^\p, X^\p,Y^\p$ are chain homotopic to the mapping cones $\cone(f),\cone(g),\cone(g\circ f)$, respectively. Define $$\begin{aligned}
    \psi:Y\oplus X\{1\}&\to Z\oplus X\{1\}\\
    \psi(y,x)&\mapsto (g(y),x)
\end{aligned}$$and$$\begin{aligned}
    \phi:Z\oplus X\{1\}&\to Z\oplus Y\{1\}\\
    \phi(z,x)&\mapsto (z,f\{1\}(x))
\end{aligned}$$The map $\psi$ is a chain map from $\cone(f)$ to $\cone(g\circ f)$ and the map $\phi$ is a chain map from $\cone(g\circ f)$ to $\cone(g)$. Since the underlying vector space of $\cone(\psi)$ is $Z\oplus X\{1\}\oplus Y\{1\}\oplus X\{2\}$, the inclusion $Z\oplus Y\{1\}\to Z\oplus X\{1\}\oplus Y\{1\}\oplus X\{2\}$ induces a map $\eta$ from $\cone(g)$ to $\cone(\psi)$, which is a chain map and makes the following diagram commute
\begin{equation*}
    \xymatrix@R=6ex{
\cone(f)\ar[rr]^{\psi}\ar[d]^{=}&& \cone(g\circ f)\ar[d]^{=}\ar[rr]^{\phi}&&\cone(g)\ar[rr]^{h\{1\}\circ l}\ar[d]^{\eta}&&\cone(f)\{1\}\ar[d]^{=}\\
    \cone(f)\ar[rr]^{\psi}&& \cone(g\circ f)\ar[rr]&&\cone(\psi)\ar[rr]&&\cone(f)\{1\}
    }
\end{equation*}
Define 
$$\begin{aligned}
    \zeta:Z\oplus X\{1\}\oplus Y\{1\}\oplus X\{2\}&\to Z\oplus Y\{1\}\\
    \zeta(z,x,y,x^\p)&\mapsto (z,y+f\{1\}(x))
\end{aligned}$$Then we can check that $\zeta\circ \eta$ is the identity map on $\cone(g)$ and $\eta\circ\zeta$ is chain homotopic to the identity on $\cone(\psi)$. Hence $\cone(f),\cone(g\circ f)$ and $\cone(g)$ form a long exact sequence.
\epf


    Note that the chain homotopies in the proof of Theorem \ref{thm: octahedral axiom} are not canonical, and hence the maps $\psi$ and $\phi$ are also not canonical. Thus, we usually cannot identify them with other given maps $\psi^\p,\phi^\p$. However, in the special case that $\phi\circ j=\phi^\p\circ j=0$, it is possible to identify $\phi$ and $\phi^\p$ by the following lemma. 

\blem\label{lem: identify phi}
Suppose $X,Y, Z,X^\p,Y^\p$ are graded vector spaces satisfying the following horizontal exact sequences. 
\begin{equation*}
\xymatrix@R=3ex{
Z\ar[rr]^{j}\ar[dd]^{=}&&Y^\p\ar[rr]^{l^\p}\ar@<-.5ex>[dd]_{\phi} \ar@<.5ex>[dd]^{\phi^{\prime}}&&X\{1\}\ar[dd]^{f\{1\}}\\
&&\\
Z\ar[rr]^{0}&&X^\p\ar[rr]^{l}&&Y\{1\}
}
\end{equation*}
Suppose $\phi:Y^\p\to X^\p$ satisfies the two commutative diagrams, i.e.\, $\phi\circ j=0$ and $f\{1\}\circ l^\p=l\circ \phi$. Suppose $\phi^\p:Y^\p\to X^\p$ satisfies the two commutative diagrams up to a unit, i.e.\, $\phi^\p\circ j=0$ and $f\{1\}\circ l^\p=c \cdot l\circ \phi^\p$ for some $c\in\mathbb{C}\backslash\{0\}$. Then we have $\phi\doteq \phi^\p$ and hence $H(\cone(\phi))\cong H(\cone(\phi^\p))$.
\elem
\bpf
By exactness at $X^\p$, we have $$\im(\phi-c\phi^\p)=\ke(l)=\im(0)=0.$$Hence $\phi=c\phi^\p$.
\epf

\section{Differentials and the large surgery formula}\label{sec: Large surgery formula}

In this section, we provide more details for constructions in Subsection \ref{subsec: Constructions parallel to Heegaard Floer theory} and prove Theorem \ref{thm: large surgery formula, main}. Most notations follow from \cite[Section 4]{LY2020}.

\subsection{Backgrounds on sutured instanton homology}\label{subsec: Backgrounds on sutured instanton homology}\quad

In this subsection, we review some basic facts of sutured instanton homology.

\bdefn[{\cite[Definition 2.2]{juhasz2006holomorphic}}]\label{defn_2: balanced sutured manifold}
A \textbf{balanced sutured manifold} $(M,\ga)$ consists of a compact 3-manifold $M$ with non-empty boundary together with a closed 1-submanifold $\ga$ on $\partial{M}$. Let $A(\ga)=[-1,1]\times\ga$ be an annular neighborhood of $\ga\subset \partial{M}$ and let $R(\ga)=\partial{M}\backslash{\rm int}(A(\ga))$, such that they satisfy the following properties.
\benu
    \item Neither $M$ nor $R(\ga)$ has a closed component.
    \item If $\partial{A(\ga)}=-\partial{R(\ga)}$ is oriented in the same way as $\ga$, then we require this orientation of $\partial{R(\ga)}$ induces the orientation on $R(\ga)$, which is called the \textbf{canonical orientation}.
    \item Let $R_+(\ga)$ be the part of $R(\ga)$ for which the canonical orientation coincides with the induced orientation on $\partial{M}$ from $M$, and let $R_-(\ga)=R(\ga)\backslash R_+(\ga)$. We require that $\chi(R_+(\ga))=\chi(R_-(\ga))$. If $\ga$ is clear in the contents, we simply write $R_\pm=R_\pm(\ga)$, respectively.
\eenu

For any balanced sutured manifold $(M,\ga)$, Kronheimer and Mrowka \cite[Section 7]{kronheimer2010knots} constructed a $\mathbb{C}$-vector space $SHI(M,\ga)$ called the \textbf{sutured instanton homology} of $(M,\ga)$. The construction was based on closures of $(M,\ga)$, i.e.\, a tuple $(Y,R,\omega)$ consisting of a closed 3-manifold $Y$, a closed surface $R\subset Y$, and a 1-cycle $\omega\subset Y$ with some admissible conditions. 

A priori, the space $SHI(M,\ga)$ only represents an isomorphism class. Later, Baldwin and Sivek \cite[Section 9]{baldwin2015naturality} dealt with the naturality issue and constructed a projectively transitive system $\shiib(M,\ga)$ (twisted version). This system records the collection of vector spaces associated to different closures of $(M,\ga)$, which are all isomorphic to $SHI(M,\ga)$, together with canonical isomorphisms relating these spaces, where these isomorphisms are well-defined up to multiplication by a unit in $\mathbb{C}$. 

In practice, when considering maps between sutured instanton homology, we can always fix closures of corresponding balanced sutured manifolds and consider linear maps between actual vector spaces, at the cost that equations between maps only hold up to multiplication by a unit. Hence if it is clear, we will not distinguish between the projectively transitive system and the vector space in the system.

To be consistent with the notations in \cite{LY2020}, we write $\shi(M,\ga)$ for the system $\shiib(M,\ga)$. Note that $\shi(M,\ga)$ represents the isomorphism class in \cite[Section 9]{baldwin2015naturality}, and we write $SHI(M,\ga)$ for the isomorphism class instead. 

There is another projectively transitive system $\shib^g(M,\ga)$ (untwisted version) constructed in \cite[Section 9]{baldwin2015naturality}. The main difference of two systems is that $\shi(M,\ga)$ corresponds to closures of $(M,\ga)$ for which the surface $R$ may have different genera and $\shib^g(M,\ga)$ corresponds to closures for which $g=g(R)$ is fixed. Many arguments for $\shi(M,\ga)$ also hold for $\shib^g(M,\ga)$ when $g$ is sufficiently large. In \cite{LY2021}, we considered $\shib^g(M,\ga)$ as a special case of formal sutured homology and calculated its graded Euler characteristic for sufficiently large $g$. By \cite[Theorem 9]{baldwin2015naturality}, the subsystem of $\shi(M,\ga)$ for closures of fixed genus $g$ is isomorphic to $\shib^g(M,\ga)$, so properties of $\shg^g(M,\ga)$ (especially about graded Euler characteristics) also apply to $\shi(M,\ga)$.

Suppose $K$ is a knot in a closed 3-manifold $Y$. Let $$Y(1)\deq Y\backslash {\rm int} B^3\aand Y(K)\deq Y\backslash {\rm int} N(K).$$Suppose $\delta$ is a simple closed curve on $\partial Y(1)\cong S^2$ and suppose $\ga_K$ is two copies of the meridian of $K$ with opposite orientations. Define\[I^\sharp(Y)\deq \shi(Y(1),\delta)\aand \khii(Y,K)\deq \shi(Y(K),\ga_K).\]Note that $I^\sharp(Y)$ also denotes the framed instanton homology of $Y$ constructed in \cite{kronheimer2011khovanov}, though it is isomorphic to $\shi(Y(1),\delta)$. So we abuse notation and do not distinguish these two definitions in this paper.

\edefn
\bdefn[{\cite[Definition 2.26]{li2019decomposition}}]\label{defn_2: admissible surfaces}
Suppose $(M,\ga)$ is a balanced sutured manifold and $S\subset (M,\ga)$ is a properly embedded surface in $M$. The surface $S$ is called an \textbf{admissible surface} if the following conditions hold.
\benu
    \item Every boundary component of $S$ intersects $\ga$ transversely and nontrivially.
    \item We require that $\frac{1}{2}|S\cap \ga|-\chi(S)$ is an even integer.
\eenu
\edefn
For an admissible surface $S\subset (M,\ga)$, there is a $\mathbb{Z}$-grading on $\shi(M,\ga)$ \cite{li2019direct,li2019decomposition}:\[\shi(M,\ga)=\bigoplus_{i\in\mathbb{Z}}\shi(M,\ga,S,i).\]From the construction of the grading, we have the following basic proposition, which implies (\ref{eq: mirror iso}).
\bprop[{\cite[Proposition 2.30]{LY2021}}]\label{prop: mirror}
For any balanced sutured manifold $(M,\ga)$ and any admissible surface $S\subset (M,\ga)$, there are canonical isomorphisms$$\shi(-M,\ga,S,i)\cong {\rm Hom}_\mathbb{C}(\shi(M,\ga,S,i),\mathbb{C})$$and$$\shi(M,\ga,-S,i)\cong \shi(M,-\ga,S,i)\cong \shi(M,\ga,S,-i).$$
\eprop

\subsection{The canonical basis on the torus boundary}\label{subsec: basic setups}\quad

In this subsection, we provide a canonical way to fix the basis on the boundary of the knot complement and introduce some notations about sutures.

Suppose $Y$ is a closed 3-manifold and $K\subset Y$ is a null-homologous knot. Let $Y(K)$ be the knot complement $Y\backslash {\rm int}(N(K))$. Any Seifert surface $S$ of $K$ gives rise to a framing on $\partial Y(K)$: the longitude $\lambda$ can be picked as $S\cap \partial Y(K)$ with the induced orientation from $S$, and the meridian $\mu$ can be picked as the meridian of the solid torus $N(K)$ with the orientation such that $\mu\cdot \lambda=-1$. The `half lives and half dies' fact for 3-manifolds implies that the following map has a 1-dimensional image:
$$\partial_*: H_2(Y(K),\partial Y(K);\mathbb{Q})\ra H_1(\partial Y(K);\mathbb{Q}).$$
Hence any two Seifert surfaces lead to the same framing on $\partial Y(K)$.

\bdefn\label{defn_4: surgery on pairs}
The framing $(\mu,\lambda)$ defined as above is called the \textbf{canonical framing} of $(Y,K)$. With respect to this canonical framing, let
$$\widehat{Y}_{q/p}=Y(K)\cup_{\phi}S^1\times D^2$$
be the 3-manifold obtained from $Y$ by a $q/p$ surgery along $K$, i.e.\, $$\phi(\{1\}\times \partial D^2)=q\mu+p\lambda.$$

We also write $\widehat{Y}_\al$ for $\widehat{Y}_{q/p}$, where $\al=\phi(\{1\}\times \partial D^2)$. When the surgery slope is understood, we also write $\widehat{Y}_{q/p}$ simply as $\widehat{Y}$. Let $\widehat{K}$ be the dual knot, i.e.\, the image of $S^1\times\{0\}\subset S^1\times D^2$ in $\widehat{Y}$ under the gluing map.
\edefn

\begin{conv}
Throughout this section, we will always assume that either ${\rm gcd}(p,q)=1$ with $q> 0$ or $(p,q)=(1,0)$ for a Dehn surgery. In particular, the original pair $(Y,K)$ can be thought of as a pair $(\widehat{Y},\widehat{K})$ obtained from $(Y,K)$ by the $1/0$ surgery. Moreover, we will always assume that the knot complement $Y(K)$ is irreducible. This is because if $Y(K)$ is not irreducible, then $Y(K)\cong  Y^\p(K^\p)\sharp Y^\pp$ for some closed 3-manifolds $Y^\p,Y^\pp$ and a null-homologous knot $K^\p\subset{Y^\p}$. By the connected sum formula \cite[Proposition 4.15]{li2018contact}, we have $$\shi(Y(K),\ga)\cong \shi(Y^\p(K^\p),\ga)\otimes I^\sharp(Y^\pp)$$for any suture $\ga$. Hence all results hold after tensoring $I^\sharp(Y^\pp)$.
\end{conv}

Next, we describe various families of sutures on the knot complement. Suppose $K\subset Y$ is a null-homologous knot and the pair $(\widehat{Y},\widehat{K})$ is obtained from $(Y,K)$ by a $q/p$ surgery. Note that we can identify the complement of $K\subset Y$ with that of $\widehat{K}\subset \widehat{Y}$, i.e.\, $\widehat{Y}(\widehat{K})=Y(K).$

On $\partial{Y(K)}$, there are two framings: one comes from $K$, and we write longitude and meridian as $\lambda$ and $\mu$, respectively. The other comes from $\widehat{K}$. Note that only the meridian $\hat{\mu}$ of $\widehat{K}$ is well-defined, and by definition, it is $\hat{\mu}=q\mu+p\lambda.$
\bdefn\label{defn_4: Ga_n-hat sutures}
If $p=0$, then $q=1$ and $\hat{\mu}=\mu$. We can take $\hat{\lambda}=\lambda$. If $(q,p)=(0,1)$, then we take $\hat{\lambda}=-\mu$. If $p,q\neq 0$, then we take $\hat{\lambda}=q_0\mu+p_0\lambda$, where $(q_0,p_0)$ is the unique pair of integers such that the following conditions are true.
\benu
\item $0\le |p_0|<|p|$ and $p_0p\le 0$.
\item $0\le |q_0|<|q|$ and $q_0q\le 0$.
\item $p_0q-pq_0=1$.
\eenu
In particular, if $(q,p)=(n,1)$, then $\hat{\lambda}=-\mu$.

For a homology class $x\lambda+y\mu$, let $\ga_{x\lambda+y\mu}$ be the suture consisting of two disjoint simple closed curves representing $\pm(x\lambda+y\mu)$ on $\partial{Y(K)}$. Furthermore, for $n\in\intg$, define
$$\widehat{\Ga}_{n}(q/p)=\ga_{\hat{\lambda}-n\hat{\mu}}=\ga_{(p_0-np)\lambda+(q_0-nq)\mu},\aand \widehat{\Ga}_{\mu}(q/p)=\ga_{\hat{\mu}}=\ga_{p\lambda+q\mu}.$$
Suppose $(q_n,p_n)\in\{\pm(q_0-nq,p_0-np)\}$ such that $q_n\ge 0$. Note that there might be a sign ambiguity of $q_0$: if $q>0$, then by term (2) above $q_0<0$; but here $n=0$ implies the new $q_0$ is the opposite number of the original $q_0$. We keep this ambiguity and use the first definition of $q_0$ only for $\hat{\lambda}$ and uses the second definition only in the formula of $q_n$.

When emphasizing the choice of $\hat{\mu}$, we also write $\widehat{\Ga}_n(\hat{\mu})$ and $\widehat{\Ga}_{\mu}(\hat{\mu})$. When $\hat{\lambda}$ and $\hat{\mu}$ are understood, we omit the slope $q/p$ and simply write $\widehat{\Ga}_n$ and $\widehat{\Ga}_{\mu}$. When $(q,p)=(1,0)$, we write $\Ga_n$ and $\Ga_\mu$ instead.
\edefn

\brem\label{rem: remark on indexing the sutures}
Since the two components of the suture must be given opposite orientations, the notations $\ga_{x\lambda+y\mu}$ and $\ga_{-x\lambda-y\mu}$ represent the same suture on the knot complement $Y(K)$. Our choice makes $q_{n+1}\le q_n$ for $n<-1$ and $q_{n+1}\ge q_n$ for $n\ge 0$.
\erem

\subsection{Bypass maps on the knot complements}\label{subsec: Sutured instanton homology for knot complements}\quad

In this subsection, we review results in \cite[Section 4]{LY2020} that are useful in this paper. 

If $(M,\ga)=(Y(K),\ga_{x\lambda+y\mu})$ and $S$ is an admissible surface obtained from a minimal genus Seifert surface (cf.\ \cite[Definition 4.10]{LY2020}, where we write $S^\tau$ for $\tau\in\{0,-1\}$), then we can calculate the maximal and minimal nontrivial gradings explicitly. Note that we assume that $Y(K)$ is irreducible, so the decompositions of $(Y(K),\ga_{x\lambda+y\mu})$ along $S$ and $-S$ are both taut (cf.\ \cite[Definition 2.6]{juhasz2006holomorphic}). Since we will use contact gluing maps later, it is more convenient to consider $(-M,-\ga)$ instead of $(M,\ga)$.
\bdefn\label{defn_4: top and bottom non-vanishing gradings}
For any integer $y\in\mathbb{N}$, define
$$i_{max}^y=\lceil \frac{y-1}{2}\rceil+g(K),{~\rm and~}i_{min}^y=\lceil -\frac{y-1}{2}\rceil-g(K),$$where $\lceil x\rceil$ is the smallest integer larger than $x$. For $\hat{\mu}=q\mu+p\lambda$ and the sutures $\widehat{\Ga}_n$ and $\widehat{\Ga}_\mu$, define
$$\hat{i}^n_{max}=i_{max}^{q_n},\hat{i}^n_{min}=i^{q_n}_{min},\aand \hat{i}^\mu_{max}=i^{q}_{min},\hat{i}^\mu_{min}=i^{q}_{min}.$$
\edefn
\blem[{\cite[Lemma 4.12]{LY2020}}]\label{lem_4: top and bottom non-vanishing gradings}
Suppose $K\subset Y$ is a null-homologous knot and $\ga_{x\lambda+y\mu}$ is a suture on $\partial Y(K)$ with $y\ge 0$. Suppose further that $S$ is a Seifert surface of $K$. Then the maximal and minimal nontrivial gradings of $\shi(-Y(K),-\ga_{(x,y)})$ associated to $S$ are $i^y_{max}$ and $i^y_{min}$, respectively. In particular, the maximal and minimal nontrivial gradings of $\shi(-Y(K),-\widehat{\Ga}_n)$ associated to $S$ are $\hat{i}^n_{max}$ and $\hat{i}^n_{min}$, respectively.
\elem
It is easy to see that $$\lim_{n\to +\infty}(\hat{i}^n_{max}-\hat{i}^n_{min})=\lim_{n\to +\infty}(2g(K)+nq-q_0-1)=+\infty.$$However, by following lemmas, there is no more information in $\shi(-Y(K),-\widehat{\Ga}_n)$ when $n$ is large. To see this, we first introduce the bypass exact triangles.

\bdefn\label{defn_3: shifting the grading}
Suppose $(M,\ga)$ is a balanced sutured manifold and $S$ is an admissible surface in $(M,\ga)$. For any $i,j\in\intg$, define
$$\shi(M,\ga,S,i)[j]=\shi(M,\ga,S,i-j).$$Moreover, let $\shi(M,\ga,S,i)\{1\}$ be obtained from $\shi(M,\ga,S,i)$ by switching the odd and the even relative $\mathbb{Z}_2$-gradings.
\edefn

\bprop[{\cite[Proposition 4.14]{LY2020}, see also \cite[Proposition 5.5]{li2019direct}}]\label{prop_4: graded bypass for Ga_n-hat sutures}
Suppose $K\subset Y$ is a null-homologous knot and suppose the pair $(\widehat{Y},\widehat{K})$ is obtained from $(Y,K)$ by a $q/p$ surgery. Suppose further that the sutures $\widehat{\Ga}_n$ and $\widehat{\Ga}_{\mu}$ are defined as in Definition \ref{defn_4: Ga_n-hat sutures} and $S$ is a Seifert surface of $K$. Then the following conditions hold, where all maps are grading-preserving.
\benu
\item For $n\in\intg$ such that $q_{n+1}=q_n+q$, i.e.\, $n\ge 0$, there are two bypass exact triangles:
\begin{equation*}\label{eq_4: positive bypass, graded}
\xymatrix{
\shi(-Y(K),-\widehat{\Ga}_{n},S)[\hat{i}_{min}^{n+1}-\hat{i}_{min}^{n}]\ar[rr]^{\quad\quad\psi^{n}_{+,n+1}}&&\shi(-Y(K),-\widehat{\Ga}_{n+1},S)\ar[dll]^{\psi^{n+1}_{+,\mu}}\\
\shi(-Y(K),-\widehat{\Ga}_{\mu},S)[\hat{i}_{max}^{n+1}-\hat{i}_{max}^{\mu}]\ar[u]^{\psi^{\mu}_{+,n}}&&
}
\end{equation*}
and
\begin{equation*}\label{eq_4: negative bypass, graded}
\xymatrix{
\shi(-Y(K),-\widehat{\Ga}_{n},S)[\hat{i}_{max}^{n+1}-\hat{i}_{max}^n]\ar[rr]^{\quad\quad\psi^{n}_{-,n+1}}&&\shi(-Y(K),-\widehat{\Ga}_{n+1},S)\ar[dll]^{\psi^{n+1}_{-,\mu}}\\
\shi(-Y(K),-\widehat{\Ga}_{\mu},S)[\hat{i}_{min}^{n+1}-\hat{i}_{min}^{\mu}]\ar[u]^{\psi^{\mu}_{-,n}}&&
}
\end{equation*}

\item For $n\in\intg$ such that $q_{n+1}=q_n-q$, i.e.\, $n< -1$, there are two bypass exact triangles:
\begin{equation*}\label{eq_4: positive bypass, graded, 2}
\xymatrix{
\shi(-Y(K),-\widehat{\Ga}_{n},S)\ar[rr]^{\psi^{n}_{+,n+1}\quad\quad}&&\shi(-Y(K),-\widehat{\Ga}_{n+1},S)[\hat{i}_{max}^{n}-\hat{i}_{max}^{n+1}]\ar[dll]^{\psi^{n+1}_{+,\mu}}\\
\shi(-Y(K),-\widehat{\Ga}_{\mu},S)[\hat{i}_{min}^{n}-\hat{i}_{min}^{\mu}]\ar[u]^{\psi^{\mu}_{+,n}}&&
}
\end{equation*}
and
\begin{equation}\label{eq_4: negative bypass, graded, 2}
\xymatrix{
\shi(-Y(K),-\widehat{\Ga}_{n},S)\ar[rr]^{\psi^{n}_{-,n+1}\quad\quad}&&\shi(-Y(K),-\widehat{\Ga}_{n+1},S)[\hat{i}_{min}^{n}-\hat{i}_{min}^{n+1}]\ar[dll]^{\psi^{n+1}_{-,\mu}}\\
\shi(-Y(K),-\widehat{\Ga}_{\mu},S)[\hat{i}_{max}^{n}-\hat{i}_{max}^{\mu}]\ar[u]^{\psi^{\mu}_{-,n}}&&
}
\end{equation}

\item For $n\in\intg$ such that $q_{n+1}+q_n=q$, i.e.\, $n=-1$, there are two bypass exact triangles:
\begin{equation}\label{eq_4: positve bypass, graded, 3}
\xymatrix{
\shi(-Y(K),-\widehat{\Ga}_{n},S)[\hat{i}_{min}^{\mu}-\hat{i}_{min}^n]\ar[rr]^{\psi^{n}_{+,n+1}}&&\shi(-Y(K),-\widehat{\Ga}_{n+1},S)[\hat{i}_{max}^{\mu}-\hat{i}_{max}^{n+1}]\ar[dll]^{\psi^{n+1}_{+,\mu}}\\
\shi(-Y(K),-\widehat{\Ga}_{\mu},S)\ar[u]^{\psi^{\mu}_{+,n}}&&
}
\end{equation}
and
\begin{equation}\label{eq_4: negative bypass, graded, 3}
\xymatrix{
\shi(-Y(K),-\widehat{\Ga}_{n},S)[\hat{i}_{max}^{\mu}-\hat{i}_{max}^{n}]\ar[rr]^{\psi^{n}_{-,n+1}}&&\shi(-Y(K),-\widehat{\Ga}_{n+1},S)[\hat{i}_{min}^{\mu}-\hat{i}_{min}^{n+1}]\ar[dll]^{\psi^{n+1}_{-,\mu}}\\
\shi(-Y(K),-\widehat{\Ga}_{\mu},S)\ar[u]^{\psi^{\mu}_{-,n}}&&
}
\end{equation}
\eenu
\eprop

\brem
The maps $\psi_{+,*}^{*}$ and $\psi_{-,*}^{*}$ are called \textbf{bypass maps}, which are contact gluing maps induced by bypass attachments on balanced sutured manifolds. The exact triangles in Proposition \ref{prop_4: graded bypass for Ga_n-hat sutures} are called \textbf{bypass exact triangles}. In this paper, we will omit the definitions and focus on their algebraic properties.
\erem

\begin{lem}[{\cite[Lemma 4.34]{LY2020}}]\label{lem: commutative diagram 2}
For any surgery slope $q/p$, consider the bypass maps $\psi_{+,*}^*$ and $\psi_{-,*}^*$ in Proposition \ref{prop_4: graded bypass for Ga_n-hat sutures}. For any $n\in\mathbb{Z}$, we have two commutative diagrams
\begin{equation}\label{eq: commutative diagram 2a}
\xymatrix@R=6ex{
\shi(-Y(K),-\widehat{\Ga}_{n})\ar[rd]_{\psi_{+,\mu}^n}\ar[rr]^{\psi_{-,n+1}^n}&&\shi(-Y(K),-\widehat{\Ga}_{n+1})\ar[ld]^{\psi_{+,\mu}^{n+1}}\\
&\shi(-\widehat{Y}(K),-\widehat{\Ga}_{\mu})&
}
\end{equation}
and
\begin{equation}\label{eq: commutative diagram 2b}
\xymatrix@R=6ex{
\shi(-Y(K),-\widehat{\Ga}_{n})\ar[rr]^{\psi_{-,n+1}^n}&&\shi(-Y(K),-\widehat{\Ga}_{n+1})\\
&\shi(-\widehat{Y}(K),-\widehat{\Ga}_{\mu})\ar[ul]^{\psi_{+,n}^\mu}\ar[ur]_{\psi_{+,n+1}^\mu}&
}
\end{equation}
Similar commutative diagrams hold if we switch the roles of $\psi_{+,*}^*$ and $\psi_{-,*}^*$.
\end{lem}

In the following lemma, we abuse the notations for bypass maps so they also denote the restrictions on some gradings associated to $S$. 
\blem[{\cite[Lemma 4.16]{LY2020}}]\label{lem_4: psi pm is isomorphism}
For any $n\in\mathbb{N}$, the map 
$$\psi_{+,n+1}^{n}: \shi(-Y(K),-\widehat{\Ga}_{n},S,i)\ra\shi(-Y(K),-\widehat{\Ga}_{n+1},S,i-\hat{i}_{min}^n+\hat{i}_{min}^{n+1})$$is an isomorphism if $i\leq \hat{i}_{max}^n-2g(K)$. Similarly, the map
$$\psi_{-,n+1}^{n}: \shi(-Y(K),-\widehat{\Ga}_{n},S,i)\ra\shi(-Y(K),-\widehat{\Ga}_{n+1},S,i-\hat{i}_{max}^n+\hat{i}_{max}^{n+1})$$is an isomorphism if $i\geq \hat{i}_{min}^n+2g(K)$.
\elem

\blem[{\cite[Lemma 4.20]{LY2020}}]\label{lem_4: q-cyclic}
Suppose $n\in\mathbb{N}$ satisfies $q_n\geq q+2g(K)$, and suppose $i,j\in \intg$ with
$$\hat{i}_{min}^n+2g(K)\leq i,j\leq \hat{i}_{max}^n-2g(K),{\rm~and~}i-j=q.$$
Then we have
$$\shi(-Y(K),-\widehat{\Ga}_n,S,i)\cong \shi(-Y(K),-\widehat{\Ga}_n,S,j).$$
\elem
Thus, we can divide $\shi(-Y(K),-\widehat{\Ga}_n)$ into three parts: the top $2g(K)$ gradings, the middle gradings, and the bottom $2g(K)$ gradings. All parts stabilize by Lemma \ref{lem_4: psi pm is isomorphism} and the spaces in the middle gradings are cyclic by Lemma \ref{lem_4: q-cyclic}. Moreover, by Proposition \ref{prop: mirror}, we have a canonical isomorphism\begin{equation*}
    \shi(-M,-\ga,S,i)\cong \shi(-M,\ga,S,-i).
\end{equation*}If $\partial M\cong T^2$, we can identify $-\ga$ with $\ga$, which induces an isomorphism
\begin{equation}\label{eq: involution}
    \iota_{\ga}:\shi(-M,-\ga,S,i)\xra{\cong} \shi(-M,\ga,S,-i)\xra{=}\shi(-M,-\ga,S,-i).
\end{equation}
Hence the spaces in the top $2g(K)$ gradings and the bottom $2g(K)$ gradings are isomorphic. The following theorems imply that spaces in the middle gradings encode information of  $I^\sharp(-\widehat{Y})$.

\blem[{\cite[Lemma 4.9]{LY2020}, see also \cite[Section 3]{li2019tau}}]\label{lem_4: surgery triangle for knots}
Suppose $K\subset Y$ is a null-homologous knot and suppose the pair $(\widehat{Y},\widehat{K})$ is obtained from $(Y,K)$ by a $q/p$ surgery. Suppose further that the sutures $\widehat{\Ga}_n$ are defined as in Definition \ref{defn_4: Ga_n-hat sutures}. Then, there is an exact triangle
\begin{equation}\label{eq_2: surgery exact triangle}
\xymatrix@R=6ex{
\shi(-Y(K),-\widehat{\Ga}_{n})\ar[rr]&&\shi(-Y(K),-\widehat{\Ga}_{n+1})\ar[dl]^{F_{n+1}}\\
&I^\sharp(-\widehat{Y})\ar[ul]^{G_n}&
}
\end{equation}
where $F_n$ are the contact gluing maps associated to the contact 2-handle attachment along $\hat{\mu}=q\mu+p\lambda \subset \partial Y(K)$. Furthermore, we have four commutative diagrams related to $\psi^{n}_{+,n+1}$ and $\psi^{n}_{-,n+1}$, respectively
\begin{equation*}\label{eq_4: commutative diagram, G, +}
\xymatrix@R=6ex{
\shi(-Y(K),-\widehat{\Ga}_{n})\ar[rr]^{\psi_{\pm,n+1}^n}&&\shi(-Y(K),-\widehat{\Ga}_{n+1})\\
&I^\sharp(-\widehat{Y})\ar[ul]^{G_n}\ar[ur]_{G_{n+1}}&
}
\end{equation*}
and
\begin{equation*}\label{eq_4: commutative diagram, F, +}
\xymatrix@R=6ex{
\shi(-Y(K),-\widehat{\Ga}_{n})\ar[rd]_{F_n}\ar[rr]^{\psi_{\pm,n+1}^n}&&\shi(-Y(K),-\widehat{\Ga}_{n+1})\ar[ld]^{F_{n+1}}\\
&I^\sharp(-\widehat{Y})&
}
\end{equation*}
\elem

\bthm[{\cite[Proposition 4.26]{LY2020}}]\label{thm: isomorphism}
Suppose $n\in\mathbb{N}$ satisfies $q_n\geq q+2g(K)$. Then there exists an isomorphism

\[F_n^\p:\bigoplus_{i=0}^{q-1}\shi(-Y(K),-\widehat{\Ga}_n,S,\hat{i}_{max}^n-2g(K)-i)\xra{\cong} I^\sharp(-\widehat{Y}),\]
where $F_n^\p$ is the restriction of $F_n$ in Lemma \ref{lem_4: surgery triangle for knots}.
\ethm
\bdefn\label{defn: spinc decomposition}
For a fixed integer $q>0$ and any integer $s\in[0,q-1]$, suppose $[s]$ is the image of $s$ in $\mathbb{Z}_q$. Define $$I^\sharp(-\widehat{Y},[s])\deq F_n^\p(\shi(-Y(K),-\widehat{\Ga}_n,S,\hat{i}_{max}^n-2g(K)-s))\subset I^\sharp(-\widehat{Y}).$$
It is well-defined by the isomorphisms in Lemma \ref{lem_4: psi pm is isomorphism} and commutative diagrams in Lemma \ref{lem_4: surgery triangle for knots}.
\edefn
\bprop[{\cite[Corollary 1.20]{LY2021}}]\label{prop: l space dim 1}
Suppose $K$ is a knot in an integral homology sphere $Y$ and suppose $n$ is an integer. Then $-Y_{-n}(K)$ is an instanton L-space if and only if for any $[s]\in\mathbb{Z}_{|n|}$, we have $$\dim_\mathbb{C}I^\sharp(-Y_{-n}(K),[s])=1.$$
\eprop
\brem\label{rem: proof of l space dim 1}
Proposition \ref{prop: l space dim 1} also follows from the special case $(M,\ga)=(Y(1),\delta)$ in \cite[Theorem 1.1]{LY2021enhanced}:$$\en(I^\sharp(Y))=\chi(\widehat{HF}(Y))=\sum_{h\in H_1(Y)} h \in\mathbb{Z}[H_1(Y)]/\pm H_1(Y),$$where $Y$ is any rational homology sphere.
\erem

\subsection{Two spectral sequences}\quad

In this subsection, we construct spectral sequences from $\khii(-\widehat{Y},\widehat{K})$ to $I^\sharp(-\widehat{Y})$ via bypass exact triangles in Proposition \ref{prop_4: graded bypass for Ga_n-hat sutures}. 

For a fixed integer $q>0$, any fixed large integer $n$, and any integer $i$, we have the following diagram of exact triangles
\begin{equation}\label{eq: useful triangles}
\xymatrix@R=6ex{
\cdots&\ar[l]\widehat{\Ga}_{n+1}^{i,+}\ar[dr]_{\psi_{+,\mu}^{n+1}}&&\widehat{\Ga}_n^{i,+}\ar[ll]_{\psi_{+,n+1}^{n}}\ar[dr]_{\psi_{+,\mu}^n}&&\widehat{\Ga}_{n-1}^{i,+}\ar[ll]_{\psi_{+,n}^{n-1}}\ar[dr]_{\psi_{+,\mu}^{n-1}}&&\widehat{\Ga}_{n-2}^{i,+}\ar[ll]_{\psi_{+,n-1}^{n-2}}&\ar[l]\cdots
\\
&\cdots&\widehat{\Ga}_\mu^{i-q}\ar[ur]^{\psi_{+,n}^\mu}\ar[dl]^{\psi_{-,n-2}^\mu}&&\widehat{\Ga}_\mu^{i}\ar[ur]^{\psi_{+,n-1}^\mu}\ar[dl]^{\psi_{-,n-1}^\mu}&&\widehat{\Ga}_\mu^{i+q}\ar[ur]^{\psi_{+,n-2}^\mu}\ar[dl]^{\psi_{-,n}^\mu}&\cdots&\\
\cdots\ar[r]&\widehat{\Ga}_{n-2}^{i,-}\ar[rr]_{\psi_{-,n-1}^{n-2}}&&\widehat{\Ga}_{n-1}^{i,-}\ar[rr]_{\psi_{-,n}^{n-1}}\ar[ul]_{\psi_{-,\mu}^{n-1}}&&\widehat{\Ga}_n^{i,-}\ar[rr]_{\psi_{-,n+1}^n}\ar[ul]_{\psi_{-,\mu}^n}&&\widehat{\Ga}_{n+1}^{i,-}\ar[ul]_{\psi_{-,\mu}^{n+1}}\ar[r]&\cdots
}
\end{equation}
where we write\[\begin{aligned}
\widehat{\Ga}_\mu^{i}&=\shi(-Y(K),-\widehat{\Ga}_\mu,S,i)\\
\widehat{\Ga}_k^{i,+}&=\shi(-Y(K),-\widehat{\Ga}_k,S,i+\hat{i}^k_{min}-\hat{i}^n_{min}+\hat{i}^n_{max}-\hat{i}^\mu_{max})\\
\widehat{\Ga}_k^{i,-}&=\shi(-Y(K),-\widehat{\Ga}_k,S,i+\hat{i}^k_{max}-\hat{i}^n_{max}+\hat{i}^n_{min}-\hat{i}^\mu_{min})
\end{aligned}\]for any $k\in\mathbb{N}$, and we abuse notation so that the maps $\psi_{+,*}^*,\psi_{-,*}^*$ also denote the restrictions on corresponding gradings. Note that $\hat{i}^*_{max}$ and $\hat{i}^*_{min}$ are the maximal and minimal nontrivial gradings of $\shi(-Y(K),-\widehat{\Ga}_*)$ associated to $S$, respectively. By direct calculation, we have \begin{equation}\label{eq: vanishing a}
    \widehat{\Ga}_{n+k}^{i,+}\cong \widehat{\Ga}_{n+k-1}^{i,+}\text{ for }k>\frac{i-\hat{i}_{min}^\mu}{q}\aand\widehat{\Ga}_{n-k}^{i,+}=0\text{ for }-k<\frac{i-\hat{i}_{max}^\mu}{q},
\end{equation}
\begin{equation}\label{eq: vanishing b}\widehat{\Ga}_{n+k}^{i,-}\cong \widehat{\Ga}_{n+k-1}^{i,-}\text{ for }k>\frac{\hat{i}_{max}^\mu-i}{q}\aand \widehat{\Ga}_{n-k}^{i,-}=0\text{ for }-k< \frac{\hat{i}_{min}^\mu-i}{q}.\end{equation}

\bthm\label{thm: spectral sequence}
There exist two spectral sequences $\{(E_{r,+},d_{r,+})\}_{r\ge 1}$ and $\{(E_{r,-},d_{r,-})\}_{r\ge 1}$ with \[E_{1,+}=E_{1,-}=\khii(-\widehat{Y},\widehat{K})\]induced by exact triangles in (\ref{eq: useful triangles}) involving $\psi_{+,*}^*$ and $\psi_{-,*}^*$, respectively. They are independent of the choice of the integer $n$. Suppose $\{(E_{r,\pm},d_{r,\pm})\}_{r\ge 1}$ converge to $\mathcal{G}_{\pm}$, respectively. Then there are isomorphisms $$\mathcal{G}_{\pm}\cong I^\sharp(-\widehat{Y}).$$
\ethm
\bpf
The proof is based on unrolled exact couples introduced in Subsection \ref{subsec: Unrolled exact couples}.

The exact triangles about $\psi_{+,*}^*$ form an unrolled exact couple in the sense of Definition \ref{defn: unrolled exact couples}. For simplicity, we consider the direct sum of the unrolled exact couples about $i=i_0+1,\dots,i_0+q$ for some $i_0$ such that $i\in [\hat{i}^\mu_{min},\hat{i}^\mu_{max}]$. Then the first page is the same as\[\khii(-\widehat{Y},\widehat{K})=\shi(-Y(K),-\widehat{\Ga}_\mu)\]Since there are only finitely many nontrivial gradings associated to $S$, this unrolled exact couple is bounded. Proposition \ref{prop: spectral sequence} provides a spectral sequence $\{(E_{r,+},d_{r,+})\}_{r\ge 1}$ with $E_{1,+}=\khii(-\widehat{Y},\widehat{K})$. 

Since \[\hat{i}^k_{max}-\hat{i}_{min}^k=kq-q_0-1+2g(K)\aand \hat{i}^\mu_{max}-\hat{i}_{min}^\mu=q-1+2g(K),\]for any integers $i\ge \hat{i}_{min}^\mu$ and $k< n-(q-1+2g(K))/q$, we have \begin{equation}\label{eq: calculation}
    \begin{aligned}
    (i+\hat{i}^k_{min}-\hat{i}^n_{min}+\hat{i}^n_{max}-\hat{i}^\mu_{max})-\hat{i}_{max}^k&=i+(\hat{i}^k_{min}-\hat{i}_{max}^k)+(\hat{i}^n_{max}-\hat{i}^n_{min})-\hat{i}^\mu_{max}\\&=i-(kq-q_0-1+2g(K))+(nq-q_0-1+2g(K))-\hat{i}^\mu_{max}\\&=i+(n-k)q-\hat{i}^\mu_{max}\\&\ge \hat{i}^\mu_{min}+(n-k)q-\hat{i}^\mu_{max}\\&=(n-k)q-(q-1+2g(K))\\&> 0.
\end{aligned}
\end{equation}
For such $k$, we have $\widehat{\Ga}^{i,+}_k=0$. Thus, by Theorem \ref{thm: convergence}, we know that $\{(E_{r,+},d_{r,+})\}_{r\ge 1}$ converges to \[\mathcal{G}_{+}=\bigoplus_{i=i_0+1}^{i_0+q}\widehat{\Ga}_{n+l}^{i,+}\subset \shi(-Y(K),-\widehat{\Ga}_{n+l})\]for some large integer $l$. The calculation in (\ref{eq: calculation}) also indicates that $\mathcal{G}_{+}$ lives in the middle gradings of $\shi(-Y(K),-\widehat{\Ga}_{n+l})$. Hence by Lemma \ref{lem_4: q-cyclic} and Theorem \ref{thm: isomorphism}, we know that $\mathcal{G}_{+}\cong I^\sharp(-\widehat{Y})$. The independence of the integer $n$ follows from Lemma \ref{lem_4: psi pm is isomorphism} and Lemma \ref{lem: commutative diagram 2}. The maps $\psi_{-,*}^*$ induce an isomorphism between spectral sequences since they induce an isomorphism between the first pages.


A similar argument applies to exact triangles involving $\psi_{-,*}^*$ and we obtain another spectral sequence $\{(E_{r,-},d_{r,-})\}_{r\ge 1}$ with $E_{1,-}=\khii(-\widehat{Y},\widehat{K})$, which converges to $$\mathcal{G}_{-}\subset\shi(-Y(K),-\widehat{\Ga}_{n+l})$$in middle gradings for some large integer $l$. Also, we have $\mathcal{G}_{-}\cong I^\sharp(-\widehat{Y})$. 
\epf

\subsection{Bent complexes}\label{subsec: the bent complex}\quad

In this subsection, we construct the bent complex and relate its homology to negative large surgeries. The construction and the name are inspired by Heegaard Floer theory (cf.\ \cite[Section 4.1]{Rasmussen2007}, \cite[Section 2.2]{Rasmussen2017}; see also \cite[Section 4]{ozsvath2004holomorphicknot}).

\bcons\label{cons: bent complex}
Suppose $\hat{\mu}=q\mu+p\lambda$. Consider the spectral sequences $\{(E_{r,+},d_{r,+})\}_{r\ge 1}$ and $\{(E_{r,-},d_{r,-})\}_{r\ge 1}$ constructed in Theorem \ref{thm: spectral sequence}. By fixing a basis of $\khii(-\widehat{Y},\widehat{K})$, Construction \ref{cons: recover filtered chain complex} provides two filtered chain complexes \[(\khii(-\widehat{Y},\widehat{K}),d_+)\aand (\khii(-\widehat{Y},\widehat{K}),d_-)\]such that the induced spectral sequences are $\{(E_{r,+},d_{r,+})\}_{r\ge 1}$ and $\{(E_{r,-},d_{r,-})\}_{r\ge 1}$, respectively. For any integer $s$, the \textbf{bent complex} is \[A_s=A_s(-\widehat{Y},\widehat{K})\deq (\bigoplus_{k\in\mathbb{Z}}\shi(-Y(K),-\widehat{\Ga}_\mu,S,s+kq),d_s),\]where for any element $x\in \shi(-Y(K),-\widehat{\Ga}_\mu,S,s+kq)$, 
\[
d_s(x)=\begin{cases}
d_+(x)&k>0,\\
d_+(x)+d_-(x)&k=0,\\
d_-(x)&k<0.
\end{cases}\]It is easy to check that $d_s\circ d_s=0$.
\econs
\brem\label{rem: ambiguity}
Since $\shi$ is a projectively transitive system, the maps $d_{r,+}$ and $d_{r,-}$ are only well-defined up to multiplication by a unit. However, the kernel and the image of a map are still well-defined, so we can still define exact sequences for projectively transitive systems. Moreover, if $f:A\to B$ and $g:A\to C$ are maps between projectively transitive systems, though the map $$f+g\deq f\oplus g=(f,g):A\to B\oplus C$$is not well-defined, its kernel ($\ke f\cap \ke g$) is well-defined, so there is no ambiguity in considering the dimension of the homology of the bent complex. Alternatively, by discussion in Subsection \ref{subsec: Backgrounds on sutured instanton homology}, we can always fix closures of corresponding balanced sutured manifolds and consider linear maps between actual vector spaces, at the cost that equations between maps only hold up to multiplication by a unit.
\erem

The main theorem of this subsection is the following.
\bthm\label{thm: homology of bent complex}
Suppose $\hat{\mu}=q\mu+p\lambda$ with $q\in\mathbb{N}_+$. For any integer $s$, let $H(A_s)$ denote the homology of the bent complex $A_s$ in Construction \ref{cons: bent complex}. For any integer $n$ satisfying $(n-1)q\ge 2g(K)$, we have an isomorphism for some integer $j_n$:\begin{equation}\label{eq: isomorphism of bent complex}
    a_{s,n}:H(A_s)\xra{\cong} \shi(-Y(K),-\ga_{2\hat{\lambda}-(2n-1)\hat{\mu}},S,s+j_n).
\end{equation}Suppose the maximal and minimal nontrivial gradings of $\shi(-Y(K),-\ga_{2\hat{\lambda}-(2n-1)\hat{\mu}})$ are $\hat{i}^\sharp_{max}$ and $\hat{i}^\sharp_{min}$, which can be calculated by Lemma \ref{lem_4: top and bottom non-vanishing gradings}. Then we have $$j_n=\hat{i}^\sharp_{min}-\hat{i}^n_{min}+\hat{i}^n_{max}-\hat{i}^\mu_{max}=\hat{i}^\sharp_{max}-\hat{i}^n_{max}+\hat{i}^n_{min}-\hat{i}^\mu_{min}.$$
\ethm

\brem
By Definition \ref{defn_4: top and bottom non-vanishing gradings}, we have $i_{max}^y-i_{min}^y=2g(K)+y-1$. Then
$$\begin{aligned}
    (\hat{i}^\sharp_{min}-&\hat{i}^n_{min}+\hat{i}^n_{max}-\hat{i}^\mu_{max})-(\hat{i}^\sharp_{max}-\hat{i}^n_{max}+\hat{i}^n_{min}-\hat{i}^\mu_{min})\\&=2(\hat{i}^n_{max}+\hat{i}^n_{min}-)-(\hat{i}^\sharp_{max}-\hat{i}^\sharp_{min})-(\hat{i}^\mu_{max}-\hat{i}^\mu_{min})\\&=2(nq-q_0-1)-((2n-1)q-2q_0-1)-(q-1)\\&=0.
\end{aligned}$$Hence $j_n$ in Theorem \ref{thm: homology of bent complex} is well-defined. 
\erem

\bpf[Proof of Theorem \ref{thm: homology of bent complex}]
We consider two cases. The first case is special, and we use the octahedral axiom to prove it. The second case is more general, and we reduce it to the first case. For the bent complex $A_s$, we fix $i=s$ in the diagram (\ref{eq: useful triangles}).

{\bf Case 1}. Suppose $\widehat{\Ga}_k^{i,+}=\widehat{\Ga}_k^{i,-}=0$ for $k\le n-2$ in the diagram (\ref{eq: useful triangles}).

In this case, higher differentials $d_{r,\pm}$ for $r\ge 2$ vanish and the maps$$\psi_{\pm,\mu}^{n-1}:\widehat{\Ga}_{n-1}^{i,\pm}\to \widehat{\Ga}_{\mu}^{i\pm q}$$are isomorphisms. Hence $$A_s=(\widehat{\Ga}_\mu^i\oplus\widehat{\Ga}_{n-1}^{i,+}\oplus \widehat{\Ga}_{n-1}^{i,-},f),$$where $$\begin{aligned}f:\widehat{\Ga}_\mu^i&\to \widehat{\Ga}_{n-1}^{i,+}\oplus \widehat{\Ga}_{n-1}^{i,-}\\f(x)&=(\beta_+(x),\beta_-(x))\end{aligned}$$is the restriction of $(\psi_{+,n-1}^\mu(x),\psi_{-,n-1}^\mu(x))$. Define $g:\widehat{\Ga}_{n-1}^{i,+}\oplus \widehat{\Ga}_{n-1}^{i,-}\to \widehat{\Ga}_{n-1}^{i,+}$ to be the projection map. Then we apply Theorem \ref{thm: octahedral axiom} to $$X=\widehat{\Ga}_\mu^i,Y=\widehat{\Ga}_{n-1}^{i,+}\oplus \widehat{\Ga}_{n-1}^{i,-},Z=\widehat{\Ga}_{n-1}^{i,+},X^\p=\widehat{\Ga}_{n-1}^{i,-},Y^\p=\widehat{\Ga}_{n}^{i,+},Z^\p=H(A_s).$$Then there exist maps $\psi$ and $\phi$ making  the following diagram commute and exact
\begin{equation*}
    \xymatrix@R=6ex{
    &&&H(A_s)\ar[dr]^{\psi}&&&
    \\
    &&\widehat{\Ga}_{n-1}^{i,+}\oplus \widehat{\Ga}_{n-1}^{i,-}\ar[dr]^{g}\ar[ur]&&\widehat{\Ga}_{n}^{i,+}\ar[ddrr]^{\phi}&&
    \\
 &&&\widehat{\Ga}_{n-1}^{i,+}\ar[drrr]^{0}\ar[ur]&&&
    \\
    \widehat{\Ga}_\mu^i\ar[uurr]^{f}\ar[urrr]^{g\circ f=\beta_+}&&&&&&\widehat{\Ga}_{n-1}^{i,-}
    }
    \end{equation*}
Thus, we obtain a long exact sequence$$\cdots\to H(A_s)\xra{\psi} \widehat{\Ga}_{n}^{i,+}\xra{\phi} \widehat{\Ga}_{n-1}^{i,-}\to H(A_s)\{1\}\to\cdots$$Let $$\al_+:\widehat{\Ga}_{n}^{i,+}\to \widehat{\Ga}_{\mu}^i$$ be the restriction of $\psi_{+,\mu}^n$. Note that$$\widehat{\Ga}_{n}^{i,+}\cong {\rm Im}(\psi_{+,n}^{n-1}:\widehat{\Ga}_{n-1}^{i,+}\to \widehat{\Ga}_{n}^{i,+}) \oplus{\rm Coker}(\psi_{+,n}^{n-1}:\widehat{\Ga}_{n-1}^{i,+}\to \widehat{\Ga}_{n}^{i,+})\cong \ke(\beta_+)\oplus\cok(\beta_+).$$We know that the maps $\phi$ and $\phi^\p\deq \beta_-\circ \al_+$ satisfy the assumption of Lemma \ref{lem: identify phi}. Thus, we have \begin{equation}\label{eq: a}
    H(A_s)\cong H(\cone(\phi))\cong H(\cone(\beta_-\circ \al_+)).
\end{equation}
Note that we assume $\hat{\mu}=q\mu+p\lambda$ for $q\ge 0$ and $\hat{\lambda}=q_0\mu+p_0\lambda$ satisfying Definition \ref{defn_4: Ga_n-hat sutures}. When $n$ is large, the coefficient of $\mu$ in$$\hat{\mu}^\p\deq n\hat{\mu}-\hat{\lambda}=(nq-q_0)\mu+(np-p_0)\lambda$$is positive. By Definition \ref{defn_4: Ga_n-hat sutures}, we set$$\hat{\lambda}^\p\deq \hat{\lambda}-(n-1)\hat{\mu}=(q_0-(n-1)q)\mu+(p_0-(n-1)p)\lambda.$$Then $$\hat{\lambda}^\p+\hat{\mu}^\p=\hat{\mu}\aand \hat{\lambda}^\p-\hat{\mu}^\p=2\hat{\lambda}-(2n-1)\hat{\mu}.$$Note that $\ga_{x\lambda+y\mu}=\ga_{-x\lambda-y\mu}$. Applying the diagram (\ref{eq: commutative diagram 2b}) with $\psi_{+,*}^-$ and $\psi_{-,*}^+$ switched to $$\widehat{\Ga}_\mu(\hat{\mu}^\p)=\ga_{\hat{\mu}^\p}=\widehat{\Ga}_n,\widehat{\Ga}_{-1}(\hat{\mu}^\p)=\ga_{\hat{\lambda}^\p+\hat{\mu}^\p}=\widehat{\Ga}_\mu,\aand \widehat{\Ga}_0(\hat{\mu}^\p)=\ga_{\hat{\lambda}^\p}=\widehat{\Ga}_{n-1},$$we obtain the following commutative diagram\begin{equation}\label{eq: commutative diagram 2b prime}
\xymatrix@R=6ex{
\shi(-Y(K),-\widehat{\Ga}_{-1}(\hat{\mu}^\p))\ar[rr]^{\psi_{+,0}^{-1}(\hat{\mu}^\p)}&&\shi(-Y(K),-\widehat{\Ga}_0(\hat{\mu}^\p))\\
&\shi(-\widehat{Y}(K),-\widehat{\Ga}_\mu(\hat{\mu}^\p))\ar[ul]^{\psi_{-,-1}^\mu(\hat{\mu}^\p)}\ar[ur]_{\psi_{-,0}^\mu(\hat{\mu}^\p)}&
}    
\end{equation}
where the notations $\hat{\mu}^\p$ in bypass maps indicate that they correspond to $\hat{\mu}^\p$. By comparing the grading shifts, we have $$\psi_{+,0}^{-1}(\hat{\mu}^\p)=\be_-\aand \psi_{-,-1}^\mu(\hat{\mu}^\p)=\al_+.$$Indeed, this can be obtained by a diagrammatic way as shown in \cite[Remark 4.15]{LY2020}.

Let $\delta:\widehat{\Ga}_{n}^{i,+}\to\widehat{\Ga}_{n-1}^{i,-}$ be the restriction of $$\psi_{-,0}^\mu(\hat{\mu}^\p):\shi(-\widehat{Y}(K),-\widehat{\Ga}_{n})\to \shi(-Y(K),-\widehat{\Ga}_{n-1}).$$Then (\ref{eq: commutative diagram 2b prime}) implies $\delta=\beta_-\circ \al_+=\phi$.

Applying the negative bypass triangle in Theorem \ref{prop_4: graded bypass for Ga_n-hat sutures} to $$\widehat{\Ga}_\mu(\hat{\mu}^\p)=\ga_{\hat{\mu}^\p}=\widehat{\Ga}_n,\widehat{\Ga}_0(\hat{\mu}^\p)=\ga_{\hat{\lambda}^\p}=\widehat{\Ga}_{n-1},\aand \widehat{\Ga}_1(\hat{\mu}^\p)=\ga_{\hat{\lambda}^\p-\hat{\mu}^\p}=\ga_{2\hat{\lambda}-(2n-1)\hat{\mu}},$$we have the following exact triangle\begin{equation}\label{eq: triangle a}
\xymatrix@R=6ex{
\shi(-Y(K),-\widehat{\Ga}_0(\hat{\mu}^\p))\ar[rr]^{\psi_{-,1}^0(\hat{\mu}^\p)}&&\shi(-Y(K),-\widehat{\Ga}_1(\hat{\mu}^\p))\ar[dl]^{\psi_{-,\mu}^1(\hat{\mu}^\p)}\\
&\shi(-\widehat{Y}(K),-\widehat{\Ga}_\mu(\hat{\mu}^\p))\ar[ul]^{\psi_{-,0}^\mu(\hat{\mu}^\p)}&
}    
\end{equation}
By grading shifts in Theorem \ref{prop_4: graded bypass for Ga_n-hat sutures}, the restriction of (\ref{eq: triangle a}) on a single grading implies \begin{equation}\label{eq: b}
    H(\cone(\delta))\cong \shi(-\widehat{Y}(K),-\ga_{2\hat{\lambda}-(2n-1)\hat{\mu}},S,j_n)
\end{equation}
Then the isomorphism in (\ref{eq: isomorphism of bent complex}) follows from (\ref{eq: a}) and (\ref{eq: b}).

{\bf Case 2}. We do not suppose $\widehat{\Ga}_k^{i,+}=\widehat{\Ga}_k^{i,-}=0$ for all $k\le n-2$ in the diagram (\ref{eq: useful triangles}). Since $(n-1)q\ge 2g(K)$ and $i\in [\hat{i}^\mu_{min},\hat{i}^\mu_{max}]$, we have $$|\frac{i-\hat{i}^\mu_{min}}{q}|, |\frac{i-\hat{i}^\mu_{max}}{q}|\le |\frac{\hat{i}^\mu_{max}-\hat{i}^\mu_{min}}{q}|=\frac{q-1+2g(K)}{q}< n.$$By (\ref{eq: vanishing a}) and (\ref{eq: vanishing b}), we have $\widehat{\Ga}_{0}^{i,\pm}=0$.

In this case, let $$A_s^\p=(\bigoplus_{k\in\mathbb{Z}\backslash\{0\}}\shi(-Y(K),-\widehat{\Ga}_\mu,S,s+kq),d_s)$$ be the subcomplex of $A_s$. The quotient $A_s/A_s^\p$ is $\shi(-Y(K),-\widehat{\Ga}_\mu,S,s)$ with no differentials. Then we have a long exact sequence$$\cdots\to H(A^\p_s)\to H(A_s)\to H(A_s/A_s^\p)\xra{\partial_*} H(A^\p_s)\{1\}\to \cdots$$Since $\widehat{\Ga}_{0}^{i,\pm}=0$, by Theorem \ref{thm: convergence}, we know that \begin{equation}\label{eq: isomorphism a}
    H(A^\p_s)\cong \widehat{\Ga}_{n-1}^{i,+}\oplus \widehat{\Ga}_{n-1}^{i,-}.
\end{equation}It is straightforward to check $\partial_*=(\beta_+,\beta_-)$ under the isomorphism (\ref{eq: isomorphism a}). Then by Case 1, we have$$H(A_s)\cong H(\cone(\partial _*))\cong H(\cone(f))\cong H(\cone(\phi))\cong \shi(-\widehat{Y}(K),-\ga_{2\hat{\lambda}-(2n-1)\hat{\mu}},S,j_n).$$

\epf

Then we prove the large surgery formula for negative surgeries.
\bthm[Theorem \ref{thm: large surgery formula, main}, $n>0$]\label{thm: large surgery formula, bent}
Suppose $\widehat{\mu}=q\mu+p\lambda$ with $q\in\mathbb{N}_+$ and suppose $\widehat{\lambda}=q_0\mu+p_0\lambda$ is defined as in Definition \ref{defn_4: surgery on pairs}. Note that when $(q,p)=(1,0)$, we have $(q_0,p_0)=(0,1)$. For a fixed integer $n$ satisfying $(n-1)q\ge 2g(K)$, suppose$$\hat{\mu}^\p=n\hat{\mu}-\hat{\lambda}=(nq-q_0)\mu+(np-p_0)\lambda.$$For any integer $s^\p$, suppose $[s^\p]$ is the image of $s^\p$ in $\mathbb{Z}_{(nq-q_0)}$. Suppose $$s_{min}=-(nq-q_0-1)-\lceil -\frac{q-1}{2}\rceil+g(K)\aand s_{max}=(nq-q_0-1)-\lceil \frac{q-1}{2}\rceil-g(K)$$and suppose we have an integer $s\in[s_{min},s_{max}].$
For such $n$ and $s$, there is an isomorphism 
\[H(A_{-s})\cong I^\sharp(-\widehat{Y}_{\hat{\mu}^\p},[s-s_{min}]).\]
\ethm
\brem
When $(n-1)q\ge 2g(K)$, there are more than $(nq-q_0)$ integers in the interval $[s_{min},s_{max}]$. Thus, the bent complexes contain all information of $I^\sharp(-\widehat{Y}_{\hat{\mu}^\p})$.
\erem
\bpf[Proof of Theorem \ref{thm: large surgery formula, bent}]
Since $(n-1)q\ge 2g(K)$, we apply Theorem \ref{thm: homology of bent complex} to obtain$$H(A_{-s})\cong  \shi(-Y(K),-\ga_{2\hat{\lambda}-(2n-1)\hat{\mu}},S,j_n-s).$$
We adapt the notations $$\hat{\lambda}^\p=\hat{\lambda}-(n-1)\hat{\mu}\aand \hat{\lambda}^\p-\hat{\mu}^\p=2\hat{\lambda}-(2n-1)\hat{\mu}=(2q_0-(2n-1)q)\mu+(2p_0-(2n-1)p)\lambda$$from the proof of Theorem \ref{thm: homology of bent complex}. Then $\widehat{\Ga}_1(\hat{\mu}^\p)=\ga_{2\hat{\lambda}-(2n-1)\hat{\mu}}$. Since $(n-1)q\ge 2g(K)$, we have$$(2n-1)q-2q_0\ge nq-q_0 + 2g(K).$$Hence, we can apply Theorem \ref{thm: isomorphism} to obtain $$I^\sharp(-\widehat{Y}_{\hat{\mu}^\p},[s])\cong \shi(-Y(K),-\ga_{2\hat{\lambda}-(2n-1)\hat{\mu}},S,\hat{i}_{max}^\sharp-2g(K)-s).$$By direct calculation, we have$$\begin{aligned}j_n-s_{min}=&\hat{i}^\sharp_{max}-\hat{i}^n_{max}+\hat{i}^n_{min}-\hat{i}^\mu_{min}-s_{min}\\=&\hat{i}_{max}^\sharp-2g(K)-(nq-q_0-1)-\lceil -\frac{q-1}{2}\rceil+g(K)-s_{min}
\\=&\hat{i}_{max}^\sharp-2g(K).
\end{aligned}$$For any $s\in[s_{min},s_{max}]$, we have 

 $$\begin{aligned}j_n-s=&\hat{i}^\sharp_{min}-\hat{i}^n_{min}+\hat{i}^n_{max}-\hat{i}^\mu_{max}-s\\=&\hat{i}_{min}^\sharp+2g(K)+(nq-q_0-1)-\lceil \frac{q-1}{2}\rceil-g(K)-s\
 \\\ge&\hat{i}_{min}^\sharp+2g(K).
 \end{aligned}$$

Thus, the isomorphism follows from Definition \ref{defn: spinc decomposition} and Lemma \ref{lem_4: q-cyclic}.
\epf

Finally, we state an instanton analog of \cite[Theorem 2.3]{Ozsvath2008integral} and \cite[Theorem 4.1]{Ozsvath2011rational}, which is an important step of the proof of the mapping cone formula (cf.\ Section \ref{sec: future}).

\bcons\label{cons: bent complex, b}
Following the notations in Construction \ref{cons: bent complex}.
 For $\circ\in\{+,-\}$, define \[B_s^\circ=B_s^\circ(-\widehat{Y},\widehat{K})\deq(\bigoplus_{k\in\mathbb{Z}}\shi(-Y(K),-\widehat{\Ga}_\mu,S,s+kq),d_\circ)\]and define\[\pi_{s}^\circ:A_s\to B_s^\circ\]by
\[\pi_s^+(x)=\begin{cases}
x&k>0,\\
0&k\le 0,\end{cases}\aand \pi_s^-(x)=\begin{cases}
0&k\ge 0,\\
0&k<0,\end{cases}\]where $x\in \shi(-Y(K),-\widehat{\Ga}_\mu,S,s+kq)$.
\econs
 Suppose $\hat{\mu}=q\mu+p\lambda$ with $q\in\mathbb{N}_+$. For $n$ and $s$ in Theorem \ref{thm: homology of bent complex}, let $H(A_s),H(B_s^+),H(B_s^-)$ be the homologies of complexes in Construction \ref{cons: bent complex} and let $(\pi_s^+)_*,(\pi_s^-)_*$ denote the induced maps on homologies. Let $j_n$ be the integer in Theorem \ref{thm: homology of bent complex} and write $\widehat{\Ga}^{s,\sharp}$ for $$\shi(-Y(K),-\ga_{2\hat{\lambda}-(2n-1)\hat{\mu}},S,j_n+s).$$By Theorem \ref{thm: homology of bent complex}, we have an isomorphism $$a_{s,n}:H(A_s)\xra{\cong} \widehat{\Ga}^{s,\sharp}$$ 

We use notations in (\ref{eq: useful triangles}) and set $i=s$. Let$$\rho_+:\widehat{\Ga}^{s,\sharp}\to \widehat{\Ga}_n^{s,+}$$be the restriction of $\psi_{-,\mu}^1(\hat{\mu}^\p)$ in the proof of Theorem \ref{thm: homology of bent complex}. Choose $l$ as in the proof of Theorem \ref{thm: spectral sequence} such that $\widehat{\Ga}_{n+l}^{s,+}\subset G_+$. Note that $H(B_s^\pm)=\widehat{\Ga}_{n+l}^{s,\pm}$ by the proof of Theorem \ref{thm: spectral sequence}. Let $$\Psi_{+,n+l}^n:\widehat{\Ga}_n^{s,+}\to\widehat{\Ga}_{n+l}^{s,+}$$ be the composition of $\psi_{+,n+k+1}^{n+k}$ for $k=0,\dots,l-1$. Similarly, let $$\rho_-:\widehat{\Ga}^{s,\sharp}\to \widehat{\Ga}_n^{s,-}$$be the restriction of $\psi_{+,\mu}^1(\hat{\mu}^\p)$ and let $$\Psi_{-,n+l}^n:\widehat{\Ga}_n^{s,-}\to\widehat{\Ga}_{n+l}^{s,-}\subset G_-$$ be the composition of $\psi_{-,n+k+1}^{n+k}$ for $k=0,\dots,l-1$. 
\bprop\label{prop: bent commute}
The following diagram commutes.
\begin{equation*}
    \xymatrix@R=6ex{
    H(A_s)\ar[rr]^{(\pi_s^\pm)_*}\ar[d]^{a_{s,n}}&& H(B_s^\pm)\ar[d]^{=}\\
    \widehat{\Ga}^{s,\sharp}\ar[rr]^{\Psi_{\pm,n+l}^n\circ\rho_\pm}&& \widehat{\Ga}_{n+l}^{s,\pm},
    }
\end{equation*}
\eprop
\bpf
The proof is straightforward by the proof of Theorem \ref{thm: homology of bent complex}.
\epf

\brem
By direct calculation, the difference of gradings of $\widehat{\Ga}_{n+l}^{s,+}$ and $\widehat{\Ga}_{n+l}^{s,-}$ is$$\begin{aligned}
(\hat{i}^{n+l}_{min}&-\hat{i}^n_{min}+\hat{i}^n_{max}-\hat{i}^\mu_{max})-(\hat{i}^{n+l}_{max}-\hat{i}^n_{max}+\hat{i}^n_{min}-\hat{i}^\mu_{min})\\&=-(\hat{i}^{n+l}_{max}-\hat{i}^{n+l}_{min})+2(\hat{i}^n_{max}-\hat{i}^n_{min})-(\hat{i}^\mu_{max}-\hat{i}^\mu_{min})\\&=-(n+l)q+q_0+2(nq-q_0)-q\\&=(n-l-1)q-q_0.
\end{aligned}$$
By Lemma \ref{lem_4: q-cyclic}, the spaces $\widehat{\Ga}_{n+l}^{s,+}$ and $\widehat{\Ga}_{n+l}^{s,-}$ correspond to $I^\sharp(-\widehat{Y},[s_0-q_0])$ and $I^\sharp(-\widehat{Y},[s_0])$ for some integer $s_0$, respectively. Note that the core knot corresponding to $\hat{\mu}=q\mu+p\lambda$ is isotopic to the curve $q_0\mu+p_0\lambda$ on $\partial Y(K)$.
\erem

\subsection{Dual bent complexes}\quad

In this subsection, we construct the dual bent complex and relate its homology to large positive surgeries. Proofs are similar to those in Subsection \ref{subsec: the bent complex}, so we only point out the difference.

\bcons\label{cons: bent complex, dual}
Following the notations in Construction \ref{cons: bent complex}. For any integer $s$, define the \textbf{dual bent complex} as \[A_s^\vee=A^\vee_s(-\widehat{Y},\widehat{K})\deq(\bigoplus_{k\in\mathbb{Z}}\shi(-Y(K),-\widehat{\Ga}_\mu,S,s+kq),d_s^\vee),\]where for any element $x\in \shi(-Y(K),-\widehat{\Ga}_\mu,S,s+kq)$, 
\[
d_s^\vee(x)=\begin{cases}
d_-(x)&k>0,\\
d_+(x)+d_-(x)&k=0,\\
d_+(x)&k<0.
\end{cases}\]
\econs

Similar to Theorem \ref{thm: homology of bent complex} and Theorem \ref{thm: large surgery formula, bent}, we have the following theorems.

\bthm\label{thm: the homology of dual bent complex}
Suppose $\hat{\mu}=q\mu+p\lambda$ with $q\in\mathbb{N}_+$. For any integer $s$, let $H(A_s^\vee)$ denote the homology of the bent complex $A_s^\vee$ in Construction \ref{cons: bent complex, dual}. For any integer $n$ satisfying $(n-1)q\ge 2g(K)$, we have an isomorphism for some integer $j_n^\vee$:\begin{equation}\label{eq: isomorphism of bent complex, dual}
    a_{s,n}^\vee:H(A_s^\vee)\xra{\cong} \shi(-Y(K),-\ga_{2\hat{\lambda}+(2n+1)\hat{\mu}},S,s+j_n^\vee).
\end{equation}Suppose the maximal and minimal nontrivial gradings of $\shi(-Y(K),-\ga_{2\hat{\lambda}+(2n+1)\hat{\mu}})$ are $\hat{i}^{\sharp,\vee}_{max}$ and $\hat{i}^{\sharp,\vee}_{min}$, which can be calculated by Lemma \ref{lem_4: top and bottom non-vanishing gradings}. Then we have $$j_n^\vee=\hat{i}^{\sharp,\vee}_{max}-\hat{i}^{-n}_{max}+\hat{i}^{-n}_{min}-\hat{i}^\mu_{min}=\hat{i}^{\sharp,\vee}_{min}-\hat{i}^{-n}_{min}+\hat{i}^{-n}_{max}-\hat{i}^\mu_{max}.$$
\ethm

\bthm[Theorem \ref{thm: large surgery formula, main}, $n<0$]\label{thm: large surgery formula, dual bent}
Suppose $\widehat{\mu}=q\mu+p\lambda$ with $q\in\mathbb{N}_+$ and suppose $\widehat{\lambda}=q_0\mu+p_0\lambda$ is defined as in Definition \ref{defn_4: surgery on pairs}. Note that when $(q,p)=(1,0)$, we have $(q_0,p_0)=(0,1)$. For a fixed integer $n$ satisfying $(n-1)q\ge 2g(K)$, suppose$$\hat{\mu}^\pp=n\hat{\mu}+\hat{\lambda}=(nq+q_0)\mu+(np+p_0)\lambda.$$For any integer $s^\p$, suppose $[s^\p]$ is the image of $s^\p$ in $\mathbb{Z}_{(nq+q_0)}$. Suppose $$s_{min}^\vee=-(nq+q_0-1)-\lceil -\frac{q-1}{2}\rceil+g(K)\aand s_{max}^\vee=(nq+q_0-1)-\lceil \frac{q-1}{2}\rceil-g(K)$$and suppose an integer $s\in[s_{min}^\vee,s_{max}^\vee].$
For such $n$ and $s$, there is an isomorphism 
\[H(A_{-s}^\vee)\cong I^\sharp(-\widehat{Y}_{\hat{\mu}^\pp},[s-s_{min}^\vee]).\]
\ethm

\bpf[Proof of Theorem \ref{thm: the homology of dual bent complex}]
Instead of using the diagram \ref{eq: useful triangles}, we use the following diagram of exact triangles from Proposition \ref{prop_4: graded bypass for Ga_n-hat sutures}:
\begin{equation}\label{eq: useful triangles, dual}
\xymatrix@R=6ex{
\cdots&\ar[l]\widehat{\Ga}_{-n+1}^{i,+}\ar[dr]_{\psi_{+,\mu}^{-n+1}}&&\widehat{\Ga}_{-n}^{i,+}\ar[ll]_{\psi_{+,-n+1}^{-n}}\ar[dr]_{\psi_{+,\mu}^{-n}}&&\widehat{\Ga}_{-n-1}^{i,+}\ar[ll]_{\psi_{+,-n}^{-n-1}}\ar[dr]_{\psi_{+,\mu}^{-n-1}}&&\widehat{\Ga}_{-n-2}^{i,+}\ar[ll]_{\psi_{+,-n-1}^{-n-2}}&\ar[l]\cdots
\\
&\cdots&\widehat{\Ga}_\mu^{i-q}\ar[ur]^{\psi_{+,-n}^\mu}\ar[dl]^{\psi_{-,-n-2}^\mu}&&\widehat{\Ga}_\mu^{i}\ar[ur]^{\psi_{+,-n-1}^\mu}\ar[dl]^{\psi_{-,-n-1}^\mu}&&\widehat{\Ga}_\mu^{i+q}\ar[ur]^{\psi_{+,-n-2}^\mu}\ar[dl]^{\psi_{-,-n}^\mu}&\cdots&\\
\cdots\ar[r]&\widehat{\Ga}_{-n-2}^{i,-}\ar[rr]_{\psi_{-,-n-1}^{-n-2}}&&\widehat{\Ga}_{-n-1}^{i,-}\ar[rr]_{\psi_{-,-n}^{-n-1}}\ar[ul]_{\psi_{-,\mu}^{-n-1}}&&\widehat{\Ga}_{-n}^{i,-}\ar[rr]_{\psi_{-,-n+1}^{-n}}\ar[ul]_{\psi_{-,\mu}^{-n}}&&\widehat{\Ga}_{-n+1}^{i,-}\ar[ul]_{\psi_{-,\mu}^{-n+1}}\ar[r]&\cdots
}
\end{equation}
where we write\[\begin{aligned}
\widehat{\Ga}_\mu^{i}&=\shi(-Y(K),-\widehat{\Ga}_\mu,S,i)\\
\widehat{\Ga}_{-k}^{i,+}&=\shi(-Y(K),-\widehat{\Ga}_{-k},S,i+\hat{i}^{-k}_{max}-\hat{i}^{-n}_{max}+\hat{i}^{-n}_{min}-\hat{i}^\mu_{min})\\
\widehat{\Ga}_{-k}^{i,-}&=\shi(-Y(K),-\widehat{\Ga}_{-k},S,i+\hat{i}^{-k}_{min}-\hat{i}^{-n}_{min}+\hat{i}^{-n}_{max}-\hat{i}^\mu_{max})
\end{aligned}\]for any $k\in\mathbb{N}$, and we abuse notation so that the maps $\psi_{+,*}^*,\psi_{-,*}^*$ also denote the restrictions on corresponding gradings. In this case, we have 
\begin{equation}\label{eq: vanishing c}
     \widehat{\Ga}_{-n-k}^{i,+}\cong \widehat{\Ga}_{-n-k-1}^{i,+}\text{ for }k>\frac{\hat{i}_{max}^\mu-i}{q}\aand\widehat{\Ga}_{-n+k}^{i,+}=0\text{ for }-k<\frac{\hat{i}_{min}^\mu-i}{q},
 \end{equation}
 \begin{equation}\label{eq: vanishing d}\widehat{\Ga}_{-n-k}^{i,-}\cong \widehat{\Ga}_{-n-k-1}^{i,-}\text{ for }k>\frac{i-\hat{i}_{min}^\mu}{q}\aand \widehat{\Ga}_{-n+k}^{i,-}=0\text{ for }-k< \frac{i-\hat{i}_{max}^\mu}{q}.\end{equation}
By Proposition \ref{prop: spectral sequence} and Theorem \ref{thm: convergence}, there exist spectral sequences from $$\bigoplus_{k\in\mathbb{Z}}\widehat{\Ga}_\mu^{i+kq}$$to $\widehat{\Ga}_{-n-l}^{i,+}$ and $\widehat{\Ga}_{-n-l}^{i,-}$ for some large $l$. By Lemma \ref{lem: commutative diagram 2}, those spectral sequences are isomorphic to $\{(E_{r,+},d_{r,+})\}_{r\ge 1}$ and $\{(E_{r,-},d_{r,-})\}_{r\ge 1}$ in Theorem \ref{thm: spectral sequence}, and hence we can define the dual bent complex by maps in (\ref{eq: useful triangles, dual}). 

By Definition \ref{defn_4: Ga_n-hat sutures}, we set $$\hat{\mu}^\pp=n\hat{\mu}+\hat{\lambda}\aand \hat{\lambda}^\pp=-\hat{\mu}.$$Then $$\hat{\lambda}^\pp-\hat{\mu}^\pp=-\hat{\lambda}-(n+1)\hat{\mu}\aand \hat{\lambda}^\pp-2\hat{\mu}^\pp=-2\hat{\lambda}-(2n+1)\hat{\mu}.$$Note that $\ga_{x\lambda+y\mu}=\ga_{-x\lambda-y\mu}$. 

Similar to the proof of Theorem \ref{thm: homology of bent complex}, we consider two cases and finally obtain that $$\begin{aligned}
    H(A_i^\vee)\cong & H(\cone(\psi_{+,\mu}^{-n}+\psi_{-,\mu}^{-n}: \widehat{\Ga}_{-n}^{i,+}\oplus \widehat{\Ga}_{-n}^{i,-}\to \widehat{\Ga}_\mu^i)\\
    \cong & H( \cone(\psi_{-,-n-1}^{\mu}\circ \psi_{+,\mu}^{-n}:\widehat{\Ga}_{-n}^{i,+}\to \widehat{\Ga}_{-n-1}^{i,-}))\\
    \cong & \shi(-Y(K),-\ga_{2\hat{\lambda}+(2n+1)\hat{\mu}},S,i+j_n^\vee).
\end{aligned}$$
\epf
\bpf[Proof of Theorem \ref{thm: large surgery formula, dual bent}]
Similar to the proof of Theorem \ref{thm: large surgery formula, bent}, the isomorphism follows from Theorem \ref{thm: isomorphism}, Definition \ref{defn: spinc decomposition}, and Lemma \ref{lem_4: q-cyclic}.
\epf

The following proposition explains the name of the `dual bent complex'.

\bprop\label{prop: dual complex}
$A^\vee_s(-Y,K)$ is the dual complex of $A_{-s}(Y,K)$. 
\eprop
\bpf
Suppose $(\bar{Y},\widebar{K})=(-Y,K)$ is the mirror of $(Y,K)$. Note that $(-\bar{Y},\widebar{K})=(Y,K)$. Suppose $S$ is the Seifert surface of $K$. Then $-S$ is the induced Seifert surface of $\widebar{K}$. By Proposition \ref{prop: mirror}, we have canonical isomorphisms

$$\begin{aligned}
    \shi(-\bar{Y}(\widebar{K}),-\widehat{\Ga}_{n},-S,i)=&\shi(Y(K),-\widehat{\Ga}_{-n},-S,i)\\\cong&\shi(Y(K),-\widehat{\Ga}_{-n},S,-i)\\\cong&{\rm Hom}_\mathbb{C}(\shi(-Y(K),-\widehat{\Ga}_{-n},S,-i),\mathbb{C})
\end{aligned}.$$

Then this proposition follows from the fact that both diagram (\ref{eq: useful triangles}) and diagram (\ref{eq: useful triangles, dual}) can be used to define the bent complex and the dual bent complex.
\epf

\subsection{Grading shifts of differentials}\label{subsec: Properties of differentials}\quad

In this subsection, we study the grading shifts of differentials $d_+$ and $d_-$ and relate the bent complex to the dual bent complex. First, it is straightforward to check from the construction that the map $d_+$ increases the $\mathbb{Z}$-grading and $d_-$ decreases the $\mathbb{Z}$-grading. So we focus on the grading shifts of $d_+$ and $d_-$ on the relative $\mathbb{Z}_2$-grading.

\begin{conv}
Throughout this subsection, `grading' means the relative $\mathbb{Z}_2$-grading and we set $M=Y(K)$ for a rationally null-homologous knot $K\subset Y$. The bypass map $\psi_{+,*}^*$ and the corresponding negative one $\psi_{-,*}^*$ are from $\shi(-M,-\ga_1)$ to $\shi(-M,-\ga_2)$ for some $\ga_1$ and $\ga_2$ consisting of two parallel simple closed curves.
\end{conv}

Since all bypass maps are homogeneous (they are constructed by cobordism maps, cf.\ the proof of \cite[Theorem 1.21]{baldwin2018khovanov}), the differentials $d_+$ and $d_-$ are also homogeneous. To study the grading shifts of $d_+$ and $d_-$, we first study the isomorphism \begin{equation}\label{eq: involution 2}
    \iota_{\ga}:\shi(-M,-\ga)\xra{\cong} \shi(-M,\ga)\xra{=}\shi(-M,-\ga)
\end{equation}
defined in (\ref{eq: involution}) more carefully.

By the construction of $\shi(-M,-\ga)$ in \cite{kronheimer2010knots,baldwin2015naturality}, we can construct a closure $(Y^\p,R,\omega)$ of $(-M,-\ga)$ with $g(R)\ge 2$ and take the $(2,2g(R)-2)$-eigenspace of $(\mu({\rm pt}),\mu(R))$ on $I^\omega(Y^\p)$. It is straightforward to check that $(Y^\p,-R,\omega)$ is a closure of $(-M,\ga)$. Hence, we can define $\shi(-M,\ga)$ by the $(2,2g(R)-2)$-eigenspace of $(\mu({\rm pt}),\mu(-R))$ on $I^\omega(Y^\p)$, which is the same as the $(2,2-2g(R))$-eigenspace of $(\mu({\rm pt}),\mu(R))$ on $I^\omega(Y^\p)$. Note that $I^\omega(Y^\p)$ has a $\mathbb{Z}_8$-grading and $\mu({\rm pt})$ and $\mu(R)$ have degree $-4$ and $-2$, respectively. The canonical isomorphism $\shi(-M,-\ga)\cong\shi(-M,\ga)$ in (\ref{eq: involution 2}) comes from the map sending $$(v_0,v_1,v_2,v_3,v_4,v_5,v_6,v_7)\in I^\omega(Y^\p)$$to $$(v_0,v_1,-v_2,-v_3,v_4,v_5,-v_6,-v_7),$$which preserves the $\mathbb{Z}_2$-grading induced by the $\mathbb{Z}_8$-grading.

Since $\ga$ and $-\ga$ are isotopic on $\partial M\cong T^2$, there is an identification $\shi(-M,-\ga)=\shi(-M,\ga)$. However, this identification may depend on the isotopy since there may be some basepoint moving map similar to Heegaard Floer theory \cite{sarkar15moving,zemke17moving}. Since we do not care about the precise identification, we omit discussion about specifying the isotopy.

\blem\label{lem: involution commute}
Suppose $\psi_{+,*}^*$ and $\psi_{-,*}^*$ are two bypass maps from $\shi(-M,-\ga_1)$ to $\shi(-M,-\ga_2)$ and suppose $\iota_{\ga_1}$ and $\iota_{\ga_2}$ are isomorphisms defined in (\ref{eq: involution}). Under some choices of isotopies of sutures, we have $$\psi_{\pm,*}^*\circ \iota_{\ga_1}=\iota_{\ga_2}\circ \psi_{\mp,*}^*.$$
\elem
\bpf
By the construction in \cite[Section 4.2]{LY2020}, the bypass arc related to $\psi_{+,*}^*$ on $(Y(K),\ga_{x\lambda+y\mu})$ is the same as the bypass arc related to $\psi_{-,*}^*$ on $(Y(K),-\ga_{x\lambda+y\mu})$. The lemma follows from the construction of the isomorphism $\iota_{\ga}$.
\epf
\bcor\label{cor: involution}
The isomorphism $\iota_{\ga}$ induces an isomorphism between spectral sequences \[\{(E_{r,+},d_{r,+})\}_{r\ge 1}~{\rm and}~\{(E_{r,-},d_{r,-})\}_{r\ge 1}\]constructed in Theorem \ref{thm: spectral sequence} and hence induces an isomorphism between the chain complexes$$(\khii(-Y,K),d_+)\aand (\khii(-Y,K),d_-).$$Moreover, it induces a canonical identification between $A_{-s}$ and $A_{s}$.
\ecor

\blem[]\label{lem: same grading}
Suppose $\psi_{+,*}^*$ and $\psi_{-,*}^*$ are two bypass maps from $\shi(-M,-\ga_1)$ to $\shi(-M,-\ga_2)$. If $x$ is a homogeneous element in $\shi(-M,-\ga_1)$, then $\psi_{+,*}^*(x)$ and $\psi_{-,*}^*(x)$ are homogeneous elements in $\shi(-M,-\ga_2)$ and they have the same grading.
\elem
\bpf
It follows from Lemma \ref{lem: involution commute} and the fact that the isomorphism $\iota_{\ga}$ preserves the grading for any $\ga\subset \partial M$. 
\epf
\bprop[]\label{prop: z2 shift}
Suppose $d_+$ and $d_-$ are differentials on $\khii(-Y,K)$ induced by spectral sequences $\{(E_{r,+},d_{r,+})\}_{r\ge 1}$ and $\{(E_{r,-},d_{r,-})\}_{r\ge 1}$ in Theorem \ref{thm: spectral sequence}. For any homogeneous element $x\in \khii(-Y,K)$, the gradings of $d_+(x)$ and $d_-(x)$ are different from the grading of $x$.
\eprop
\bpf
We prove the statement only for $d_+(x)$. The proof for $d_-(x)$ is similar. We adopt the notations in diagram (\ref{eq: useful triangles}). Without loss of generality, suppose $x\in \widehat{\Ga}_\mu^i$. Consider the projection $y$ of $d_+(x)$ on $\widehat{\Ga}_\mu^{i+kq}$ for some $k\in\mathbb{N}_+$. By the construction of $d_+$, there exist homogeneous elements $z\in \widehat{\Ga}_{n-1}^{i,+}$ and $w\in\widehat{\Ga}_{n-k}^{i,+}$ such that$$y=\psi_{+,\mu}^{n-k}(w)\aand z=\psi_{+,n-1}^\mu(x)=\psi_{+,n-1}^{n-2}\circ\cdots\circ\psi_{+,n-k+1}^{n-k}(w).$$By Lemma \ref{lem: same grading}, the element $$z^\p\deq \psi_{-,n-1}^{n-2}\circ\cdots\circ\psi_{-,n-k+1}^{n-k}(w)$$has the same grading as $z$. By Lemma \ref{lem: commutative diagram 2}, we have $$\psi_{+,\mu}^{n-1}(z^\p)=y.$$Define $$u\deq\psi_{+,n}^{n-1}(z^\p)\aand u^\p\deq\psi_{-,n}^{n-1}(z^\p).$$By Lemma \ref{lem: same grading}, they have the same grading. By \ref{lem: commutative diagram 2}, we have $$\psi_{+,\mu}^{n}(u^\p)=y.$$Let ${\rm gr}_2(x)$ denote the grading of $x$ and let ${\rm gr}_2(\psi_{+,*}^*)$ denote the grading shift of $\psi_{+,*}^*$. Then we have$$\begin{aligned}
    {\rm gr}_2(y)-{\rm gr}_2(x)&=({\rm gr}_2(y)-{\rm gr}_2(u^\p))+({\rm gr}_2(u)-{\rm gr}_2(z^\p))+({\rm gr}_2(z)-{\rm gr}_2(x))\\&={\rm gr}_2(\psi_{+,\mu}^{n})+{\rm gr}_2(\psi_{+,n}^{n-1})+{\rm gr}_2(\psi_{+,n-1}^{\mu})\\&=1,
\end{aligned}
    $$where the last equation follows from the fact that the bypass exact triangle shifts the grading (the bypass exact triangle comes from the surgery exact triangle, cf.\ the proof of \cite[Theorem 1.21]{baldwin2018khovanov}). Since any projection of $d_+(x)$ has different grading from $x$, we know that $d_+(x)$ has different grading from $x$.
\epf

\section{Vanishing results about contact elements}\label{sec: Vanishing results for contact elements}
In this section, we study contact elements in Heegaard Floer theory and instanton theory. In particular, we prove Theorem \ref{thm: contact support}, Theorem \ref{thm: SHI giroux torsion}, and a vanishing result for cobordism maps. We only need Corollary \ref{cor: important vanishing} in the remaining sections. 

\subsection{Contact elements in Heegaard Floer theory}\label{subsec: Vanishing results about Giroux torsion}\quad

In this subsection, we review the strategy to prove the vanishing result about Giroux torsion in \cite{girouxtorsion1}. 

Suppose $(N,\xi)$ is a contact 3-manifold with convex boundary and dividing set $\Ga$ on $\partial N$. Honda, Kazez, and Mati\'{c} \cite{honda2009contact} defined an element $c(N,\Ga,\xi)$ in sutured Floer homology $SFH(-N,-\Ga)$, called the \textbf{contact element} of $(N,\xi)$. When $(N,\xi)$ is obtained from a closed contact 3-manifold $(Y,\xi^\p)$ by removing a 3-ball, the element $$c(N,\Ga,\xi)\in SFH(-N,-\Ga)\cong \widehat{HF}(-Y)$$ recovers the contact element $c(Y,\xi^\p)\in\widehat{HF}(-Y)$ defined by Ozsv\'{a}th and Szab\'{o} \cite{OS05contact}.

Consider the Giroux torsion defined in Definition \ref{defn: giroux torsion}. We have the following vanishing result.

\bthm[{\cite[Theorem 1]{girouxtorsion1}}]\label{thm: SFH giroux torsion}
If a closed contact 3-manifold $(Y, \xi)$ has Giroux torsion, then its contact element $c(Y, \xi) \in\widehat{HF}(-Y)$ vanishes.
\ethm
\brem\label{rem: Z coeff}
The statement of Theorem \ref{thm: SFH giroux torsion} in \cite{girouxtorsion1} is about $\mathbb{Z}$ coefficients. However, since the naturality of $SFH$ is only proved for $\ft$ coefficients \cite{Juhasz2012}, the contact element in $\mathbb{Z}$ coefficient is not well-defined. Some progress on the naturality for $\mathbb{Z}$ coefficients is made in \cite{2019naturality}.
\erem
\brem
There are many partial results and applications of Theorem \ref{thm: SFH giroux torsion}. See the introduction of \cite{girouxtorsion1}.
\erem
Following the notations in \cite[Section 5.2]{honda2000classification}, consider a
basic slice $N_0 =(T^2 \times I, \bar{\xi})$ with the dividing set $\Ga_*$ on $T^2\times \{i\}$ for $i=0,1$ consisting of two parallel curves of slopes $s_0=\infty$ and $s_1=0$. There are two possible choices of tight structures on $N_0$ corresponding to two bypasses $\psi_{+,0}^\mu$ and $\psi_{-,0}^\mu$. They are both positively co-oriented but have different orientations. Hence, the relative Euler classes differ by signs. Let $\bar{\xi}$ be the tight structure on $N_0$ corresponding to $\psi_{+,0}^\mu$. Let $N_{\frac{n\pi}{2}}$ be obtained from $N_0$ by rotating counterclockwise by $\frac{n\pi}{2}$. Note that $N_{\pi}$ is the basic slice corresponding to $\psi_{-,0}^\mu$ and $N_{\frac{n\pi}{2}+2\pi}=N_{\frac{n\pi}{2}}$. Define $$(N_*,\zeta_1^+)=N_0\cup N_{\frac{\pi}{2}}\cup N_{\pi}\cup N_{\frac{3\pi}{2}}\cup N_{2\pi}\aand (N_*,\zeta_1^-)=N_\pi\cup N_{\frac{3\pi}{2}}\cup N_{2\pi}\cup N_{\frac{5\pi}{2}}\cup N_{3\pi}.$$
Then Theorem \ref{thm: SFH giroux torsion} follows from the following three lemmas.

\blem[{\cite[Lemma 5]{girouxtorsion1}}]\label{lem: giroux torsion}
A contact closed 3-manifold $(Y, \xi)$ has Giroux torsion if and only if there exists an embedding of $(N_*,\Ga_*,\zeta_1^+)$ or $(N_*,\Ga_*,\zeta_1^+)$ into $(Y, \xi)$.
\elem
\brem\label{rem: change slope}
In the definition of Giroux torsion, there is no condition on the orientation of the contact structure. By construction, the contact structures $\zeta_1^+$ and $\zeta_1^-$ differ in orientations. In \cite{girouxtorsion1}, the authors did not deal with these two contact structures separately (cf.\ the definition of $\zeta_0$ in \cite{girouxtorsion1}) since the proofs are almost identical. Also, in the original statement of \cite[Lemma 5]{girouxtorsion1}, the slopes of the dividing set on $\partial N_*$ are $-1$ and $-2$, respectively. However, there is a diffeomorphism of $T^2\times I$ sending the slopes to $\infty$ and $0$, respectively. Note that under this diffeomorphism, the slope $\infty$ is sent to $-1$.
\erem
\blem[{\cite[Theorem 4.5]{honda2009contact}}]\label{lem: vanishing piece SFH}
Let $(Y, \xi)$ be a closed contact 3-manifold and $N \subset Y$ be a compact submanifold (without any closed components) with convex boundary and dividing set $\Ga$. If $c(N, \Ga, \xi|_N) = 0$, then $c(Y, \xi) = 0$.
\elem
\blem[{From the proof of \cite[Theorem 1]{girouxtorsion1}}]\label{lem: vanishing result SFH}
The elements $c(N_*,\Ga_*,\zeta_1^+)$ and $c(N_*,\Ga_*,\zeta_1^-)$ vanish.
\elem

\subsection{Construction of instanton contact elements}\label{subsec: The instanton contact element}\quad

In \cite{baldwin2016instanton}, Baldwin and Sivek constructed a contact invariant in sutured instanton theory which we call the \textbf{instanton contact element}. In this subsection, we review the construction and prove Theorem \ref{thm: contact support}.

\bdefn
Suppose $(M,\ga)$ is a balanced sutured manifold. A contact structure $\xi$ on $M$ is said to be \textbf{compatible} if $\partial M$ is convex and $\ga$ is the dividing set on $\partial M$.

A \textbf{contact handle} is a $3$-ball $B^3$ with the standard tight contact structure. The attachment of $B^3$ to a balanced sutured manifold $(M,\ga)$ is called a \textbf{contact $i$-handle attachment} in the following cases:
\benu
\item $i=0$ when the resulting manifold is a disjoint union $(M,\ga)\cup B^3$.
\item $i=1$ when $B^3$ is attached to $(M,\ga)$ along two points on the suture $\ga$.
\item $i=2$ when $B^3$ is attached to $(M,\ga)$ along a simple closed curve $\delta$ on $\partial M$ with $|\delta\cap\ga|=2$.
\eenu
\edefn

Suppose $(M,\ga)$ is a balanced sutured manifold. Let $(M',\ga')$ be the resulting manifold after attaching a contact $i$-handle. Baldwin and Sivek \cite[Section 3]{baldwin2016instanton} constructed a map 
\begin{equation}\label{eq: contact handle attaching map}
C: SHI(-M,-\ga)\ra SHI(-M',-\ga').
\end{equation}
We sketch the construction as follows.
\benu
\item When $i=0$ or $1$, we can construct the same closure for $(M,\ga)$ and $(M',\ga')$ and define $C$ to be the identity map.
\item When $i=2$, suppose $\delta\subset \partial M$ is the attaching curve of the contact handle. Then a closure of $(M',\ga')$ can be obtained from a closure of $(M,\ga)$ by performing a $0$-surgery along $\delta$, with respect to the framing from $\partial M$. Then $C$ is induced by the corresponding cobordism between closures.
\eenu

Suppose $(M,\ga)\subset(M^\p,\ga^\p)$ is a proper inclusion of balanced sutured manifolds and suppose $\xi$ is a contact structure compatible with $(M^\p\backslash{\rm int} M,\ga^\p\cup(-\ga))$. Based on maps associated to contact handle attachments, we can construct a \textbf{contact gluing map}\[\Phi_\xi:SHI(-M,-\ga)\ra SHI(-M^\p,-\ga^\p).\]The first author \cite{li2018gluing} showed that the contact gluing map is functorial, i.e.\, it is independent of the contact handle decompositions, and gluing two contact structures induces composite maps.

For a balanced sutured manifold $(M,\ga)$ and a compatible contact structure $\xi$, there are a few ways to decompose $\xi$ \cite{honda2009contact,baldwin2016instanton}.

{\bf Partial open book decomposition}. A partial open book decomposition is a triple $(S,P,h)$ where $S$ is a compact surface with non-empty boundary, $P\subset S$ a subsurface, and $h:P\ra S$ an embedding such that $h$ is the identity on $\partial P\cap \partial S$. 

{\bf Contact cellular decomposition}. A contact cellular decomposition of $\xi$ over $(M,\ga)$ is, roughly speaking, a Legendrian graph $\mathcal{K}\subset M$ such that $\partial \mathcal{K}\subset \ga$ and $M\backslash {\rm int}N(\mathcal{K})$ is diffeomorphic to a product $[-1,1]\times F$ for some surface $F$ withboundary and $\xi$ restricts to the $[-1,1]$-invariant contact structure on $M\backslash {\rm int}N(\mathcal{K})\cong [-1,1]\times F$.

{\bf Contact handle decomposition}. A contact handle decomposition is a decomposition of $(M,\ga,\xi)$ into contact $0$-, $1$-, and $2$-handles described above. 

These three decompositions can be related to each other as follows.

Suppose we have a contact cellular decomposition, i.e.\, a Legendrian graph $\mathcal{K}\subset M$ such that $M\backslash {\rm int}N(\mathcal{K})$ is a product manifold equipped with the product contact structure. Then $M\backslash {\rm int}N(\mathcal{K})$ equipped with the restriction of $\xi$ can be decomposed into a contact $0$-handle and a few contact $1$-handles. Furthermore, each edge of the Legendrian graph $\mathcal{K}$ corresponds to a contact $2$-handle attached along a meridian of the edge. This gives rise to a contact handle decomposition of $(M,\ga,\xi)$.

Suppose we have a contact handle decomposition of $(M,\ga,\xi)$; we can obtain a partial open book decomposition as follows. All $0$- and $1$- handles form a product sutured manifold $([-1,1]\times S,\{0\}\times\partial S)$. Suppose $2$-handles are attached along curves $\delta_1$,..., $\delta_n$. Let $P\subset \{1\}\times S$ be a neighborhood of $(\delta_1\cup\cdots\cup \delta_n)\cap \{1\}\times S$. Isotope $(\delta_1\cup\cdots\cup \delta_n)\cap \{-1\}\times S$ through $[-1,1]\times S$ onto $\{1\}\times S$. Let $h:P\ra S$ be the embedding such that $h|_{\partial S\cap\partial P}$ is the identity and $\delta_i\cap \{1\}\times S$ is sent to the image of $\delta_{i}\times \{-1\}\times S$ under the isotopy for $i=1,\dots,n$. Then $(S,P,h)$ is a partial open book decomposition of $(M,\ga,\xi)$.

Suppose we have a partial open book decomposition $(S,P,h)$ of $(M,\ga,\xi)$. We know that $([-1,1]\times S,\{0\}\times \partial S)$ is a product sutured manifold that admits a product contact structure $\xi_0$. This can be decomposed into a contact $0$-handle and a few contact $1$-handles. Let $a_1,...,a_n$ be a collection of disjoint properly embedded arcs on $S$ such that $a_i\subset P$ and $S-(a_1\cup\cdots\cup a_n)$ retracts to $S-P$. Let $\delta_i$ be the union of $a_i$ and $h(a_i)$. Then $(M,\ga,\xi)$ is obtained from $([-1,1]\times S,\{0\}\times \partial S,\xi_0)$ by attaching contact $2$-handles along all $\delta_i$.

\bdefn[\cite{baldwin2016instanton}]\label{defn: contact element in SHI}
Suppose $(M,\ga)$ is a balanced sutured manifold and $\xi$ is a compatible contact structure. Suppose $\xi$ has a partial open book decomposition $(S,h,P)$. Let $\delta_1$,..., $\delta_n$ be the attaching curves of the contact $2$-handles such that $(M,\ga,\xi)$ is obtained from $([-1,1]\times S,\{0\}\times \partial S)$ as above. Suppose the element $\mathbf{1}$ is the generator of 
$$\shi(-[-1,1]\times S,-\{0\}\times  \partial S)\cong \mathbb{C}.$$Then the {\bf instanton contact element} of $(M,\ga,\xi)$ is
$$\theta(M,\ga,\xi)\deq C_{\delta_n}\circ\cdots\circ C_{\delta_1}(\mathbf{1})\in SHI(-M,-\ga),$$
where $C_{\delta_i}$ is the contact gluing map associated to the contact 2-handle attachment along $\delta_i$.
\edefn

\bthm[Baldwin and Sivek \cite{baldwin2016instanton}]
Suppose $(M,\ga)$ is a balanced sutured manifold, and $\xi$ is a compatible contact structure. Then the instanton contact element $\theta(M,\ga,\xi)\in SHI(-M,-\ga)$ is independent of the choice of the partial open book decomposition and is well-defined up to a unit. In particular, the non-vanishing of the instanton contact element is an invariant property for the contact structure.
\ethm

Then we prove the main theorem of this subsection.

\bpf[Proof of Theorem \ref{thm: contact support}]
First, we prove the instanton contact element is homogeneous with respect to the $\intg$-grading of $SHI(-M,-\ga)$ associated to $S$. From \cite[Theorem 1.1]{honda2009contact}, any triple $(M,\ga,\xi)$ admits a contact cell decomposition. Hence, there exists a Legendrian graph $\mathcal{K}$, such that $(M\backslash{\rm int}N(\mathcal{K}),\xi|_{M\backslash{\rm int}N(\mathcal{K})})$ is contactomorphic to $([-1,1]\times F,\xi_0)$ for some surface $F$ with boundary and the product contact structure $\xi_0$. Let $\delta_1$,..., $\delta_n$ be a set of meridians of $K$, one for each edge of $\mathcal{K}$. Then we can obtain the original $\xi$ on $M$ from $([-1,1]\times F,\xi_0)$ by attaching contact $2$-handles along $\delta_1,\dots,\delta_n$. As discussed above, this gives rise to a contact handle decomposition and hence a partial open book decomposition. From Definition \ref{defn: contact element in SHI}, we know that
$$\theta(M,\ga,\xi)=C_{\delta_{n}}\circ\cdots\circ C_{\delta_1}(\mathbf{1})\in \shi(-M,-\ga),$$
where $C_{\delta_i}$ is the contact gluing map associated to the contact 2-handle attachment along $\delta_i$.

Suppose $S\subset (M,\ga)$ is an admissible surface. We can isotope $S$ so that it intersects $\mathcal{K}$ transversely and is disjoint from all $\delta_i$. Write $$S_{\mathcal{K}}=S\cap(M\backslash {\rm int}N(\mathcal{K})).$$We can consider it as a surface inside the product sutured manifold $([-1,1]\times S,\{0\}\times \partial S)$. Note that all components of $\partial S_{\mathcal{K}}\backslash \partial S$ are meridians of ${\mathcal{K}}$ and, by construction, each meridian of ${\mathcal{K}}$ has two intersections with the dividing set on $\partial (M\backslash {\rm int}N(\mathcal{K}))$, which is also identified with $$\{0\}\times\partial S\subset [-1,1]\times\{0\}\times S.$$So $S_{\mathcal{K}}$ is also admissible inside $([-1,1]\times S,\{0\}\times \partial S)$. Since
$$\shi(-[-1,1]\times S,-\{0\}\times \partial S)\cong\mathbb{C},$$
we know that there exists $i_0\in\intg$ such that
$$\mathbf{1}\in \shi(-[-1,1]\times S,-\{0\}\times \partial S,S_{\mathcal{K}},i_0).$$From \cite[Proposition 4.6]{LY2021}, we know that all maps $C_{\delta_i}$ preserve the gradings associated to $S_{\mathcal{K}}$ and $S$, respectively. Thus, we conclude that
$$\theta(M,\ga,\xi)=C_{\delta_{n}}\circ\cdots\circ C_{\delta_1}(\mathbf{1})\in \shi(-M,-\ga,i_0).$$

Then we need to figure out $i_0$. Since $\shi(-[-1,1]\times S,-\{0\}\times \partial S)$ is one-dimensional, the integer $i_0$ is determined by its graded Euler characteristic (we fix the closure to resolve the ambiguity of $\pm H$). By \cite[Proposition 4.3 and Corollary 3.45]{LY2021} (see also \cite[Theorem 3.26]{Baldwin2020}), it suffices to calculate $i_0$ when replacing $\shi$ by $SFH$. Note that the contact element of any contact structure $\xi$ compatible with $(M,\ga)$ lives in $SFH(-M,-\ga,\mathfrak{s}_\xi)$, where $\mathfrak{s}_\xi$ is the relative spin$^c$ structure corresponding to $\xi$. The formula of $i_0$ then follows from \cite[Proposition 4.5]{honda2000classification}.
\epf

\subsection{Vanishing results about Giroux torsion}\quad

Instanton contact elements share similar properties with the contact elements in $SFH$. To prove the vanishing result about Giroux torsion for the instanton contact element (Theorem \ref{thm: SHI giroux torsion}), we need to prove lemmas similar to Lemma \ref{lem: vanishing piece SFH} and Lemma \ref{lem: vanishing result SFH}.

The analog of Lemma \ref{lem: vanishing piece SFH} follows directly from the following proposition.
\bprop[{\cite[Corollary 1.4]{li2018gluing}, see also \cite[Theorem 1.2]{baldwin2016instanton}}]\label{prop: functoriality of contact element}
Consider the notations as above. If the contact structure $\xi$ on $M^\p\backslash{\rm int} M$ is a restriction of a contact structure $\xi^\p$ on $M^\p$, then we have$$\Phi_\xi(\theta(M,\ga,\xi^\p|_M))=\theta(M^\p,\ga^\p,\xi^\p)\in \shi(-M^\p,-\ga^\p).$$
\eprop
\bcor\label{cor: vanishing piece SHI}
Let $(Y, \xi)$ be a closed contact 3-manifold and $N \subset Y$ be a compact submanifold (without any closed components) with convex boundary and dividing set $\Ga$. If $\theta(N, \Ga, \xi|_N) = 0$, then $\theta(Y, \xi) = 0$.
\ecor

The following proposition is the analog of Lemma \ref{lem: vanishing result SFH}.
\bprop\label{prop: vanishing result SHI}
The instanton contact elements $\theta(N_*,\Ga_*,\zeta_1^+)$ and $\theta(N_*,\Ga_*,\zeta_1^-)$ vanish.
\eprop
\bpf

Since instanton contact elements share most properties with contact elements, we can apply the proof of Lemma \ref{lem: vanishing result SFH} with mild changes. We sketch the proof and point out the main difference. For simplicity, we only consider $\theta(N_*,\Ga_*,\zeta_1^+)$. The proof for $\theta(N_*,\Ga_*,\zeta_1^-)$ is almost identical. 

Take a copy $T_\varepsilon=T^2\times \{\varepsilon\}\subset {\rm int}N_*$ with the dividing set consisting of two curves of slope $\infty$. Let $L$ be a Legendrian ruling curve on $T_\varepsilon$ with slope $-1$ (cf.\ Remark \ref{rem: change slope}). The Legendrian curve $L$ has twisting number $-1$ with respect to the framing coming from $T_\varepsilon$. Let $(N^\p,\Ga^\p,(\zeta_1^+)^\p)$ be obtained from $(N_*,\Ga_*,\zeta_1^+)$ by a contact $(+1)$-surgery along $L$. By \cite[Theorem 4.6]{baldwin2016instanton}, the cobordism map $\Phi$, corresponding to the contact $(+1)$-surgery, sends $\theta(N_*,\Ga_*,\zeta_1^+)$ to $\theta((N^\p,\Ga^\p,(\zeta_1^+)^\p))=0$. By \cite[Lemma 7]{girouxtorsion1}, the resulting contact structure $(\zeta_1^+)^\p$ is overtwisted. Hence, by \cite[Theorem 1.3]{baldwin2016instanton}, we have $\theta((N^\p,\Ga^\p,(\zeta_1^+)^\p))=0$. It remains to show $\Phi$ is injective (at least on the subspace generated by $\theta(N_*,\Ga_*,\zeta_1^+)$).

Write $(N_*,\Ga_*,\zeta_0^+)$ for $N_0$. In the proof of Lemma \ref{lem: vanishing result SFH}, by considering the relative spin$^c$ structure, the authors of \cite{girouxtorsion1} showed that $c(N_*,\Ga_*,\zeta_0^+)$ and $c(N_*,\Ga_*,\zeta_1^+)$ lie in the same $\mathbb{F}_2$ summand of $SFH(-N_*,-\Ga_*)\cong \mathbb{F}_2^4$ (we replace $\mathbb{Z}$-summand by $\mathbb{F}_2$ summand for the naturality issue, cf.\ Remark \ref{rem: Z coeff}). The contact structure $\zeta_0^+$ and the contact structure $(\zeta_0^+)^\p$ after the contact $(+1)$-surgery along $L$ can be embedded into $S^3$ and $S^1\times S^2$ with standard tight contact structures, respectively, which are both Stein fillable. Then both $c(N_*,\Ga_*,\zeta_0^+)$ and $c(N^\p,\Ga^\p,(\zeta_0^+)^\p)$ are non-vanishing. Thus, the map $\Phi$ is injective on the $\ft$ summand generated by $c(N_*,\Ga_*,\zeta_0^+)$.

For sutured instanton homology, the analog of the (nontorsion) relative spin$^c$ decomposition is the decomposition associated to admissible surfaces, constructed in \cite{li2019decomposition,li2019direct}. We can use two annuli $$A_0=S^1\times\{{\rm pt}\}\times I, A_1=\{{\rm pt}\}\times S^1 \times I\subset T^2\times I$$to construct the decomposition, where the $S^1$ factors correspond to curves of slopes $\infty$ and $0$ parallel to the dividing sets, respectively. Since $|\partial A_i\cap \Ga_*|=2$ for $i=0,1$, by \cite[Theorem 2.21]{LY2020} there are only two nontrivial gradings for $A_i$, corresponding to the sutured manifold decomposition along $A_i$ and $-A_i$. It is straightforward to check that sutured manifold decomposition along $\pm A_0\cup\pm A_1$ gives a 3-ball with a connected suture, whose $\shi$ is 1-dimensional. Thus, $$\dim_\mathbb{C}\shi(-N_*,-\Ga_*)=4.$$




By Proposition \ref{thm: contact support}, we know that $\theta(N_*,\Ga_*,\zeta_1^+)$ and $\theta(N_*,\Ga_*,\zeta_0^+)$ live in the same grading. Since $\shi$ is 1-dimensional in any nontrivial grading, the elements $\theta(N_*,\Ga_*,\zeta_1^+)$ and $\theta(N_*,\Ga_*,\zeta_0^+)$ are linearly dependent. By \cite[Corollary 1.6]{baldwin2016instanton} and the Stein fillability, both $\theta(N_*,\Ga_*,\zeta_0^+)$ and $\theta(N^\p,\Ga^\p,(\zeta_0^+)^\p)$ are non-vanishing. Then $\Phi$ is injective on the subspace generated by $\theta(N_*,\Ga_*,\zeta_0^+)$, and $\Phi(\theta(N_*,\Ga_*,\zeta_1^+))=0$ implies $\theta(N_*,\Ga_*,\zeta_1^+)=0$. 

\epf
\bpf[Proof of Theorem \ref{thm: SHI giroux torsion}]
 This follows from Lemma \ref{lem: giroux torsion}, Corollary \ref{cor: vanishing piece SHI}, and Proposition \ref{prop: vanishing result SHI}. Note that Lemma \ref{lem: giroux torsion} is only about contact topology, so we can apply it without change. 
\epf

\subsection{Vanishing results about cobordism maps}\quad

Suppose $(M,\ga)\subset(M^\p,\ga^\p)$ is a proper inclusion of balanced sutured manifolds and suppose $\xi$ is a contact structure compatible with $(M^\p\backslash{\rm int} M,\ga^\p\cup(-\ga))$. By Corollary \ref{cor: vanishing piece SHI}, if $$\theta(M^\p\backslash{\rm int} M,\ga^\p\cup(-\ga),\xi)=0,$$then the contact gluing map $\Phi_\xi$ vanishes on the subspace of $\shi(-M,-\ga)$ generated by instanton contact elements. Indeed, we can prove a stronger result by the functoriality of $\Phi_\xi$. The proof of the following proposition is due to Ian Zemke.

\bprop\label{prop: vanishing cobordism}
Suppose $(M,\ga)\subset(M^\p,\ga^\p)$ is a proper inclusion of balanced sutured manifolds and suppose $\xi$ is a contact structure compatible with $$(M_0,\ga_0)\deq (M^\p\backslash{\rm int} M,\ga^\p\cup(-\ga)).$$If the contact element $\theta(M_0,\ga_0,\xi)$ vanishes, then the map $\Phi_\xi$ vanishes on $\shi(-M,-\ga)$.
\eprop
\bpf
We have inclusions$$(M,\ga)\subset (M,\ga)\sqcup(M_0,\ga_0)\subset (M^\p,\ga^\p),$$where $\sqcup$ denotes the disjoint union. The manifold $$M^\p\backslash{\rm int}(M\sqcup M_0)$$is contactomorphic to $\partial M\times I$. Let $\xi_0$ be the product contact structure on $\partial M\times I$. By the instanton analog of \cite[Proposition 6.5]{kronheimer2010knots}, we have$$\shi(-M\sqcup (-M_0),-\ga\sqcup (-\ga_0))\cong \shi(-M,-\ga)\otimes \shi(-M_0,-\ga_0).$$By functoriality, the map $\Phi_\xi$ is the composition of the following maps$$\begin{aligned}
    \shi(-M,-\ga)\to  &\shi(-M,-\ga)\otimes \shi(-M_0,-\ga_0)\to \shi(-M^\p,-\ga^\p)\\
    x\quad\quad \mapsto &\quad\quad\quad \quad \quad  x\otimes \theta(M_0,\ga_0,\xi)\quad \quad \mapsto \Phi_{\xi_0}(x\otimes \theta(M_0,\ga_0,\xi)).
\end{aligned}$$If $\theta(M_0,\ga_0,\xi)=0$, then $\Phi_\xi=0$.
\epf
\brem
For a general balanced sutured manifold $(M,\ga)$, instanton contact elements do not generate $\shi(-M,-\ga)$ because the number of tight contact structures compatible with $(M,\ga)$ is less than $\dim_\mathbb{C}\shi(M,\ga)$. See \cite[Section 4.3]{li2019direct} and \cite{honda2000classification} for a discussion about contact structures on the solid torus.
\erem
The following vanishing result is used in the rest of the paper.
\bcor\label{cor: important vanishing}
Suppose $(M,\ga)\subset(M^\p,\ga^\p)$ is a proper inclusion of balanced sutured manifolds. Suppose we have $$(M^\p\backslash{\rm int} M,\ga^\p\cup(-\ga),\xi)=(N_*,\Ga_*,\zeta_1^+)\text{ or }(N_*,\Ga_*,\zeta_1^-)$$as defined in Subsection \ref{subsec: Vanishing results about Giroux torsion}. Then, $\Phi_\xi=0$.
\ecor
\bpf
This follows from Proposition \ref{prop: vanishing result SHI} and Proposition \ref{prop: vanishing cobordism}
\epf

\section{Instanton L-space knots}\label{sec: Knots with instanton L-space surgeries}\quad

In this section, we study the instanton knot homology of an instanton L-space knot $K\subset Y$. In particular, we prove Theorem \ref{thm: main 2}, Theorem \ref{main: L space knot}, and Theorem \ref{thm: floer simple khi}. For technical reasons, we only deal with the case $H_1(Y(K))\cong \mathbb{Z}$.



\subsection{The dimension in each grading}\label{subsec: Dimensions of instanton knot homology in all gradings}\quad

In this subsection, we prove the following theorem. The main input is the large surgery formula and the vanishing result from Corollary \ref{cor: important vanishing}.
\bthm\label{thm: KHI for L-space knots}
Suppose $Y$ is an integral homology sphere with $I^{\sharp}(Y)\cong \mathbb{C}$. Suppose $K\subset Y$ is a knot and $S$ is the Seifert surface of $K$. If there is a positive integer $n$ such that $Y_{-n}(K)$ is an instanton L-space, then for any $i\in\mathbb{Z}$, we have
$$\dim_\mathbb{C}\khii(-Y,K,S,i)\le 1.$$
\ethm

Since $Y$ is an integral homology sphere, $K$ is always null-homologous and $\hat{\mu}=\mu,\hat{\lambda}=\lambda$ in Subsection \ref{subsec: basic setups}. By Definition \ref{defn_4: Ga_n-hat sutures}, we have $(q,p)=(1,0)$ and $(q_0,p_0)=(0,1)$. Then we have $$\widehat{\Ga}_\mu=\Ga_\mu=\ga_{\mu},\quad\widehat{\Ga}_n=\Ga_n=\ga_{\lambda-n\mu}.$$Note that in the proof of Theorem \ref{thm: homology of bent complex}, an auxiliary slope $\hat{\mu}^\p=n\hat{\mu}-\hat{\lambda}$ is used. Here we set $\hat{\mu}^\p=n\mu-\lambda$. Since $n$ is not fixed, this slope is also not fixed. 

For simplicity, we write $\ga_{(x,y)}$ for $\ga_{x\lambda+y\mu}$ in Definition \ref{defn_4: Ga_n-hat sutures}. Also, we omit $S$ in the notation $\shi(-Y(K),\ga,i)$ for any $\ga$.

Then we make the following definition.
\bdefn\label{defn: top and bottom gradings}
For any integers $n$ and $i$ with $|i|\le g(K)$, define
\begin{equation*}
\begin{aligned}
    T_{n,i}&=\shi(-Y(K),-\Ga_n,i+\lceil\frac{n-1}{2}\rceil),\\
    B_{n,i}&=\shi(-Y(K),-\Ga_n,i-1+\lceil-\frac{n-1}{2}\rceil).
\end{aligned}
\end{equation*}
For $i>g(K)$ and any $n$, define $T_{n,i}=0$. For $i<-g(K)$ and any $n$, define $B_{n,i}=0.$
\edefn
\brem\label{rem: iso TB}
The notations `T' and `B' mean `top' and `bottom'. If we use the notations after the diagram (\ref{eq: useful triangles}) and suppose $g=g(K)$, then for any integers $n$ and $i$ with $|i|\le g(K)$, we have$$T_{n,i}=\widehat{\Ga}_n^{i,+}\aand B_{n,i}=\widehat{\Ga}_{n-1}^{i,-}.$$By Lemma \ref{lem_4: psi pm is isomorphism}, we have
$$\psi_{-,n+1}^n:T_{n,i}\xra{\cong} T_{n+1,i}~{\rm~and}~\psi_{+,n+1}:B_{n,i}\xra{\cong} B_{n+1,i}$$
for $n\geq 2g(K)+1$ and $|i|\le g(K)$.
\erem
The following proposition follows from the large surgery formula.
\bprop\label{prop: structure of gamma_2n+1}
Suppose $Y$ is an integral homology sphere with $I^{\sharp}(Y)\cong \mathbb{C}$. Suppose $K\subset Y$ is a knot. Suppose $n$ is an integer such that $n\geq 2g(K)+1$ and $Y_{-n}(K)$ is an instanton L-space. Then we have the following.
\begin{equation*}
\shi(-Y(K),-\ga_{(2,1-2n)},i)\cong\begin{cases}
	T_{n,i-n+1}&n-g\leq i\leq n-1+g\\
	\mathbb{C}&-n+g+1\leq i\leq n-g-1\\
	B_{n,i+n-1}&-n+1-g\leq i\leq-n+g 
\end{cases}
\end{equation*}
\eprop
\bpf
The isomorphism of the top and bottom $2g$ gradings of $\shi(-Y(K),-\ga_{(2,1-2n)})$ follows from applying Lemma \ref{lem_4: psi pm is isomorphism} to $\hat{\mu}^\p$. Since $Y_{-n}(K)$ is an instanton L-space, by (\ref{eq: mirror iso}), the manifold $-Y_{-n}(K)$ is also an instanton L-space. The isomorphism of the middle gradings follows from Proposition \ref{prop: l space dim 1}, Lemma \ref{lem_4: q-cyclic}, and Theorem \ref{thm: large surgery formula, main}.
\epf

Note that in the proof of Theorem \ref{thm: homology of bent complex} (more precisely, in the triangle (\ref{eq: triangle a})), we have a map $\psi_{-,0}^\mu(\hat{\mu}^\p)$ from the space associated to $\widehat{\Ga}_n$ to the space associated to $\widehat{\Ga}_{n-1}$. We write this map as $\psi_{-,n-1}^n$. We also write $\psi_{-,n}^{2n-1}$ and $\psi_{-,2n-1}^{n-1}$ for $\psi_{-,\mu}^1(\hat{\mu}^\p)$ and $\psi_{-,1}^0(\hat{\mu}^\p)$ in (\ref{eq: triangle a}), respectively. Similarly we write $\psi_{+,n-1}^n,\psi_{+,n}^{2n-1},$ and $\psi_{+,2n-1}^{n-1}$ for maps in the positive bypass triangle. We abuse notation so that bypass maps also denote their restrictions to a single grading. Then the following proposition follows from the vanishing results established in Section \ref{sec: Vanishing results for contact elements}.
\bprop\label{prop: composition of bypasses is zero}
Suppose $K\subset Y$ is a null homologous knot. For any integer $n\in\intg$ with $n\geq 2g(K)+1$ and any integer $i$ with $|i|\leq g(K)$, we have
$$\psi_{+,n}^{n+1}\circ\psi_{-,n+1}^{n+2}=0:T_{n+2,i}\ra T_{n,i}$$
and
$$\psi_{-,n}^{n+1}\circ\psi_{+,n+1}^{n+2}=0:B_{n+2,i}\ra B_{n,i}.$$
\eprop
\bpf
By Remark \ref{rem: iso TB}, it suffices to prove$$\Psi_T\deq\psi_{-,n+3}^{n+2}\circ\psi_{-,n+2}^{n+1}\circ\psi_{-,n+1}^n\circ\psi_{+,n}^{n+1}\circ\psi_{-,n+1}^{n+2}=0:T_{n+2,i}\ra T_{n+3,i}$$
and
$$\Psi_B\deq\psi_{+,n+3}^{n+2}\circ\psi_{+,n+2}^{n+1}\circ\psi_{+,n+1}^n\circ\psi_{-,n}^{n+1}\circ\psi_{+,n+1}^{n+2}=0:B_{n+2,i}\ra B_{n+3,i}.$$
By classification of tight contact structures on $T^2\times I$ \cite{honda2000classification}, we know that the contact structures corresponding to $\Psi_T$ and $\Psi_B$ are contactomorphic to either $(N_*,\Ga_*,\zeta_1^+)$ or $(N_*,\Ga_*,\zeta_1^-)$ defined in Subsection \ref{subsec: Vanishing results about Giroux torsion}. Then the lemma follows from Corollary \ref{cor: important vanishing}.
\epf
\begin{prop}\label{prop: large enough surgeries are all L-spaces}
Suppose $Y$ is an integral homology sphere with $I^{\sharp}(Y)\cong \mathbb{C}$. Suppose $K\subset Y$ is a knot. Suppose $n_0$ is a positive integer such that $Y_{-n_0}(K)$ is an instanton L-space. Then for any integer $n$ such that $n>n_0$,  $Y_{-n}(K)$ is also an instanton L-space.
\end{prop}
\bpf
This proposition follows immediately from $\chi(I^\sharp(Y_{-n}(K))=|H_1(Y_{-n}(K))|$, the equation $$|H_1(Y_{-n-1}(K))|=|H_1(Y_{-n}(K))|+|H_1(Y)|,$$and the following surgery exact triangle (\cite[Section 4.2]{baldwin2018stein}, see also \cite{scaduto2015instanton})
\begin{equation*}
\xymatrix{
I^{\sharp}(Y_{-n-1}(K))\ar[rr]&&I^{\sharp}(Y_{-n}(K))\ar[dl]\\
&I^{\sharp}(Y)\ar[ul]&
}	
\end{equation*}

\epf

By Proposition \ref{prop: structure of gamma_2n+1} and Proposition \ref{prop: composition of bypasses is zero}, the proof of Theorem \ref{thm: KHI for L-space knots} follows from similar algebraic lemmas in \cite[Section 3]{Ozsvath2005}. We reprove them in our setting.
\begin{lem}\label{lem:  Knots admitting L-space surgeries, case 1}
Suppose $Y$ is an integral homology sphere with $I^{\sharp}(Y)\cong \mathbb{C}$. Suppose $K\subset Y$ is a knot. Suppose $n_0$ is a positive integer such that $Y_{-n_0}(K)$ is an instanton L-space. Suppose further that for a large enough integer $n$ and some integer $m$ with $|m|\leq g(K)$, we have
$T_{n,m+1}=0$. Then one of the following two cases happens. 
\benu
\item $\khii(-Y,K,m)\cong\mathbb{C}$ and $B_{n,m-1}=0$,
\item $\khii(-Y,K,m)=0$ and $T_{n,m}=0$.
\eenu
\end{lem}

\bpf
By Proposition \ref{prop: large enough surgeries are all L-spaces}, we can take any arbitrary large enough integer $n$, since they are all L-space surgery slopes. From Proposition \ref{prop_4: graded bypass for Ga_n-hat sutures}, we have the following exact triangle
\begin{equation*}
\xymatrix{
T_{n-1,m+1}\ar[rr]&&T_{n,m}\ar[dl]\\
&\khii(-Y,K,m)\ar[lu]&
}	
\end{equation*}
From Remark \ref{rem: iso TB} and the assumption $T_{n,m+1}=0$, we know that
$$T_{n-1,m+1}\cong T_{n,m+1}=0~{\rm and~}B_{n,m-1}\cong B_{n-1,m-1}.$$
Hence there exists some $k\in\mathbb{N}$ such that 
$$T_{n,m}\cong \khii(-Y,K,m)\cong \mathbb{C}^k.$$
Also, from Proposition \ref{prop_4: graded bypass for Ga_n-hat sutures}, we have the following exact diagram
\begin{equation*}
	\xymatrix{
	&\shi(-Y(K),-\ga_{(2,1-2n)},m)\ar[d]&\\
	&T_{n,m}\ar[d]^{\psi_{-,n-1}^{n,m}}&\\
	\shi(-Y(K),-\ga_{(2,3-2n)},m-1)\ar[r]^{\quad\quad\quad\quad\psi_{+,n-1}^{2n-3,m-1}}&B_{n-1,m-1}\ar[r]^{\psi_{+,n-2}^{n-1,m-1}}&T_{n-2,m}
	}
\end{equation*}
where $\psi_{-,n-1}^{n,m}$ is the map $\psi_{-,n-1}^{n}$ restricted to the graded part $T_{n,m}$ and other notations are defined similarly. Since $|m|\leq g(K)$, Proposition \ref{prop: structure of gamma_2n+1} implies that
$$\shi(-Y(K),-\ga_{(2,1-2n)},m)\cong \shi(-Y(K),-\ga_{(2,3-2n)},m-1)\cong \mathbb{C}.$$
Hence the above diagram can be rewritten as
\begin{equation}\label{eq: diagram a}
	\xymatrix{
	&&\mathbb{C}\ar[d]&\\
	&&T_{n,m}\cong \mathbb{C}^k\ar[d]^{\psi_{-,n-1}^{n,m}}&\\
	\mathbb{C}\ar[rr]^{\psi_{+,n-1}^{2n-3,m-1}\quad\quad}&&B_{n-1,m-1}\ar[rr]^{\psi_{+,n-2}^{n-1,m-1}}&&T_{n-2,m}\cong\mathbb{C}^{k}
	}
\end{equation}
We consider the following two cases.

{\bf Case 1}. $\psi_{+,n-1}^{2n-3,m-1}$ is trivial. Then from the exactness of the horizontal sequence in (\ref{eq: diagram a}), we know that $B_{n-1,m-1}\cong \mathbb{C}^{k-1}$
and $\psi_{+,n-2}^{n-1,m-1}$ is injective. Also, we conclude from the exactness of the vertical sequence in (\ref{eq: diagram a}) that $\psi_{-,n-1}^{n,m}$ is surjective. However, from Proposition \ref{prop: composition of bypasses is zero} we know that
$$\psi_{+,n-2}^{n-1,m-1}\circ\psi_{-,n-1}^{n,m}=0.$$
Hence the only possibility is that $k=1$, and this concludes that $T_{n,m}\cong \khii(-Y,K,m)\cong \mathbb{C}$, and $B_{n,m-1}\cong B_{n-1,m-1}=0$, which is the first case in the statement of the lemma.

{\bf Case 2}. $\psi_{+,n-1}^{2n-3,m-1}$ is nontrivial. Then from the exactness of the horizontal sequence in (\ref{eq: diagram a}), we know that $B_{n-1,m-1}\cong \mathbb{C}^{k+1}$ and $\psi_{+,n-2}^{n-1,m-1}$ is surjective. From the above discussion and the bypass exact triangle from Proposition \ref{prop_4: graded bypass for Ga_n-hat sutures}, we have another exact diagram
\begin{equation}\label{eq: diagram b}
	\xymatrix{
	&&\shi(-Y(K),-\ga_{(2,5-2n),m})\cong\mathbb{C}\ar[d]\\
	B_{n-1,m-1}\cong\mathbb{C}^{k+1}\ar[rr]^{\psi_{+,n-2}^{n-1,m-1}}&&T_{n-2,m}\cong\mathbb{C}^{k}\ar[d]^{\psi_{-,n-3}^{n-2,m}}\\
	&&B_{n-3,m-1}\cong \mathbb{C}^{k+1}
	}
\end{equation}
The exactness of the vertical sequence in (\ref{eq: diagram b}) implies that the map $\psi_{-,n-3}^{n-2,m}$ is injective. However, from Proposition \ref{prop: composition of bypasses is zero}, we have 
$$\psi_{-,n-3}^{n-2,m}\circ\psi_{+,n-2}^{n-1,m-1}=0.$$
Hence the only possibility is that $k=0$. Thus, we conclude that $T_{n,m}\cong \khii(-Y,K,m)=0$, which is the second case in the statement of the lemma.
\epf

\blem\label{lem:  Knots admitting L-space surgeries, case 2}
Suppose $Y$ is an integral homology sphere with $I^{\sharp}(Y)\cong \mathbb{C}$. Suppose $K\subset Y$ is a knot. Suppose $n_0$ is a positive integer such that $Y_{-n_0}(K)$ is an instanton L-space. Suppose further that for a large enough integer $n$ and some integer $m$ with $|m|\leq g(K)$, we have $B_{n,m}=0$. Then one of the following two cases happens.
\benu
\item $\khii(-Y,K,m)\cong\mathbb{C}$ and $T_{n,m}=0$,
\item $\khii(-Y,K,m)=0$ and $B_{n,m-1}=0$.
\eenu
\elem
\bpf
The proof is similar to that of Lemma \ref{lem:  Knots admitting L-space surgeries, case 1}. From Proposition \ref{prop_4: graded bypass for Ga_n-hat sutures}, we have the following triangle
\begin{equation*}
\xymatrix{
B_{n-1,m-1}\ar[rr]&&B_{n,m}\ar[dl]\\
&\khii(-Y,K,m)\ar[lu]&
}	
\end{equation*}
Hence there exists some $k\in\mathbb{N}$ such that
$$B_{n-1,m-1}\cong \khii(-Y,K,m)\cong \mathbb{C}^k.$$

Also from Proposition \ref{prop_4: graded bypass for Ga_n-hat sutures}, we have the following  exact diagram
\begin{equation}\label{eq: diagram c}
	\xymatrix{
	&&\mathbb{C}\ar[d]&&\\
	&&T_{n,m}\ar[d]^{\psi_{-,n-1}^{n,m}}&&\mathbb{C}\ar[d]\\
	\mathbb{C}\ar[rr]^{\psi_{+,n-1}^{2n-3,m-1}\quad\quad}&&B_{n-1,m-1}\cong \mathbb{C}^k\ar[rr]^{\psi_{+,n-2}^{n-1,m-1}}&&T_{n-2,m}\ar[d]^{\psi_{-,n-3}^{n-2,m-1}}\\
	&&&&B_{n-3,m-1}\cong \mathbb{C}^k
	}
\end{equation}
We consider the following two cases.

{\bf Case 1}. $\psi_{+,n-1}^{2n-3,m-1}$ is trivial. Then from the exactness of the horizontal sequence in (\ref{eq: diagram c}), we know that $T_{n-2,m}\cong \mathbb{C}^{k-1}$ and $\psi_{+,n-2}^{n-1,m-1}$ is surjective. Also, we conclude from the exactness of the second vertical sequence in (\ref{eq: diagram c}) that $\psi_{-,n-3}^{n-2,m}$ is injective. However, from Proposition \ref{prop: composition of bypasses is zero} we know that
$$\psi_{-,n-3}^{n-2,m}\circ\psi_{+,n-2}^{n-1,m-1}=0.$$
Hence the only possibility is that $k=1$. Hence, we conclude that $\khii(-Y,K,m)\cong \mathbb{C}$ and $T_{n,m}\cong T_{n-2,m}=0$, which is the first case in the statement of the lemma.

{\bf Case 2}. $\psi_{+,n-1}^{2n-3,m-1}$ is nontrivial. Then from the exactness of the horizontal sequence in (\ref{eq: diagram c}), we know that $T_{n,m}\cong T_{n-2,m}\cong \mathbb{C}^{k+1}$ and $\psi_{+,n-2}^{n-1,m-1}$ is injective. Also, we conclude from the exactness of the first vertical sequence that $\psi_{-,n-1}^{n,m}$ is surjective. However, from Proposition \ref{prop: composition of bypasses is zero} we know that
$$\psi_{+,n-2}^{n-1,m-1}\circ\psi_{-,n-1}^{n,m}=0.$$
Hence the only possibility is that $k=0$, and this concludes that
$$B_{n,m-1}\cong B_{n-1,m-1}\cong \khii(-Y,K,m)\cong \mathbb{C}^k,$$ which is the second case in the statement of the lemma.
\epf

\bpf[Proof of Theorem \ref{thm: KHI for L-space knots}]
By Definition \ref{defn: top and bottom gradings} and Lemma \ref{lem_4: top and bottom non-vanishing gradings}, we know that
$$T_{n,g(K)+1}=0\aand  \khii(-Y,K,g(K)+1)=0.$$
We apply an induction that decreases the integer $i$: assuming that for $i+1$, we have
$$\khii(-Y,K,i+1)\cong \mathbb{C}~{\rm~or}~0$$
and either $T_{n,i+1}=0$ or $B_{n,(i+1)-1}=0$, we then want to prove the same results for $i$. When $T_{n,i+1}=0$, from Lemma \ref{lem:  Knots admitting L-space surgeries, case 1}, we have either $\khii(-Y,K,i)\cong\mathbb{C}$ and $B_{n,i-1}=0$ or $\khii(-Y,K,i)=0$ and $T_{n,i}=0$. When $B_{n,(i+1)-1}=0$, from Lemma \ref{lem:  Knots admitting L-space surgeries, case 2}, we have either $\khii(-Y,K,i)\cong\mathbb{C}$ and $T_{n,i}=0$ or $\khii(-Y,K,i)=0$ and $B_{n,i-1}=0$. Hence, the inductive step is completed, and we conclude that
$$\khii(-Y,K,i)\cong \mathbb{C}~{\rm~or~}0.$$
for all $i\in\intg$ such that $|i|\leq g(K)$. From Lemma \ref{lem_4: top and bottom non-vanishing gradings}, we know that
$$\khii(-Y,K,i)\cong 0$$
for all $i\in\intg$ with $|i|>g(K)$. Hence, we conclude the proof of Theorem \ref{thm: KHI for L-space knots}.
\epf

% ...

\subsection{Coherent chains}\quad

In this subsection, we prove the instanton analog of \cite[Lemma 3.2]{Rasmussen2017} with more assumptions. First, we introduce the analog of \cite[Definition 3.1]{Rasmussen2017} in instanton theory.

\bdefn\label{defn: positive chains}
Suppose $K$ is a knot in a rational homology sphere $Y$ and suppose $\hat{\mu}$ is the meridian of $K$. Suppose the knot complement $Y(K)$ satisfies $H_1(Y(K))\cong\mathbb{Z}$ so that we can identify $[\hat{\mu}]\in H_1(Y(K))$ as an integer $q$. Indeed, if a Seifert surface $S$ of $K$ is chosen, we can set $q=S\cdot \hat{\mu}$. For any integer $s$ and its image $[s]\in\mathbb{Z}_q$, define

$$\khii(-Y,K,[s])\deq\bigoplus_{k\in\mathbb{Z}}\khii(-Y,K,S,s+kq).$$
It is called a \textbf{positive chain} if it is generated by elements $$x_1,\dots, x_l, y_1,\dots, y_{l-1},$$each of which lives in a single grading associated to $S$ and a single $\mathbb{Z}_2$-grading, and the differentials $d_+$ and $d_-$ satisfy $$d_-(y_i) \doteq  x_{i+1},\quad d_+(y_i) \doteq  x_i, \aand d_-(x_i) = d_+(x_i) = 0 \text{ for all }i,$$where $\doteq$ means equal up to multiplication by a unit. The space $\khii(-Y,K,[s])$ is called a \textbf{negative chain} if there exist similar generators so that $$d_-(x_i) \doteq  y_{i},\quad d_+(x_i) \doteq  y_{i-1}, \aand d_-(y_i) = d_+(y_i) = 0 \text{ for all }i.$$We say $\khii(-Y,K)$ \textbf{consists of positive chains} if $\khii(-Y,K,[s])$ is a positive chain for any $[s]\in \mathbb{Z}_q$ and \textbf{consists of negative chains} if $\khii(-Y,K,[s])$ is a negative chain for any $[s]\in \mathbb{Z}_q$. We say $\khii(-Y,K)$ \textbf{consists of coherent chains} if $\khii(-Y,K)$ either consists of positive chains or consists of negative chains
\edefn
\brem\label{rem: chain mirror}
By Definition \ref{defn: positive chains}, the space $\khii(-Y,K,[s])$ is both a positive chain and a negative chain if and only if $\dim_\mathbb{C}\khii(-Y,K,[s])=1$. By the proof of Proposition \ref{prop: dual complex}, the space $\khii(-Y,K)$ consists of positive chains if and only if $\khii(Y,K)$ consists of negative chains.
\erem


The main theorem in this subsection is the following.
\bthm[]\label{thm: positive chain}
Suppose $K\subset Y$ is a knot as in Definition \ref{defn: positive chains}. Note that $H_1(Y(K))\cong\mathbb{Z}$. Suppose $Y$ is an instanton L-space and suppose $n\in\mathbb{N}_+$. Suppose the basis $(\hat{\mu},\hat{\lambda})$ of $\partial Y(K)$ is from Definition \ref{defn_4: Ga_n-hat sutures}. If $Y_{-n}(K)$ is an instanton L-space, then $\khii(-Y,K)$ consists of positive chains. If $Y_{n}(K)$ is an instanton L-space, then $\khii(-Y,K)$ consists of negative chains. 
\ethm

For simplicity, we only provide details of the proof for a special case of Theorem \ref{thm: positive chain}. The proof for the general case is similar. The main input is Theorem \ref{thm: KHI for L-space knots}.

\bdefn
We adopt the notations in Subsection \ref{subsec: Dimensions of instanton knot homology in all gradings} and Construction \ref{cons: bent complex}. For any integer $s$, suppose $B_{\ge s}^{+}$ is the subcomplex of $B_s^+$ with the underlying space $$\bigoplus_{k\ge s}\shi(-Y(K),-\widehat{\Ga}_\mu,S,s+kq)$$and suppose $B_{<s}^-$ is the subcomplex of $B_s^-$ with the underlying space $$\bigoplus_{k< s}\shi(-Y(K),-\widehat{\Ga}_\mu,S,s+kq).$$Let $H(B_{\ge s}^{+})$ and $H(B_{<s}^-)$ be the corresponding homologies. 
\edefn
\blem[]\label{lem: tb iso}
For any integers $n$ and $i$ with $|i|\le g(K)$, we have$$T_{n,i}\cong H(B_{\ge i}^{+})\aand B_{n,i}\cong H(B_{<i}^-).$$

\elem
\bpf
This follows from Remark \ref{rem: iso TB}, equations (\ref{eq: vanishing a}) and (\ref{eq: vanishing b}), and Theorem \ref{thm: convergence}.
\epf
\bthm[]\label{thm: positive chain, integral homology sphere}
Suppose $K$ is a knot in an integral homology sphere $Y$ with $\dim_\mathbb{C}I^\sharp(Y)=1$. If there is a positive integer $n$ such that $Y_{-n}(K)$ is an instanton L-space, then $\khii(-Y,K)$ consists of positive chains in the sense of Definition \ref{defn: positive chains}.
\ethm

\bpf
By Theorem \ref{thm: KHI for L-space knots}, for any integer $i$, we have$$\dim_\mathbb{C}\khii(-Y,K,i)\le 1.$$Then we have integers
$$n_1>n_{2}>\dots>n_k$$ such that $$\dim_\mathbb{C}\khii(-Y,K,i)=\begin{cases}1&\text{if }i=n_j\text{ for }j\in[0,k];\\0 &\text{else}. \end{cases}$$
Suppose $x_i$ is the generator of $\khii(-Y,K,n_{2i-1})$ and $y_i$ is the generator of $\khii(-Y,K,n_{2i})$. We verify that these $x_i$ and $y_i$ satisfy the positive chain condition, i.e.\, for any integer $i$, we have
\begin{equation}\label{eq: positive chain condition}
    d_-(y_i)\doteq x_{i+1},d_+(y_i)\doteq x_i,\aand d_-(x_i)=d_+(x_i)=0,
\end{equation}where $\doteq$ means the equation holds up to multiplication by a unit. We prove this condition by induction. We only consider the condition about the differential $d_+$. The proof for $d_-$ is similar. The gradings in the following arguments mean the gradings associated to the Seifert surface $S$. Note that by the proof of Theorem \ref{thm: KHI for L-space knots}, we have $$T_{n,n_{2l}}=B_{n,n_{2l-1}+1}=0\text{ for any }l.$$Hence by Lemma \ref{lem: tb iso}, we have$$T_{n,i}\cong H(B_{\ge n_{2l}}^{+})= H(B_{<2l-1}^-).$$

First, suppose $i=1$. Since $x_1$ lives in the top grading of $\khii(-Y,K)$ and $d_+$ increases the $\mathbb{Z}$-grading, we must have $d_+(x_1)=0$. Since $H(B_{\ge n_2}^+)=0$ and there are only two generators $x_1$ and $y_1$ in $B_{\ge n_2}^+$, we must have $d_+(y_1)\doteq x_1$.

Then we assume the condition (\ref{eq: positive chain condition}) holds for $i\le l-1$ and prove it also holds for $i=l$. Since $$H(B^+_{\ge n_{2l}})=H(B^+_{\ge n_{2l-2}})=0,$$we know the quotient complex $B^+_{\ge n_{2l}}/B^+_{\ge n_{2l-2}}$ also has trivial homology. Since it is generated by $x_l$ and $y_l$, the coefficient of $x_l$ in the expression of $d_+(y_l)$ must be nontrivial. Hence, $y_l$ is not in the $(n_{2l-1}-n_{2l}+1)$-page of the spectral sequence associated to $d_+$. Since the other generators $x_1,\dots,x_{l-1},y_1,\dots,y_{l-1}$ have smaller gradings than $x_l$, we know by the construction of $d_+$ in Construction \ref{cons: recover filtered chain complex} that the coefficients of those generators in the expression of $d_+(y_l)$ are all zero. Hence, $d_+(y_l)\doteq x_l$. Since $d_+\circ d_+=0$, we have $d_+(x_l)=0$. Thus, we prove the condition holds for $i=l$.

\epf

\bpf[Proof of Theorem \ref{thm: positive chain}]
If $Y_{-n}(K)$ is an instanton L-space, then the proof is similar to that of Theorem \ref{thm: positive chain, integral homology sphere}. To prove a generalization of Theorem \ref{thm: KHI for L-space knots}, we need to remove the integral homology sphere assumption in Proposition \ref{prop: l space dim 1} and Proposition \ref{prop: large enough surgeries are all L-spaces}. The corresponding proofs follow from Remark \ref{rem: proof of l space dim 1} and the proof of \cite[Proposition 4]{boyer13lspace}. If $Y_{n}(K)$ is an instanton L-space, by Remark \ref{rem: chain mirror}, we can consider the mirror knot to obtain the result.
\epf

\subsection{A graded version of the K\"{u}nneth formula}\quad

In this subsection, we prove the following graded version of K\"{u}nneth formula for the connected sum of two knots.

\bprop\label{prop: kunneth formula for connected sum}
Suppose $Y_1$ and $Y_2$ are two irreducible closed $3$-manifolds and $K_1\subset Y_1$, $K_2\subset Y_2$ are two rationally null-homologous knots such that $Y_1(K_1)$ and $Y_2(K_2)$ are both irreducible. Suppose $$(Y^\p,K^\p)=(Y_1\sharp Y_2,K_1\sharp K_2)$$is the connected sum of two knots. Then there is a minimal genus Seifert surface $S$ of $K^\p$ with the following properties.
\benu
\item There is a $2$-sphere $\Sigma\subset Y^\p$ intersecting the knot $K^\p$ in two points and intersecting $S$ in arcs. 

\item If we cut $S$ along $S\cap \Sigma^2$, then $S$ decomposes into two surfaces $S_1\subset Y_1$ and $S_2\subset Y_2$ such that $S_i$ is a union of some copies of Seifert surfaces of $K_i$ for $i=1,2$.
\item There is an isomorphism
\begin{equation}\label{eq: iso with grading}
\khii(Y^\p, K^\p,S,k)\cong\bigoplus_{i+j=k}\khii(Y_1,K_1,S_1,i)\otimes \khii(Y_2,K_2,S_2,j).
\end{equation}
\eenu
\eprop
\brem
In the first arXiv version of this paper, we assume $Y_1$ and $Y_2$ are rational homology spheres. Indeed in the proof, we only need the fact that $K_1$ and $K_2$ are rationally null-homologous knots. So we modify the assumption in later versions.
\erem
\brem
If the surface $S$ is the union of $n$ parallel copies of another surface $S^\p$, then $$khii(Y,K,S,i)=\khii(Y,K,S^\p,ni)$$by definition of the gradings.
\erem
\bpf[Proof of Proposition \ref{prop: kunneth formula for connected sum}]
Let $S$ be a minimal genus Seifert surface of $K^\p$ and let $\Sigma\subset Y^\p$ be a $2$-sphere such that $\Sigma$ intersects $K^\p$ in two points. We can choose $\Sigma$ such that
$$\Sigma\cap\partial Y^\p(K^\p)=\mu_1\cup\mu_2,$$
where $\mu_1$ and $\mu_2$ are two meridians of $K^\p$. Write 
$$A=\Sigma\cap Y^\p(K^\p).$$
From now on, we also regard $S$ as a surface inside the knot complement $Y^\p(K^\p)$. We can isotope $S$ so that $S$ intersects $A$ transversely and $S$ has minimal intersections with both $\mu_1$ and $\mu_2$. Now we argue that we can further isotope $S$ so that $S$ intersects $A$ in arcs. Suppose 
$$S\cap A=\al_1\cup\cdots\cup \al_n\cup\be_1\cup\cdots\cup\be_m,$$
where $\al_i$ are arcs and $\be_j$ are closed curves. Observe that each component of $A\backslash(\al_1\cup\cdots\cup \al_n)$ is a disk. Then using the arguments in the proof of \cite[Chapter 5, Theorem A14]{rolfsen2003knot}, we could further assume that $m=0$, i.e.\, $S$ intersects $A$ in arcs. When we cut the knot complement $Y^\p(K^\p)$ along $A$, we obtain the disjoint union of the knot complements $Y_1(K_1)$ and $Y_2(K_2)$, and the surface $S$ decomposes into $S_1\subset Y_1(K_1)$ and $S_2\subset Y_2$. Note that $S_1$ and $S_2$ must be the union of (possibly more than one) copies of Seifert surfaces of the corresponding knots. Then we prove the isomorphism (\ref{eq: iso with grading}).

First, we prove
\begin{equation}\label{eq: iso without grading}
	\khii(Y^\p, K^\p)\cong \khii(Y_1,K_1)\otimes \khii(Y_2,K_2).
\end{equation}
To do so, we pick a meridian $\mu_i'$ of $K_i$ for $i=1,2$ and pick suitable orientations such that $(Y^\p(K^\p),\mu_1'\cup\mu_2')$ is a balanced sutured manifold. Then we can decompose it along the annulus $A$:
$$(Y^\p(K^\p),\mu_1'\cup\mu_2')\leadsto(Y_1(K_1),\mu_1\cup\mu_1')\sqcup (Y_2(K_2),\mu_2\cup\mu_2').$$
From \cite[Proposition 6.7]{kronheimer2010knots}, this annular decomposition leads to the isomorphism (\ref{eq: iso without grading}). To study the grading behavior of this isomorphism, we sketch the construction of the isomorphism as follows. Pick a connected oriented compact surface $T$ such that
$$\partial T=-\mu_1\cup-\mu_2.$$
Pick an annulus $T'$ such that
$$\partial T^\p=-\mu_1^\p\cup -\mu_2^\p.$$
One could think of $T'$ as a copy of the annulus $A$.

In \cite[Section 7]{kronheimer2010knots}, Kronheimer and Mrowka constructed closures of $$(Y_1(K_1),\mu_1\cup\mu_1')\sqcup (Y_2(K_2),\mu_2\cup\mu_2')$$ as follows. First, glue $[-1,1]\times(T\cup T')$ to $Y_1(K_1)\sqcup Y_2(K_2)$ using the boundary identifications as above to obtain a pre-closure
\begin{equation}\label{eq: preclosure}
\widetilde{M}=(Y_1(K_1)\sqcup Y_2(K_2))\cup [-1,1]\times(T\cup T').
\end{equation}
The boundary of $\widetilde{M}$ has two components
$$\partial\widetilde{M}=R_+\cup R_-,$$
where
$$R_{\pm}=R_{\pm}(\mu_1\cup\mu_1')\cup R_{\pm}(\mu_2\cup\mu_2')\cup\{\pm1\}\times(T\cup T').$$
Second, choose an orientation preserving diffeomorphism
$$h:R_+\ra R_-$$
and use $h$ to close up $\widetilde{M}$ and obtain a closed $3$-manifold $Y$ with a distinguishing surface $R$. The pair $(Y,R)$ is a closure of $(Y_1(K_1),\mu_1\cup\mu_1')\sqcup (Y_2(K_2),\mu_2\cup\mu_2')$.

\brem
In \cite[Section 7]{kronheimer2010knots}, we also need to choose a simple closed curve in $Y$, either transversely intersecting $R$ at one point or is non-separating on $R$, to achieve the irreducibility condition for related instanton moduli spaces. In the current proof, the choices of simple closed curves are straightforward, so we omit them from the discussion.
\erem

Note that gluing $[-1,1]\times T_1$ to $(Y_1(K_1),\mu_1\cup\mu_1')\sqcup (Y_2(K_2),\mu_2\cup\mu_2')$ is the inverse operation of decomposing $(Y^\p(K^\p),\mu_1'\cup\mu_2')$ along the annulus $A$. As a result, $(Y,R)$ is clearly a closure of $(Y^\p(K^\p),\mu_1'\cup\mu_2')$ as well. The identification of the closures induces the isomorphism in (\ref{eq: iso without grading}). More precisely, we can pick the surface $T$ with large enough genus and pick a simple closed curve $\theta\subset T$ such that $\theta$ separates $T$ into two parts, both of large enough genus, and with $-\mu_1'$ and $-\mu_2'$ sitting in different parts. We also pick a core $\theta'$ of the annulus $T'$. When choosing the gluing diffeomorphism $h:R_+\ra R_-$, we can choose one such that
\begin{equation}\label{eq: choice of gluing diffeomorphism}
h(\{1\}\times\theta)=\{-1\}\times \theta,~{\rm and~}h(\{1\}\times\theta')=\{-1\}\times \theta'.
\end{equation}
Hence, inside $Y$, there are two tori $S^1\times\theta$ and $S^1\times \theta'$. 
If we cut $Y$ open along these two tori and reglue, then we obtain two connected $3$-manifolds $(Y_1,R_1)$ and $(Y_2,R_2)$, which are closures of $(Y_1(K_1),\mu_1\cup\mu_1')$ and $(Y_2(K_2),\mu_2\cup\mu_2')$, respectively. The Floer's excision theorem in \cite[Section 7.3]{kronheimer2010knots} then provides the desired isomorphism.

To study the gradings, recall that
$$S\cap A=\al_1\cup\cdots\cup \al_n$$
where $\al_i$ are arcs connecting $\mu_1$ to $\mu_2$ on $A$. We can also regard those arcs as on the annulus $T'$. Assume that $\partial S$ intersects each of $\mu_1'$ and $\mu_2'$ in $n$ points as well. Note that we have assumed that $T$ has a large enough genus. Then there are arcs $\delta_1$,..., $\delta_n$ such that the following holds. Recall we have chosen $\theta\subset T$ in previous discussions above.
\benu
	\item We have
	$\partial (\delta_1\cup\cdots\cup\delta_n)=S\cap(\mu_1'\cup\mu_2').$
	\item For $i=1,..,n$, the arc $\delta_i$ intersects $\theta_1$ transversely once.
	\item The surface $S\backslash(\delta_1\cup\cdots\cup\delta_n\cup\theta_1)$ also has two components.
	\item Let $\widetilde{S}=S\cup [-1,1]\times(\al_1\cup\cdots\cup\al_n)$ be a properly embedded surface inside the pre-closure $\widetilde{M}$ as in (\ref{eq: preclosure}), then we can choose a gluing diffeomorphism $h:R_+\ra R_-$ satisfying the condition (\ref{eq: choice of gluing diffeomorphism}) and the following extra condition
	$$h(\partial \widetilde{S}\cap R_+)=\partial \widetilde{S}\cap R_-.$$
\eenu
Hence, the surface $S$ extends to a closed surface $\bar{S}\subset Y$ that induces the desired $\intg$-grading on $\khii(Y^\p, K^\p)$. When we cut $Y$ open along $S^1\times\theta$ and $S^1\times \theta'$ and reglue, the surface $\bar{S}$ is also cut and reglued to form two closed surfaces $\bar{S}_1\subset Y_1$ and $\bar{S}_2\subset Y_2$. They are the extensions of the Seifert surface $S_1$ of $K_1$ and the Seifert surface $S_2$ of $K_2$ in the corresponding closures. Hence, the Floer's excision theorem in \cite[Section 7.3]{kronheimer2010knots} provides the desired isomorphism (\ref{eq: iso with grading}).
\epf

\subsection{Proofs of theorems in the introduction}\quad

In this subsection, we prove Theorem \ref{thm: main 2}, Theorem \ref{main: L space knot}, and Theorem \ref{thm: floer simple khi}.

\bpf[Proof of Theorem \ref{thm: main 2}]
By Remark \ref{rem: su2 abundant}, we may assume $S_n^3(K)$ is an instanton L-space for some $n\in\mathbb{N}_+$. Then by Theorem \ref{thm: positive chain}, the space $\khii(S^3,K)$ consists of coherent chains. Then arguments about $KHI(S^3,K,S,i)$ follow from Definition \ref{defn: positive chains} and Proposition \ref{prop: z2 shift}. 

To prove $K$ is a prime knot, we can apply the proof of \cite[Corollary 1.4]{baldwin18fibred} to $KHI$, replacing \cite[Theorem 1.1]{baldwin18fibred} by \cite[Theorem 1.7]{baldwin2018khovanov}. Note that we need the graded version of K\"{u}nneth formula for $KHI$ in Proposition \ref{prop: kunneth formula for connected sum}.
\epf

\bpf[Proof of Theorem \ref{thm: floer simple khi}]
By (\ref{eq: mirror iso}), a knot $K\subset Y$ is instanton Floer simple if and only if the mirror knot $(-Y,K)$ is instanton Floer simple. Note that by spectral sequences in Theorem \ref{thm: spectral sequence}, we always have $$\dim_\mathbb{C}\khii(-Y,K,[s])\ge \dim_\mathbb{C}I^\sharp(-Y,[s])\ge 1.$$By Remark \ref{rem: chain mirror}, we know that $(-Y,K)$ is instanton Floer simple if and only if the space $\khii(-Y,K)$ consists of both positive chains and negative chains. 

By Theorem \ref{thm: large surgery formula, bent} and Theorem \ref{thm: large surgery formula, dual bent}, if $K$ is instanton Floer simple, then for any large integer $n$, the manifolds $Y_{n}(K)$ and $Y_{-n}(K)$ are both instanton L-spaces. By the similar argument in the proof of \cite[Proposition 4]{boyer13lspace}, the manifold $Y_r(K)$ is an instanton L-space for any $|r|\ge n$.

Conversely, if for any $r$ with $|r|$ sufficiently large, the manifold $Y_r(K)$ is an instanton L-space, then for any large integer $n$, the manifolds $Y_{n}(K)$ and $Y_{-n}(K)$ are both instanton L-spaces. By Proposition \ref{thm: positive chain}, the space $\khii(-Y,K)$ consists of both positive chains and negative chains. Hence, $K$ is an instanton Floer simple knot.
\epf

Finally, we prove Theorem \ref{main: L space knot}. Suppose $K\subset Y$ is a knot with $H_1(Y(K))\cong \mathbb{Z}$, and suppose $\hat{\mu}$ is the meridian of $K$ with $q =S\cdot \hat{\mu}$, where $S$ is the Seifert surface of $K$. We choose a basis $(\hat{\mu},\hat{\lambda})$ of $H_1(\partial Y(K))$ as in Definition \ref{defn_4: Ga_n-hat sutures} and identify the slope with rational numbers. Then we have the following lemma.

\blem[{\cite[Lemma 2.7]{Rasmussen2017}}]\label{lem: splicing}
Consider the setting as above. If $r=u/v$, the manifold $Y_r(K)$ is obtained from $Y^\p=Y\sharp L(v,-u)$ by some integral surgery on $K^\p=K\sharp K(v,-u,1)$, where $K(v,-u,1)$ is the unique knot in $L(v,-u)$ such that the complement is diffeomorphic to $S^1\times D^2$. Moreover, we have $$H_1(Y^\p(K^\p))\cong H_1(Y(K))\cong H_1(S^1\times D^2)/(\hat{\mu},\mu^\p),$$where $\mu^\p$ is the meridian of $K(v,-u,1)$. Hence, $H_1(Y^\p(K^\p))\cong \mathbb{Z}$ if and only if $\gcd(q,v)=1$.
\elem

\bpf[Proof of Theorem \ref{main: L space knot}]
By \cite[Lemma 3.2]{Rasmussen2017}, for a Heegaard Floer L-space knot $K\subset Y$, the space $\widehat{HFK}(Y,K)$ satisfies similar coherent chain condition as in Definition \ref{defn: positive chains}. Consider the $\mathbb{Z}$-grading on $\widehat{HFK}(Y,K)$ induced by pairing the first Chern class of the spin$^c$ structure with $S$. Since $H_1(Y(K))\cong \mathbb{Z}$, the $\mathbb{Z}$-grading encodes all information in the spin$^c$ decomposition and the coherent chain condition implies$$\dim_\ft \widehat{HFK}(Y,K,S,i)\le 1.$$Hence the dimension is determined by the graded Euler characteristic.

If $v=1$ and $r\in\mathbb{Z}$, then by a similar discussion as above, Theorem \ref{thm: positive chain} implies that $\khii(Y,K,S,i)$ is determined by the graded Euler characteristic.  Hence the theorem follows from (\ref{eq: chi same}). 

If $v\neq 1$, then by Lemma \ref{lem: splicing}, we can apply the proof for $v=1$ to $$(Y^\p,K^\p)=(Y\sharp L(v,-u),K\sharp K(v,-u,1)).$$Note that simple knots are instanton Floer simple knots by \cite[Proposition 1.9]{LY2020}. Then the theorem follows from the graded K\"{u}nneth formula for $\khii$ (Proposition \ref{prop: kunneth formula for connected sum}) and $\widehat{HFK}$ (\cite[Section 5]{Ozsvath2011rational}). We do not need to consider the irreducible condition due to the convention in Subsection \ref{subsec: basic setups}. 
\epf

\section{Dehn surgeries along genus-one knots}\label{sec: Dehn surgery along genus-one knots}

In this section, we study the framed instanton Floer homology of Dehn surgeries along knots that satisfy the following conditions:
\benu
\item The genus of the knot is $1$, i.e.\, $g(K)=1$.
\item The instanton knot homology of the knot is determined by the Alexander polynomial, i.e.\,
$$\Delta_{K}(t)=a_1t+a_0+a_{-1}\aand \dim_{\mathbb{C}}\khii(S^3,K,i)=|a_i|\text{ for }i\in\mathbb{Z}.$$
\eenu
Such knots include all genus-one Khovanov-thin knots (in particular, genus-one quasi-alternating knots; see \cite[Corollary 1.6]{kronheimer2011khovanov}). In Table \ref{crosing}, we list all genus-one alternating knots with crossings $\le 12$; these are also all known examples of genus-one quasi-alternating knots. The data are from KnotInfo \cite{knotinfo}. Note that we normalize the Alexander polynomial by (\ref{eq: alex conditions}). The first knot for each crossing number in the table is a twisted knot. The reader can compare this table with examples in \cite{baldwin2020concordance}.


\begin{table}[htbp]
\caption{Genus-one alternating knots with crossings $\le 12$\label{crosing}}
\begin{tabular}{|c|c|c|c|c|c|}
\hline  
No.&Name&4-ball genus&Signature&Two-bridge notation&Alexander polynomial\\
\hline  
$1$&$3_1$&$1$&$-2$&$3/1$&$t-1+t^{-1}$\\
$2$&$4_1$&$1$&$0$&$5/2$&$-t+3-t^{-1}$\\
$3$&$5_2$&$1$&$-2$&$7/3$&$2t-3+2t^{-1}$\\
$4$&$6_1$&$0$&$0$&$9/7$&$-2t+5-2t^{-1}$\\
$5$&$7_2$&$1$&$-2$&$11/5$&$3t-5+3t^{-1}$\\
$6$&$7_4$&$1$&$-2$&$15/11$&$4t-7+4t^{-1}$\\
$7$&$8_1$&$1$&$0$&$13/11$&$-3t+7-3t^{-1}$\\
$8$&$8_3$&$1$&$0$&$17/4$&$-4t+9-4t^{-1}$\\
$9$&$9_2$&$1$&$-2$&$15/7$&$4t-7+4t^{-1}$\\
$10$&$9_5$&$1$&$-2$&$23/17$&$6t-11+6t^{-1}$\\
$11$&$9_{35}$&$1$&$-2$&&$7t-13+7t^{-1}$\\
$12$&$10_1$&$1$&$0$&$17/15$&$-4t+9-4t^{-1}$\\
$13$&$10_3$&$0$&$0$&$25/6$&$-6t+13-6t^{-1}$\\
$14$&$11a_{247}$&$1$&$-2$&$19/17$&$5t-9+5t^{-1}$\\
$15$&$11a_{343}$&$1$&$-2$&$31/27$&$8t-15+8t^{-1}$\\
$16$&$11a_{362}$&$1$&$-2$&&$10t-19+10t^{-1}$\\
$17$&$11a_{363}$&$1$&$-2$&$35/29$&$9t-17+9t^{-1}$\\
$18$&$12a_{803}$&$1$&$0$&$21/2$&$-5t+11-5t^{-1}$\\
$19$&$12a_{1166}$&$1$&$0$&$33/4$&$-8t+17-8t^{-1}$\\
$20$&$12a_{1287}$&$1$&$0$&$37/6$&$-9t+19-9t^{-1}$\\
\hline 
\end{tabular}

\end{table}

From the conditions in (\ref{eq: alex conditions}), there are two possibilities for the Alexander polynomial:
\benu
\item[(i)] $\Delta_K(t)=at-(2a-1)+at^{-1}$ for some $a\in\mathbb{N}_+$;
\item[(ii)] $\Delta_K(t)=-at+(2a+1)-at^{-1}$ for some $a\in\mathbb{N}_+$.
\eenu
We treat these two cases separately in the following two subsections.

\begin{conv}
For simplicity, we write $\khii(K)$ for $\khii(-S^3,K)$ and $\khii(K,i)$
for $\khii(-S^3,K,S,i)$, where $S$ is a Seifert surface of $K$. Recall that we write $\widebar{K}$ for the mirror knot of $K$. We will write $H(C)$ for the homology of a complex $C$ and write $f_*$ for the induced map between homologies.
\end{conv}
Recall the results from Section \ref{subsec: the bent complex}. In this case, we have $A_s=(\khii(K),d_s)$ for any $s$, and
\[
d_s(x)=\begin{cases}
d_+(x)&{\rm gr}(x)>s,\\
d_+(x)+d_-(x)&{\rm gr}(x)=s,\\
d_-(x)&{\rm gr}(x)<s.
\end{cases}\]
where ${\rm gr}(x)$ is the grading of $x\in \khii(K)$ associated to the Seifert surface. We can further decompose the differentials as follows:
$$d_{+}=\sum_{i<j}d_{j}^{i}\aand d_{-}=\sum_{i>j}d_{j}^{i},\text{ where }d_{j}^{i}:\khii(K,i)\ra \khii(K,j).$$Since $g(K)=1$, the $-3$-surgery is a large surgery in the sense of Theorem \ref{thm: large surgery formula, main}. Hence, we have
$$I^{\sharp}(-S^3_{-3}(K))\cong\bigoplus_{s=-1}^{1}H(A_s,d_s),$$where
$$H(A_1,d_1)\cong H(\khii(K),d_-)\cong I^\sharp(-S^3)\cong \mathbb{C},$$and$$H(A_{-1},d_{-1})\cong H(\khii(K),d_+)\cong I^\sharp(-S^3)\cong \mathbb{C}.$$
Hence we know that
\begin{equation}\label{eq: dim formula}
    \dim_{\mathbb{C}}I^{\sharp}(-S^3_{-3}(K))=2+\dim_{\mathbb{C}}H(A_0,d_0).
\end{equation}
Since $a_1=a_{-1}$, by the graded Euler characteristic of $KHI$ \cite{Lim2009,kronheimer2010instanton}, we know that the parities of $\khii(K,1)$ and $\khii(K,-1)$ are the same under the $\intg_2$ grading. By Proposition \ref{prop: z2 shift}, we know that there is no $d^1_{-1}$ or $d^{-1}_1$ differentials. Hence, we know that 
$$d_0=d^0_1+d^0_{-1}.$$

\subsection{The case of \texorpdfstring{$(2a+1)$}{(2a+1)}}\quad

In this case we know that
\begin{equation*}
\khii(K,i)\cong
\begin{cases}
	\mathbb{C}^{a}&i=\pm 1,\\
	\mathbb{C}^{2a+1}&i=0,\\
	0&\text{else.}
\end{cases}	
\end{equation*}
We have the following.

\blem\label{lem: 2a+1, d10 is injective}
The differential
$$d^1_0:\khii(K,1)\ra\khii(K,0)$$
is injective and the differential
$$d^0_{-1}:\khii(K,0)\ra\khii(K,-1)$$
is surjective.
\elem

\bpf
Since$$\dim_{\mathbb{C}}\ke(d^0_{-1})\geq \dim_{\mathbb{C}}\khii(K,0)-\dim_{\mathbb{C}}\khii(K,-1)=a+1$$
and
$$\dim_{\mathbb{C}}\im(d^1_0)\leq \dim_{\mathbb{C}} \khii(K,1)=a,$$
we know
$$1\leq \dim_{\mathbb{C}}(\ke(d^0_{-1})\slash\im(d^1_0))\leq \dim_{\mathbb{C}}H(A_1,d_1)=1.$$
We conclude that
$$\dim_{\mathbb{C}}\ke(d^0_{-1})=a+1$$
which means that $d^{0}_{-1}$ is surjective. Also, we must have 
$$\dim_{\mathbb{C}}\im(d^1_0)=a+1$$ which means that $d^1_0$ is injective.
\epf                                                              

\blem\label{lem: 2a+1, ker d01= ker d0-1}
We have $\ke(d^0_1)=\ke(d^0_{-1})\cong\mathbb{C}^{a+1}$.
\elem
\bpf
Applying the argument in Lemma \ref{lem: 2a+1, d10 is injective} to the bent complex $(A_{-1},d_{-1})$, we also conclude that $$\dim_\mathbb{C}\ke(d^0_{-1})={a+1}.$$Hence $\ke(d^0_1)\cong \ke(d_{-1}^0)$. Then we show they are indeed the same space. Suppose
$$x\in\ke(d_{-1}^0)\text{ such that } x\notin \ke(d^0_1).$$
Then we know that
$$d^1_0\circ d^0_1(x)\neq 0.$$
Since $x\in \ke(d_{-1}^0)$ and $\im(d^1_0)\subset \ke(d_{-1}^0)$, the map
$$(d^1_0\circ d^0_1)_*:H(\khii(K,0)\xra{d^0_{-1}}\khii(K,-1))\ra H(\khii(K,0)\xra{d^0_{-1}}\khii(K,-1))$$
is non-trivial. By Lemma \ref{lem: tb iso}, we can identify the map $(d^1_0\circ d^0_1)_*$ between bent complexes with the composition of bypass maps
$$\psi_{-,n}^{n+1}\circ\psi_{+,n+1}^{n+2}=0:B_{n+2,i}\ra B_{n,i}.$$
By Proposition \ref{prop: composition of bypasses is zero}, this map is zero, which is a contradiction. Hence, we conclude that
$$\ke(d_{-1}^0)\subset \ke(d_{1}^0).$$
Since they have the same dimension, they must be the same vector space.
\epf

\bprop\label{prop: case of 2a+1}
Suppose $K$ is a genus-one knot such that
$$\Delta_K(t)=at+(2a+1)+at^{-1}\text{ for }a\in\posi\aand \dim_{\mathbb{C}}\khii(K)=4a+1.$$Then for any $u,v\in\intg$ with $u\neq 0, v>0$ and $\gcd(u,v)=1$, we have
$$\dim_{\mathbb{C}}I^{\sharp}(S^3_{u/v}(K))=2av+|u|.$$
\eprop

\bpf
Applying Lemma \ref{lem: 2a+1, ker d01= ker d0-1} to $K$, we have
$$\dim_{\mathbb{C}}H(A_0,d_0)={2a+1}.$$
By (\ref{eq: dim formula}), we conclude that
$$\dim_{\mathbb{C}}I^{\sharp}(-S^3_{-3}(K))=2+\dim_{\mathbb{C}}H(A_0,d_0)=2a+3.$$
The same argument applies to the mirror $\widebar{K}$ of $K$, so we know that
$$\dim_{\mathbb{C}}I^{\sharp}(-S^3_{3}(K))=\dim_{\mathbb{C}}I^{\sharp}(-S^3_{-3}(\widebar{K}))=2a+3.$$
Then the proposition follows from \cite[Theorem 1.1]{baldwin2020concordance}.
\epf

\brem\label{rem: case of 2a+1}
Under the terminologies in \cite{baldwin2020concordance}, we know that $r_0(K)=2a$ and $\nu^{\sharp}(K)=0$ under the assumption of Proposition \ref{prop: case of 2a+1}. However, we do not know if $K$ is $V$-shaped or $W$-shaped in the sense of \cite[Definition 3.6]{baldwin2020concordance}. If $K$ is slice, then \cite[Theorem 3.7]{baldwin2020concordance} implies that it is $W$-shaped. If we knew the shape, then \cite[Theorem 1.1]{baldwin2020concordance} would also tell us $\dim_\mathbb{C}I^{\sharp}(S^3_0(K))$. 
\erem

\subsection{The case of \texorpdfstring{$(2a-1)$}{(2a-1)}}\quad

In this case we know that
\begin{equation*}
\khii(K,i)\cong
\begin{cases}
	\mathbb{C}^{a}&i=\pm 1,\\
	\mathbb{C}^{2a-1}&i=0,\\
	0&\text{else.}
\end{cases}	
\end{equation*}
Since
$$\ke(d^1_0)\subset H(A_1,d_1)\cong\mathbb{C}$$
hence we must have
$$\dim_{\mathbb{C}}\ke(d^1_0)\le 1.$$
Hence we have the following two subcases.
\benu
\item $\dim_{\mathbb{C}}\ke(d^1_0)=0$.
\item $\dim_{\mathbb{C}}\ke(d^1_0)=1$.
\eenu
\blem\label{lem: 2a-1, case 1, ker d01= ker d0-1}
We have $\ke(d^0_1)=\ke(d^0_{-1})\cong\mathbb{C}^{a}$ in Case (1).
\elem
\bpf
The condition $\dim_{\mathbb{C}}\ke(d^1_0)=0$ implies $d^1_0$ is injective and $\dim_{\mathbb{C}}\im(d^1_0)=a$. Since $\im(d^1_0)\subset \ke(d_{-1}^0)$, we know that $\dim_{\mathbb{C}}\ke(d_{-1}^0)\geq a$ and hence $\dim_{\mathbb{C}}\im(d^0_{-1})\leq a-1$. Since
$$\khii(K,-1)\slash(\im(d_{-1}^0))\subset H(A_1,d_1)\cong\mathbb{C},$$
we must have $\dim_{\mathbb{C}}\ke(d_{-1}^0)=a$.

Since $d^1_0$ is injective, by the proof of Lemma \ref{lem: 2a+1, ker d01= ker d0-1}, we know that $\ke(d^0_{-1})\subset \ke(d^0_1)$. Hence, we know that $\dim_{\mathbb{C}}\ke(d_{-1}^0)\geq a$ and hence $\dim_{\mathbb{C}}\im(d^0_{-1})\leq a-1$. Since
$$\khii(K,1)\slash(\im(d_{1}^0))\subset H(A_{-1},d_{-1})\cong\mathbb{C},$$
we must have $\dim_{\mathbb{C}}\ke(d_{1}^0)=a=\dim_{\mathbb{C}}\ke(d_{-1}^0)$ and hence $\ke(d^0_1)=\ke(d^0_{-1})$. 
\epf

To distinguish the bent complexes of $K$ and its mirror $\widebar{K}$, we write $A_s(K)$ and $A_s(\bar{K)}$, respectively. We write $\bar{d}^i_j$ for the component of differentials in $A_s(\widebar{K})$.

\blem\label{lem: 2a-1, case 2, ker d01= ker d0-1}
We have $\ke(d^0_1)=\ke(d^0_{-1})\cong\mathbb{C}^{a-1}$ in Case (2).
\elem
\bpf
Note that $\ke(d^1_0)\subset H(A_1,d_1)\cong\mathbb{C}.$
This means that
$$\ke(d^0_{-1})=\im(d^1_0)~{\rm and~}\im(d^0_{-1})=\khii(K,-1).$$

Consider the bent complex of the mirror knot. By Proposition \ref{prop: dual complex} and Corollary \ref{cor: involution}, we have a duality between $d_{j}^i$ and $\bar{d}_i^j$. In particular, we have 
$$\ke(\bar{d}^{-1}_0)\cong\cok(d^0_{-1})=0.$$
So we can apply Lemma \ref{lem: 2a-1, case 1, ker d01= ker d0-1} to $A_s(\widebar{K})$ and conclude that $\ke(\bar{d}^0_{1})=\ke(\bar{d}^0_{-1}).$
Using the duality again, we have $\im(d_0^{1})=\im(d^{-1}_{0})\cong \mathbb{C}^{a-1}.$ Hence $\ke(d^{-1}_0)\cong \mathbb{C}$. Since $$\ke(d^{-1}_0)\subset H(A_{-1},d_{-1}),$$
we conclude that 
$$\ke(d^0_{-1})=\im(d^{1}_0)=\im(d^{-1}_0)=\ke(d^0_1).$$
\epf
The following corollary is straightforward from the above discussion.
\bcor\label{cor: 2a-1, switch between case 1 and 2}
For a knot $K\subset S^3$, its bent complex $A_s(K)$ falls into Case (1) if and only if $A_S(\widebar{K})$ falls into Case (1).
\ecor

\bprop\label{prop: case of 2a-1}
Suppose $K$ is a genus-one knot such that
$$\Delta_K(t)=at+(2a-1)+at^{-1}\text{ for }a\in\posi\aand \dim_{\mathbb{C}}\khii(K)=4a-1.$$Then for any $u,v\in\intg$ with $u\neq 0, v>0$ and $\gcd(u,v)=1$, one and exactly one of the following two cases happens.
\benu
	\item[(a)]
$\dim_{\mathbb{C}}I^{\sharp}(S^3_{u/v}(K))=(2a-1)v+|u-v|.$
	\item[(b)]
$\dim_{\mathbb{C}}I^{\sharp}(S^3_{u/v}(K))=(2a-1)v+|u+v|.$
\eenu
\eprop
\bpf
When $A_s(K)$ falls into Case (1), from Lemma \ref{lem: 2a-1, case 1, ker d01= ker d0-1} we know that
$$\dim_{\mathbb{C}}H(A_0,d_0)={2a+1}.$$
Hence by (\ref{eq: dim formula}), we conclude that
$$\dim_{\mathbb{C}}I^{\sharp}(-S^3_{-3}(K))=2+\dim_{\mathbb{C}}H(A_0,d_0)=2a+3.$$
Furthermore, by Corollary \ref{cor: 2a-1, switch between case 1 and 2}, we know that $A_s(\widebar{K})$ falls into Case (1). By Lemma \ref{lem: 2a-1, case 2, ker d01= ker d0-1}, it follows that
$$\dim_{\mathbb{C}}I^{\sharp}(-S^3_{3}(K))=\dim_{\mathbb{C}}I^{\sharp}(-S^3_{-3}(\widebar{K}))=2a+1.$$ Then from \cite[Theorem 1.1]{baldwin2020concordance} we know that Case (a) holds. When $A_s(K)$ of $K$ falls into Case (1), by a similar proof, we know that Case (b) holds.
\epf

\brem
Note that $K$ satisfies Case (a) in Proposition \ref{prop: case of 2a-1}, if and only if $\widebar{K}$ satisfies Case (b) in Proposition \ref{prop: case of 2a-1}. The hypothesis of Proposition \ref{prop: case of 2a-1} only involves the genus, the Alexander polynomial, and the total dimension of the instanton knot homology of the knot, which are all impossible to be used to distinguish $K$ from its mirror.
\erem

\brem\label{rem: case of 2a-1}
The two cases of Proposition \ref{prop: case of 2a+1} correspond to the two cases where $\nu^{\sharp}(K)=1$ and $\nu^{\sharp}(K)=-1$, respectively. For genus-one alternating knots, from \citep[Corollary 1.10]{baldwin2020concordance} we know that 
$$\tau^{\sharp}(K)=-\frac{1}{2}\sigma(K),|\sigma(K)|\le 2,$$
$$2\tau^{\sharp}(K)-1\leq\nu^{\sharp}(K)\leq 2\tau^{\sharp}(K)+1,$$
and hence
$$-1\leq\nu^{\sharp}(K)\leq 1.$$
If we suppose further that the Alexander polynomial is of the form
$$\Delta_K(t)=at+(2a-1)+a^{-1},$$
then we have $\sigma(K)\neq0$ and hence $\tau^{\sharp}(K)=\nu^{\sharp}(K)=-\sigma(K)/2$. Thus, for genus-one alternating knots, which case of Proposition \ref{prop: case of 2a-1} happens depends on the signature of $K$.
\erem

\bpf[Proof of Theorem \ref{thm: alternating}]
The result in instanton theory is a combination of Proposition \ref{prop: case of 2a+1}, Proposition \ref{prop: case of 2a-1}, Remark \ref{rem: case of 2a+1}, and Remark \ref{rem: case of 2a-1}. The result in Heegaard Floer theory follows from \cite[Proposition 15]{Hanselman20surgery}.
\epf

\section{Examples of SU(2)-abundant knots}\label{sec:Examples of SU(2)-abundant knots}\quad

In this subsection, we provide many examples of $SU(2)$-abundant knots. 
\bprop\label{prop: classification of L-space knots}
Instanton L-space knots in $S^3$ are classified in the following cases
\benu
\item An alternating knot is an instanton L-space knot if and only if it is the torus knots $T(2,2n+1)$.
\item A Montesinos knot (in particular, a pretzel knot) is an instanton L-space knot if and only if it is the torus knot $T(2,2n+1)$, the pretzel knot $P(-2, 3, 2n + 1)$ for $n\in\mathbb{N}_+$, and their mirrors.
\item Knots that are closures of 3-braids are not instanton L-space knots except for the twisted torus knots $K(3,q;2,p)$ with $pq>0$ and their mirrors.
\eenu
\eprop
\bpf
Note that torus knots admit lens space surgeries \cite{Moser1971} and pretzel knots $P(-2, 3, 2n + 1)$ admit Seifert fibered L-space surgeries \cite{moore16pretzel}. Hence, they are instanton L-space knots.

Theorem \ref{thm: main SU2} provides many necessary conditions for instanton L-space knots. By \cite[Proposition 4.1]{Ozsvath2005}, if an alternating knot satisfies term (1) in Theorem \ref{thm: main SU2}, then it is the $T(2,2n+1)$ torus knot. Hence, hyperbolic alternating knots are not instanton L-space knots.

In \cite{baker18montesinos}, there is a classification of (Heegaard Floer) L-space knots for Montesinos knots. From \cite[Section 3.1]{baker18montesinos}, the proof of this classification only depends on term (1) in Theorem \ref{thm: main SU2}, the inequality (\ref{eq: det genus}), the fiberedness, and the strongly quasi-positive condition \cite[Theorem 1.5]{baldwin2019lspace}. Hence, the classification also works for instanton L-space knots.

In \cite{vafaee21braid}, it is shown that all closures of 3-braids except $K(3,q;2,p)$ do not satisfy term (1) and term (2) in Theorem \ref{thm: main SU2}. Hence, they are not instanton L-space knots.
\epf
\brem
For pretzel knots, there is another approach \cite{moore16pretzel} to classify L-space knots, which only depends on term (1) in Theorem \ref{thm: main SU2}, the inequality (\ref{eq: det genus}), the fiberedness, and the direct calculation on $\widehat{HFK}(S^3,P(3,-5,3,-2))$. However, it is hard to calculate $KHI(S^3,P(3,-5,3,-2))$ directly, so we use the approach in \cite{baker18montesinos}.
\erem
\brem
Note that $K=K(3,q;2,p)$ with $pq>0$ is a $(1,1)$-L-space knot from the proof of \cite[Theorem 3.1(a)]{vafaee14twisted}. By \cite[Corollary 1.7]{LY2020}, we know that $\dim_\mathbb{C}KHI(S^3,K)=\dim_{\ft}\widehat{HFK}(S^3,K)$. However, we do not know if $K$ is an instanton L-space knot because \cite[Theorem 3.1(a)]{vafaee14twisted} depends on the calculation of the chain complex $CFK^-(S^3,K)$ by working with a genus one doubly-pointed Heegaard diagram.
\erem
\bpf[Proof of Corollary \ref{cor: abundant knot su2}]
This follows directly from Proposition \ref{prop: classification of L-space knots} and Remark \ref{rem: su2 abundant}.
\epf

\brem
There is a family of twisted torus knots $K(p,q;2,m)$ with some conditions in \cite[Theorem 5]{morton06twist} whose Alexander polynomials do not satisfy term (1) in Theorem \ref{thm: main SU2}. Thus, those knots are also not instanton L-space knots and hence $SU(2)$-abundant. In general, the classification of L-space knots for twisted torus knots is still open; see \cite{vafaee14twisted,motegi16twisted,baker19twisted} for some special cases.
\erem

Then we consider satellite knots and cable knots. There are some useful theorems.

\bdefn[\cite{sivek20su2}]
A knot $K\subset S^3$ is called \textbf{$SU(2)$-averse} if there are infinitely many $r\in\mathbb{Q}\{0\}$ such that all representations $\pi_{1}(S_r^3(K))\to SU(2)$ have abelian images.
\edefn
\brem
If $b_1(Y) = 0$, then an $SU(2)$ representation of $Y$ has abelian image if and only if it has cyclic image.
\erem
\bthm[{\cite[Theorem 1.8]{sivek20su2}}]
Let $K\subset S^3$ be a nontrivial knot, and suppose that some satellite $P (K)$ with winding number $w$ is $SU(2)$-averse. Then we have the following.
\benu
\item If $P (U)$ is not the unknot $U$, then it is also $SU(2)$-averse.
\item If $w=0$, then $K$ is $SU(2)$-averse.
\eenu
\ethm

\bthm[{\cite[Theorem 10.6]{sivek20su2}}]
Let $K\subset S^3$ be a nontrivial knot, and let $p,q\in\mathbb{Z}$ satisfy $\gcd(p,q)=1$ and $q\ge 2$. If cable knot $K_{p,q}$ of $K$ is $SU(2)$-averse, then $K$ is also $SU(2)$-averse.
\ethm

\bthm[{\cite[Lemma 8.5]{baldwin2019lspace}}]
Let $K\subset S^3$ be a nontrivial knot, and let $p,q\in\mathbb{Z}$ satisfy $\gcd(p,q)=1,q\ge 2$, and $p/q>2g(K)-1$. Then the cable knot $K_{p,q}$ is a positive instanton L-space knot if and only if $K$ is an instanton L-space knot.
\ethm
\bdefn
A $K\subset S^3$ is called a \textbf{distinguished knot} if it is an alternating knot, a Montesinos knot, or a knot from a 3-braid except for the unknot, $T(2,2n+1)$, $P(-2, 3, 2n + 1)$ with $n\in\posi$, $K(3,q;2,p)$ with $pq>0$, and their mirrors.
\edefn
Note that distinguished knots are not instanton L-space knots and hence not $SU(2)$-averse by Remark \ref{rem: su2 abundant}. Then we have the following corollaries.
\bcor
Suppose $P(K)\subset S^3$ is a satellite knot with winding number $w\ge 0$ of the pattern $P\subset S^1\times D^2$. If one of the following holds, then $P(K)$ is not $SU(2)$-averse:
\benu
\item $P(U)$ is a distinguished knot;
\item $w\neq 0$ and $K$ is a distinguished knot.
\eenu
\ecor

\bcor
Let $K\subset S^3$ be a distinguished knot, and let $p,q\in\mathbb{Z}$ satisfy $\gcd(p,q)=1,q\ge 2$, and $p/q>2g(K)-1$. Then the cable knot $K_{p,q}$ is $SU(2)$-abundant.
\ecor

Finally, we strengthen a result in \cite[Theorem 1.8]{baldwin2019lspace}.
\begin{cor}
Suppose $K\subset S^3$ is a nontrivial knot and suppose $S^3_{3}(K)$ does not have irreducible $SU(2)$ representations. Then $K$ is a prime, fibered, strongly quasi-positive knot of genus two, and its instanton knot homology has the form
\begin{equation}\label{eq: special KHI}
\dim_\mathbb{C}KHI(S^3,K,S,i)=
\begin{cases}
	1&|i|\leq 2,\\
	0&{\rm else.}\\
\end{cases}	
\end{equation}
\end{cor}
\bpf
By Remark \ref{rem: su2 abundant}, we know that $K$ is an instanton L-space knot. Then Theorem \ref{thm: main 2} applies. By \cite[Theorem 1.8]{baldwin2019lspace} we know $K$ is fibered, so \cite[Theorem 1.11]{baldwin2018khovanov} applies and we obtain (\ref{eq: special KHI}).
\epf
\brem
In \cite{Liang2021lspace}, the first author and Liang proved that if $KHI(S^3,K)$ has the form (\ref{eq: special KHI}) for some knot $K\subset S^3$, then $K$ must be an instanton L-space knot. Then by \cite[Theorem 1.15]{baldwin2019lspace}, we know that $S^3_3(K)$ must be an instanton L-space. However, it is not enough to figure out whether $S^3_3(K)$ has irreducible $SU(2)$ representations.
\erem

\section{Further directions}\label{sec: future}\quad

In this section, we discuss some further directions of techniques introduced in this paper.

First, in Heegaard Floer homology, Ozsv\'{a}th and Szab\'{o} \cite{Ozsvath2008integral,Ozsvath2011rational} introduced a mapping cone formula. Roughly speaking, for a null-homologous knot $K$ in a closed 3-manifold $Y$, the homology $\widehat{HF}(Y_r(K))$ for any slope $r$ can be computed by the filtrations on $\widehat{CF}(Y)$ induced by $K$ and $-K$. The large surgery formula is the first step of their proof, which is recovered in instanton theory by Theorem \ref{thm: large surgery formula, main}. To prove an analog of the mapping cone formula in instanton theory, we need to further recover the following structures.
\bfa
Suppose $K$ is a null-homologous knot in a closed 3-manifold $Y$. For any integer $n$, suppose $W_n(K)$ is the cobordism from $Y$ to $Y_n(K)$ induced by attaching a 4-dimensional 2-handle and suppose $W^\p_n(K)$ is the cobordism from $Y_n(K)$ to $Y$ obtained from $-W_n(K)$ by turning around the two ends. We have the following structures in Heegaard Floer theory:
\benu
\item There is a spin$^c$ decomposition of the cobordism map:$$\widehat{HF}(W_n(K))=\sum_{\mathfrak{s}\in\spin(W_n(K))} \widehat{HF}(W_n(K),\mathfrak{s}):\widehat{HF}(Y)\to\widehat{HF}(Y_n(K)).$$Also, there is a spin$^c$ decomposition of $\widehat{HF}(W_n^\p(K))$.

\item For a large enough $n$, the spin$^c$ decomposition of $\widehat{HF}(W_n^\p(K))$ is compatible with some maps constructed by the filtrations on $\widehat{CF}(Y)$ from $K$ and $-K$.

\item For any integer $n$ and any positive integer $m$, there is a generalized surgery exact triangle\begin{equation*}
\xymatrix@R=6ex{
\widehat{HF}(Y_n(K))\ar[rr]&&\widehat{HF}(Y_{n+m}(K))\ar[dl]^{F}\\
&\bigoplus_{i=1}^m\widehat{HF}(Y)\ar[ul]&
}    
\end{equation*}where the map $F$ is related to the spin$^c$ decomposition of $\widehat{HF}(W_n^\p(K))$.
\eenu
\efa
Baldwin and Sivek constructed an analog of the term (1) in instanton theory when $b_1(W_n(K))=0$. The assumption of $b_1$ is due to the proof of some structure theorem for the cobordism map. If $b_1\ge 1$, then it is harder to prove the structure theorem. Also, in their construction, the closures used to define $I^\sharp(Y)$ and $I^\sharp(Y_n(K))$ are special (the connected sum with $T^3$). It is unknown how to extend the decomposition of the cobordism map to general closures of balanced sutured manifolds.

For term (2), we can still use the lifts of two spectral sequences to recover filtrations. However, without the decomposition of the cobordism map, it is impossible to write down a precise statement.

For term (3), we expect that the proof \cite{donaldson95triangle,scaduto2015instanton} of the usual exact triangle between $I^\sharp(Y),I^\sharp(Y_n(K))$, and $I^\sharp(Y_{n+1}(K))$ can be applied to the generalized triangle with some modifications.
\begin{conj}
Consider the manifolds defined above. For any integer $n$ and any positive integer $m$, there is an exact triangle\begin{equation*}
\xymatrix@R=6ex{
I^\sharp(Y_n(K))\ar[rr]&&I^\sharp(Y_{n+m}(K))\ar[dl]^{F^\p}\\
&\bigoplus_{i=1}^m I^\sharp(Y)\ar[ul]&
}    
\end{equation*}where the map $F$ is related to the cobordism $W_n^\p(K)$.
\end{conj}


Second, for any quasi-alternating knot $K\subset S^3$, Petkova \cite[Section 3]{Petkova2009thin} proved that the chain complex $CFK^-(S^3, K)$ is determined by $\Delta_K(t)$ and the signature $\sigma(K)$. The essential observation is that in this case, $CFK^-(S^3,K)$ is chain homotopic to $$(\widehat{HFK}(S^3,K)\otimes \ft[U],\partial_z+U\partial_w).$$where $\partial_z$ and $\partial_w$ shift the Alexander grading only by one. Then the result follows from the equation $\partial_z\circ\partial_w=\partial_w\circ \partial_z$ and algebraic lemmas. We can regard $d_+$ and $d_-$ on $\khii(-S^3,K)$ as analogs of $\partial_w$ and $\partial_z$ in instanton theory, respectively. If the following conjecture was proven, then we could apply algebraic lemmas in \cite[Section 3]{Petkova2009thin} to determine the differentials $d_+$ and $d_-$ by $\Delta_K(t)$ and $\sigma(K)$. By the large surgery formula, we could compute $I^\sharp(-S_{-n}^3(K))$ for $|n|\ge 2g(K)+1$. By results in \cite{baldwin2020concordance}, we might calculate $I^\sharp(-S_r^3(K))$ for any quasi-alternating knot, which generalizes Theorem \ref{thm: alternating}.
\begin{conj}
Suppose $K\subset S^3$ is a quasi-alternating knot and suppose the maps $d_+$ and $d_-$ are on $\khii(-S^3,K)$. Then the maps shift the grading associated to the Seifert surface by one, and the following equation holds$$d_+\circ d_-\doteq d_-\circ d_+,$$where $\doteq$ means it holds up to multiplication by a unit in $\mathbb{C}$. 
\end{conj}

\bibliographystyle{alpha}
\newcommand{\etalchar}[1]{$^{#1}$}

\end{document}